\documentclass{amsart}

\usepackage{latexsym,enumerate}
\usepackage{amsmath,amsthm,amsopn,amstext,amscd,amsfonts,amssymb}
\usepackage[ansinew]{inputenc}
\usepackage{verbatim}
\usepackage{psfig,labelfig} %% the command \centerline around \include graphics stretches the grid to textwidth
\usepackage{epsfig} % Para gr\'aficos
\usepackage{graphicx,psfrag} % Mejor \'este para gr\'aficos: es m\'as moderno
\usepackage{subfigure}
\usepackage{epstopdf} % Para complilar figuras eps con PDFLatex
\usepackage[mathscr]{eucal}
\usepackage[T3,T1]{fontenc}
\DeclareSymbolFont{tipa}{T3}{cmr}{m}{n}
\DeclareMathAccent{\invbreve}{\mathalpha}{tipa}{16}
\usepackage[normalem]{ulem}

\usepackage{color}
\definecolor{darkred}{rgb}{0.78,0.15,0.10}
\definecolor{darkblue}{rgb}{0.05, 0.15, 0.76}
\definecolor{darkgreen}{rgb}{0.12,0.40,0.10}
\definecolor{darkyellow}{rgb}{0.99,0.67,0.00}
\newcommand{\NN}{\mathbb{N}}

\newcommand{\RR}{\mathbb{R}}
\newcommand{\ZZ}{\mathbb{Z}}
\newcommand{\EE}{\mathbb{E}}

\renewcommand{\SS}{\mathbb{S}}
\newcommand{\HH}{\mathbb{H}}
\newcommand{\sL}{\mathscr{L}}
\newcommand{\sN}{\mathscr{N}}
\newcommand{\sR}{\mathscr{R}}
\newcommand{\sQ}{\mathscr{Q}}
\newcommand{\sC}{\mathscr{C}}
\newcommand{\sM}{\mathscr{M}}
\newcommand{\sK}{\mathscr{K}}
\newcommand{\sA}{\mathscr{A}}
\newcommand{\sE}{\mathscr{E}}
\newcommand{\sH}{\mathscr{H}}

\newtheorem{theorem}{Theorem}[section]
\newtheorem{lemma}[theorem]{Lemma}
\newtheorem{proposition}[theorem]{Proposition}
\newtheorem{corollary}[theorem]{Corollary}

\newtheorem{remark}[theorem]{Remark}
\theoremstyle{definition}
\newtheorem{upremark}[theorem]{Remark}
\newtheorem{definition}[theorem]{Definition}

\DeclareMathOperator{\cosech}{cosech}

\DeclareMathOperator{\arcsinh}{arcsinh}
\DeclareMathOperator{\arccosh}{arccosh}
\DeclareMathOperator{\arctanh}{arctanh}%Bu
\DeclareMathOperator{\dist}{dist}

\newcommand{\spb}[1]{\smallskip}
\newcommand{\mpb}[1]{\medskip}
\newcommand{\bpb}[1]{\bigskip}

\newcommand{\p}{\partial}

\renewcommand{\a}{\alpha}
\renewcommand{\b}{\beta}
\newcommand{\e}{\varepsilon}
\renewcommand{\d}{\delta}
\newcommand{\D}{\Delta}
\newcommand{\g}{\gamma}
\newcommand{\G}{\Gamma}
\renewcommand{\th}{\theta}
\renewcommand{\l}{\lambda}
\renewcommand{\L}{\Lambda}
\renewcommand{\O}{\Omega}
\newcommand{\s}{\sigma}
\newcommand{\abs}[1]{\vert #1 \vert}
\usepackage{accents}
\newcommand*{\bdot}[1]{%
 \accentset{\mbox{\footnotesize\bfseries .}}{#1}}

\newcommand{\Pmove}[2]{\!\!\overset{\ \scriptscriptstyle{#1\to #2}}{P}\!}

\numberwithin{equation}{section}

\begin{document}

\title[THE COLLAR THEOREM FOR COMPLETE SURFACES]{TOPONOGOV COMPARISON AND THE COLLAR THEOREM FOR COMPLETE SURFACES\\
WITH AN APPENDIX ON THE LEVEL SETS OF DISTANCE FUNCTIONS}

\author[Peter Buser]{Peter Buser}
\address{\'Ecole Polytechnique F\'ed\'erale de Lausanne, Switzerland}
\email{peter.buser@epfl.ch}

\author[Jos\'e M. Rodr{\'\i}guez]{Jos\'e M. Rodr{\'\i}guez}
\address{Universidad Carlos III de Madrid, ROR: https://ror.org/03ths8210, Departamento de Matem\'aticas, Avenida de la Universidad, 30 (edificio Sabatini), 28911 Legan\'es (Madrid), Spain, ORCID: 0000-0003-2851-7442}
\email{jomaro@math.uc3m.es}

\maketitle{}

\vspace{-12pt}
\begin{center}
 \today
\end{center}
\vspace{1pt}

\begin{abstract}
In the 1970s, the collar theorem was proven, establishing the existence of uniform tubular neighborhoods of simple closed geodesics on compact surfaces, whose widths depend only on the lengths of the geodesics and the lower bound of the curvature, but not on the surface.
In this paper, we improve this result by eliminating the compactness hypothesis.
To achieve this result, we needed to prove new Toponogov-type triangle comparison theorems.
We also add a new theorem to the literature on the rectifiabilty of the level sets of the distance function, with the corollary that on thin infinite cylinders with geodesic boundary \emph{all} sets of constant distance to the boundary are simple closed Lipschitz curves.
\end{abstract}



\section{Introduction}
\label{sec:Intro}

Let $S$ be a complete orientable Riemannian surface with curvature $K \ge -k^2$, where $k$ is a positive constant.
For any simple closed geodesic $\g$ in $S$ and $d>0$ we denote by
\begin{equation*}
C(d)
= \big\{ x\in S \,|\, \dist(x, \g) < d \, \big\}
\end{equation*}
the tubular neighborhood of radius $d$ of $\g$. If $d > 0$ is sufficiently small then $C(d)$ is doubly connected, i.e.\ topologically equivalent to $\SS^1 \times (-1,1)$. If $S$ has non-positive curvature and $C(d)$ is doubly connected then the same holds for all $C(d')$ with $0<d'< d$. In the general case this is not always so. For $\g$ a simple closed geodesic of length $L(\g)$ in $S$ we thus define $d_\g$ to be the supremum of all $d>0$ for which $C(d)$ and also all $C(d')$ with $0<d'<d$ are topologically equivalent to $\SS^1 \times (-1,1)$.
The cylinder $C(d_\g)$ of width $d_\g$ is called the \emph{collar about} $\g$.

If $S$ is compact with constant curvature $K = -k^2$, Randol proved in \cite{Ra} that
\begin{equation} \label{eq:a0}
d_\g
> \frac1k \, \arccosh \coth \Big( \frac{k}2\, L(\g) \Big)
= \frac1k \, \arcsinh \cosech \Big( \frac{k}2\, L(\g) \Big) .
\end{equation}
Chavel and Feldman proved in \cite{CF} that
\eqref{eq:a0} also holds if $S$ is compact and its curvature satisfies $-k^2 \le K \le -c^2$, with $c>0$.

After that, an alternative approach in \cite{Bu} (see also \cite{Bu2}) allowed the first author to show that \eqref{eq:a0} holds if $S$ is compact with $K \ge -k^2$.

While the previous arguments of Randol and Chavel and Feldman work for complete surfaces,
the argument in \cite{Bu} strongly uses compactness.

Our goal in this paper is to generalize the argument in \cite{Bu} for complete surfaces:

\begin{theorem} \label{t:a}
Let $S$ be a complete orientable Riemannian surface different from a closed torus and assume that it has curvature $K \ge -k^2$.
If $\g$ is a homotopically nontrivial simple closed geodesic in $S$, then
\begin{equation} \label{eq:a}
d_\g
\ge \frac1k \, \arccosh \coth \Big( \frac{k}2\, L(\g) \Big)
= \frac1k \, \arcsinh \cosech \Big( \frac{k}2\, L(\g) \Big) .
\end{equation}
\end{theorem}

\smallskip

As usual, we assume that every surface is connected.
Otherwise, Theorem \ref{t:a} can be applied to each connected component.

\smallskip

The arguments in the proof of \eqref{eq:a} also provide another inequality:

\begin{theorem} \label{t:b}
Let $S$ be a complete orientable Riemannian surface with curvature $K \ge -k^2$.
If $\g$ and $\mu$ are homotopically nontrivial simple closed geodesics in $S$,
and $\mu$ is not homotopic to and not intersecting $\g$, then the collars of $\g$ and $\mu$ of width $d_\g$ and $d_\mu$,
respectively, are disjoint and so,
\begin{equation} \label{eq:b}
\dist(\g,\mu)
\ge \frac1k \, \arccosh \coth \Big( \frac{k}2\, L(\g) \Big)
+ \frac1k \, \arccosh \coth \Big( \frac{k}2\, L(\mu) \Big) .
\end{equation}
\end{theorem}

The following result is a direct consequence of Theorem \ref{t:a}:

\begin{theorem} \label{t:c}
Let $S$ be a complete orientable Riemannian surface with curvature $K \ge -k^2$ different from a torus.
If $\g$ and $\eta$ are homotopically nontrivial simple closed geodesics in $S$,
and $\eta$ is not homotopic to but intersecting $\g$, then
\begin{equation} \label{eq:c}
\sinh \Big( \frac{k}2\, L(\g) \Big)
\sinh \Big( \frac{k}2\, L(\eta) \Big)
\ge 1 .
\end{equation}
\end{theorem}

The last section in \cite{Bu} contains examples showing that inequalities
\eqref{eq:a}, \eqref{eq:b} and \eqref{eq:c}
are sharp.

\smallskip

Although it is well known that the conclusion of Theorem \ref{t:a} does not hold for
homotopically trivial simple closed geodesics (e.g.\ if the geodesic bounds a half sphere with small radius), we still have a partial result for these geodesics:

Let $S$ be a complete orientable Riemannian surface which is not homeomorphic to a sphere.
If $\g$ is a homotopically trivial simple closed geodesic in $S$ and $D_\g$ is the topological open disk in $S$ with $\p D_\g = \g$, let us define
$d_\g^*$ as the largest $d$ for which the set
\begin{equation*}
C^*(d)
= \big\{ x\in S \,|\, \dist(x, \g) < d \, \big\} \setminus D_\g
= C(d) \setminus D_\g
\end{equation*}
and also all $C^*(d')$ with $0<d'<d$ are topologically equivalent to $\SS^1 \times [0,1)$.
The cylinder $C^*(d_\g^*)$ of width $d_\g^*$ is called the \emph{half-collar about} $\g$.

\begin{theorem} \label{t:d}
Let $S$ be a complete orientable Riemannian surface different from a sphere, with curvature $K \ge -k^2$.
If $\g$ is a homotopically trivial simple closed geodesic in $S$, then
\begin{equation*}
d_\g^*
\ge \frac1k \, \arccosh \coth \Big( \frac{k}2\, L(\g) \Big)
= \frac1k \, \arcsinh \cosech \Big( \frac{k}2\, L(\g) \Big) .
\end{equation*}
\end{theorem}

Some corollaries of these theorems follow in Section~\ref{sec:Improvements} and a particular version is proven in Section~\ref{sec:Thin} for the case where a simple closed geodesic bounds a thin infinitely long cylinder (Theorem \ref{t:cactus1}).

For a generalization of the Collar Theorem to higher dimensions and to Hitchin representations, see \cite{Ba} and \cite{LZ}, respectively. We also would like to point out the recent results on collars in complex hyperbolic manifolds by Basmajian and Kim \cite{BK}.

As a basic tool from Riemannian geometry we shall use Toponogov's triangle comparison theorem \cite{To}, \cite{K}, \cite{M}. For compact surfaces we can apply it in its standard form (see \cite{Bu}). For the non compact case, however,  we need some extensions. In particular, we shall need additional angle and side comparisons that do not hold in general, but are valid for right angled geodesic triangles (Propositions \ref{p:CompRightTriangles}, \ref{p:CompRightTriangle2}). As we could not find these extensions in the literature we insert an account on them in Section ~\ref{sec:Additional}.

We also insert an appendix on distance sets. It is well known (see Remark \ref{r:biblio} for some literature) that for a compact subset on a Riemannian surface almost all sets that have constant distance to this subset are rectifiable curves. In the special case of Section ~\ref{sec:Thin}, we shall use that the sets of constant distance $r$ to the boundary are rectifiable (in fact, Lipschitzian) for \emph{all} $r>0$. We shall base this on a more general theorem (Theorem~\ref{t:rectA} and its extensions Theorems~\ref{t:rectB},~\ref{t:rectC})  that we could not find in the literature and for which we shall give an elementary proof.

Theorems \ref{t:a}-\ref{t:d} are proved in Sections \ref{sec:Proof} and \ref{Proofd}.

\section{Background and technical results}
\label{sec:Background}

Since the products $k d_\g$, $\frac{k}2\, L(\g)$, etc. occurring in Theorems \ref{t:a}, \ref{t:b}, \ref{t:c}, \ref{t:d} are scaling invariant we shall, in most parts of what follows, replace the curvature bound $K \ge -k^2$ by $K \geq -1$.

\medskip

We begin with some general considerations concerning non-compact Riemannian surfaces.
As usual, we say that an end $E$ in a surface $S$ is \emph{collared} if $E$ has a neighborhood homeomorphic to $\SS^1 \times (0,\infty)$.
We say that a simple closed curve $\g$ \emph{bounds} a collared end $E$ in $S$ if one of the arc components $E'$ of $S \setminus \g$ is a neighborhood of $E$ homeomorphic to $\SS^1 \times (0,\infty)$. By abuse of language we shall also say that $E'$ and its closure $\overline{E'} = \gamma \cup E'$ are \emph{collared ends} of $S$.

\begin{definition}\label{d:minimizer1}
Given a Riemannian surface $S$
and a closed curve $\a$ in $S$, we define the
\emph{length of the free homotopy class} $[\a]$ as
\begin{equation*}
L([\a]) := \inf \big\{ L(\s) : \, \s \in [\a] \big\} \,.
\end{equation*}
The curve $\a$ is \emph{minimizing} if $L(\a) = L([\a])$.
\end{definition}

In Definition \ref{d:minimizer2} we shall also consider minimizing geodesics in other equivalence classes of curves on $S$.

\smallskip

Collars and the quest of their sizes arise in connection with the decomposition of complete Riemannian surfaces along geodesics into building blocks: cylinders, Y-pieces and half planes. These blocks will also be used as a tool in the proofs of the collar theorems and are defined as follows (see \cite{PRT}). The terminology is in analogy with the one used for surfaces of constant negative curvature.

\begin{definition}
Let $S$ be a complete bordered or unbordered Riemannian surface.

A \emph{halfplane} in $S$ is a complete bordered subsurface of $S$ that is simply connected
and whose border is a complete non closed simple geodesic.

A \emph{generalized funnel} is a bordered subsurface of $S$ whose border is a minimizing (with respect to $S$) simple closed geodesic and whose interior is a neighborhood of a collared end of $S$.

A \emph{generalized cusp} is a collared end of $S$ whose fundamental
group is generated by a simple closed curve $\s$ and there is no
minimizing closed geodesic in $[\s]$.

%A \emph{Y-piece} is a compact bordered Riemannian surface which is topologically a sphere without three open disks, whose boundary consists of three pairwise disjoint simple closed geodesics.

A \emph{generalized Y-piece} is a bordered or unbordered complete Riemannian surface which is topologically a sphere without three open disks, and for which there exist integers $n,m\ge 0$ with $n+m=3$, so that the boundary consists of $n$ pairwise disjoint simple closed geodesics and the surface has $m$  generalized cusps.

\end{definition}

%\begin{figure}[h]
%    \centering
%    {\epsfig{figure=imgs/GnFunnel.eps,height=3.5cm}}\qquad\qquad %
%    \subfigure[Generalized funnels]{\epsfig{figure=imgs/GnFunnel2.eps,width=4cm,height=3.5cm}}\qquad\qquad %
%    {\epsfig{figure=imgs/Funnel.eps,height=3.5cm}}
%    %\caption{Generalized funnels}
%\end{figure}
%
%\begin{figure}[h]
%    \centering
%    {\epsfig{figure=imgs/GnPunct.eps,height=4.0cm}}\qquad\qquad %
%    \subfigure[Generalized cusps]{\epsfig{figure=imgs/GnPunct2.eps,width=3.5cm,height=4cm}}\qquad\qquad %
%    {\epsfig{figure=imgs/Puncture.eps,height=3.5cm}}%
%    %\caption{Generalized cusps}
%\end{figure}

%\begin{figure}[h!]
%\centering
%\includegraphics[width=0.675\linewidth]{imgs/a_ThreeGnFunnels}%[width=0.7\linewidth]%0.7 perfekt
%\vspace{-6pt}
%\caption{Generalized funnels}
%\label{fig:CompRightTriangles}
%\end{figure}
%
%\begin{figure}[h!]
%\centering
%\rule{16pt}{0pt}\includegraphics{a_ThreeGnCusps}%[width=0.66\linewidth]%0.66 perfekt
%\vspace{-6pt}
%\caption{Generalized cusps}
%\label{fig:CompRightTriangles}
%\end{figure}

\begin{upremark}\label{r:remDef}
We comment the preceding definitions and introduce some additional terminology.

(1) In contrast to the funnels which are defined here as subsurfaces, the cusps are defined only as ends, and there is no geometrically privileged cusp neighborhood (i.e.\ collared neighborhood of the end), in general. It may, however, occur that the homotopy class $[\s]$ in the above definition contains \emph{non} minimizing geodesics. In that case, any such geodesic $\d$ cuts away a cusp neighborhood $S_{\d}$ from $S$ and we call its closure $S_{\d} \cup \d$ a \emph{cusp with geodesic boundary}.

(2) A (bordered or unbordered) surface is \emph{doubly connected} or a \emph{cylinder}
if its fundamental group is isomorphic to $\ZZ$ or, equivalently, if it is homeomorphic to $\SS^1 \times I$, where $I \subset \RR$ is an interval.
Every generalized funnel and every generalized cusp with geodesic boundary is a doubly connected surface homeomorphic to $\SS^1 \times [0, 1)$.

(3) If a generalized Y-piece is compact and thus has $n=3$ simple closed boundary geodesics we sometimes (but only under this hypothesis) omit the attribute \guillemotleft generalized\guillemotright\  and just call it a Y-piece. A clear description of these ``ordinary'' Y-pieces and their use is given in \cite[chapter X.3]{C}.

(4) A generalized Y-piece $Y$ is viewed as a \emph{subsurface of $S$} if its interior ($Y$ minus the boundary geodesics) is isometrically embedded in $S$. The embedding isometrically extends to each boundary geodesic of $Y$ individually  but it can occur that two boundary geodesics of $Y$ with the same length are mapped to the same geodesic of $S$.

(5) If a generalized Y-piece $Y$ is a subsurface of $S$ (in the sense of (4)) and has a cusp, say $E$, then by the completeness of the Riemannian metric on $Y$ and $S$ any cusp neighborhood $V$ of $E$ is at the same time a collared neighborhood of an end $E'$ of $S$ and $E'$ is independent of the choice of $V$. We may thus identify the end $E$ of $Y$ with the end $E'$ of $S$.

(6) In (5) it may occur that the end $E$ of $Y$ is a cusp with respect to the Riemannian surface $Y$ while the end $E'$ with respect to $S$ is not. This conceptual flaw plays no role in the following because cusps will only be looked at with respect to generalized Y-pieces. Furthermore, the discrepancy does not occur if all boundary geodesics of $Y$ that are not contained in the boundary of $S$ are minimizing geodesics in $S$.

\end{upremark}

\smallskip

In \cite[Theorem 4.3]{PRT} appears the following result generalizing the decomposition of compact surfaces into Y-pieces, and a previous result for surfaces with constant negative curvature in \cite{AR}.

\begin{theorem} \label{t:descomposicion}
Every complete unbordered orientable Riemannian surface
with non-abelian fundamental group
is the union (with pairwise disjoint interiors) of generalized Y-pieces,
generalized funnels and halfplanes.
\end{theorem}

Note that the only surfaces which are left out in Theorem \ref{t:descomposicion} are the simply and doubly connected surfaces and the tori.

These kinds of results about the decomposition of non-compact surfaces have shown to be useful in the study of escaping geodesics,
see \cite{FM} and \cite{MRT}.

\smallskip

In \cite[Theorem 4.12]{PRT} there is also a version for bordered surfaces of the following kind: we shall say that a Riemannian surface is \emph{bordered by simple closed geodesics} if the boundary is a non empty (finite or infinite) union of pairwise disjoint simple closed geodesics.

The statement in \cite[Theorem 4.12]{PRT} is in a refined form inasmuch as the generalized Y-pieces in the decomposition have minimizing boundary geodesics. This requires to use additional building blocks in the form of compact geodesically bordered cylinders, each of them having one boundary geodesic on the boundary of $S$ and the other on the boundary of one of the generalized Y-pieces. By pasting each such cylinder to the corresponding generalized Y-piece we obtain the theorem in  the following form.

\begin{theorem} \label{t:descomposicion2} %\label{c:descomposicion3}
Every complete orientable Riemannian surface bordered by simple closed geodesics which is neither simply nor doubly connected
is the union (with pairwise disjoint interiors) of generalized Y-pieces, generalized funnels and halfplanes.
Furthermore, the decomposition can be made such that all boundary geodesics of the generalized Y-pieces of the decomposition of $S$  that are not on the boundary of $S$ are minimizing.
\end{theorem}

A feature used for the above theorem is that distinct minimizing closed geodesics in the same homotopy class do not intersect each other. A similar fact holds also in other situations:

\begin{definition}\label{d:minimizer2}
Let $S$ be a complete Riemannian surface.
\begin{enumerate}
\item[(a)] For disjoint subsets $A, B \subset S$ a curve $c$ with initial point on $A$ and endpoint on $B$ is \emph{minimizing} if $c$ is a shortest curve with this property.
\item[(b)] Let $A \subset S$ be a subset of $S$ and $E$ an end of $S$. A curve $c$ with initial point on $A$ that escapes into $E$ is \emph{minimizing} if for any point $x$ on $c$ the arc on $c$ from $A$ to $x$ is minimizing.
\end{enumerate}
\end{definition}

\begin{lemma} \label{l:nocut}
Let $S$ be a complete Riemannian surface, $A \subset S$ a subset of $S$ and $F_1, F_2$ among the subsets and ends of $S$. Let $c_i$ for $i=1,2$ be a geodesic in $S$ with initial point on $A$ and endpoint on $F_i$ in the case where $F_i$ is a subset respectively, escaping into $F_i$ in the case where $F_i$ is an end. Furthermore, assume that, except possibly at its endpoints, $c_i$ does not meet the boundary of $S$. If $c_1$, $c_2$ are minimizing and distinct, then they intersect each other at most at their end points.
\end{lemma}

\begin{proof}Assume, for a contradiction, that $c_1$ intersects $c_2$ in some point $p$ different from the endpoints. Then $p$ separates $c_i$ into a segment $c_i'$ from its initial point on $A$ to $p$ and into the remaining part $\tilde{c}_i$. In the case where $F_i$ is a subset of $S$ we let $p_i$ be the endpoint of $\tilde{c}_i$ on $F_i$ and set $c_i'' = \tilde{c_i}$; in the case where $F_i$ is an end we let $p_i$ be the point on $\tilde{c}_i$ that has, say, distance 1 from $p$ and let $c_i''$ be the segment on $\tilde{c}_i$ from $p$ to $p_i$.

By the minimizing property of the two geodesics we have $L(c_1') + L(c_1'') \leq L(c_2') + L(c_1'')$ and  $L(c_2') + L(c_2'') \leq L(c_1') + L(c_2'')$. Hence, $L(c_1') = L(c_2')$ and so the curve $c_2' \cup c_1''$ is also minimizing from $A$ to $p_1$. But, the geodesic arcs $c_2'$, $c_1''$ meet at some angle different from $\pi$ at $p$ and $p$ is an interior point of $S$, which allows us to deform $c_2' \cup c_1''$ into a shorter curve from $A$ to $p_1$, a contradiction.
\end{proof}

\smallskip

We continue with some technical tools from \cite{Bu}. All surfaces are again assumed to be endowed with a complete Riemannian metric.

\begin{lemma}\cite[Lemma (2)]{Bu} \label{l:broken}
Let $G$ be a complete surface which is bordered by broken geodesics, the inner angles being less than $\pi$.
Then, any two disjoint compact subsets of $G$ can be joined by a geodesic arc of minimal length.
\end{lemma}

Note that Lemma \ref{l:broken} is proved in \cite{Bu} for compact surfaces, but the argument in the proof also works for complete surfaces. The Lemma shall frequently be used tacitly.

\begin{lemma}\cite[Lemma (3)]{Bu} \label{l:sinh}
Let $G$ be a compact right-angled geodesic pentagon with curvature $K \ge -1$,
and let $a_1,a_2$ be two sides of $G$ with a common vertex. Then,
$\sinh L(a_1) \sinh L(a_2) > 1$.
\end{lemma}
Lemma \ref{l:sinh}    has a continuation in Lemma \ref{l:sinh2} and Lemma \ref{l:openquadri}.

%In order to prove Theorem \ref{t:c4} below, we will need the following technical result in \cite[p.355]{Bu}:

\begin{lemma}\cite[Lemma (8)]{Bu} \label{l:B}
For every $0 \le |\d| \le x$ and $0 \le t \le y$, one has
\begin{equation*}
\sinh x \sinh y
\ge \min \big\{
\sinh (x+\d) \sinh (y-t),\,
\sinh (x-\d) \sinh (y+t)
\big\} .
\end{equation*}
\end{lemma}

%%%

The argument in the proof of Lemma \ref{l:B} in \cite[p.355]{Bu} also gives the following:

\begin{remark} \label{r:B}
The inequality in Lemma \ref{l:B} is strict for every $0 < |\d| \le x$ and $0 < t \le y$.
\end{remark}%%ok

\begin{corollary} \label{c:B}
For every $u_1,u_2,v_1,v_2\ge 0$, one has
\begin{equation*}
\sinh \frac{u_1+u_2}2 \,\sinh \frac{v_1+v_2}2
\ge \min \big\{
\sinh u_1 \sinh v_1,\,
\sinh u_2 \sinh v_2
\big\} .
\end{equation*}
This inequality is an equality if and only if one of the following three cases holds : $u_1 = u_2= 0$; $v_1=v_2=0$; $u_1=u_2$ and $v_1=v_2$.
\end{corollary}

\begin{proof}
The equality in the three cases is clear. Now assume that none of them holds.
By symmetry, we may assume that $v_2 > v_1$.
If we take
\begin{equation*}
x = \frac{u_1+u_2}2 \,,
\quad
y = \frac{v_1+v_2}2 \,,
\quad
\d = \frac{u_1-u_2}2 \,,
\quad
t = \frac{v_2-v_1}2 \,,
\end{equation*}
then Lemma \ref{l:B} gives the inequality.

If $u_1 \neq u_2$, then Remark \ref{r:B} implies that the inequality is strict.

If $u_1 = u_2$, then the inequality is
\begin{equation*}
\sinh u_1 \sinh \frac{v_1+v_2}2
\ge \sinh u_1 \min \big\{  \sinh v_1,\, \sinh v_2 \big\}
= \sinh u_1 \sinh v_1,
\end{equation*}
which is strict since $v_2 > v_1$.
\end{proof}

\begin{lemma} \label{l:sinh2}
Let $G$ be a compact geodesic pentagon with curvature $K \ge -1$, four angles less than or equal to $\pi/2$ and an angle $\a\in (0,\pi)$.
Let $a_1$ be the opposite side of the angle $\a$
and let $a_2$ be a side of $G$ with a common vertex with $a_1$. Then,
$\sinh L(a_1) \sinh L(a_2) > \sin \a$.
Furthermore, if $\a\in (0,\pi/2]$, then $\sinh L(a_1) \sinh L(a_2) > 1$.
\end{lemma}

%
%%% Figure 1 %%%
\begin{figure}[!b]
 \begin{center}
 \leavevmode
%\ShowGrid
\SetLabels
(.278*.995) $a_4''$\\
(.719*.975) $a_3''$\\
(.489*.88) $\overline{\alpha}$\\
(.350*.765) $\overline{a}_4$\\
(.705*.79) $\overline{a}_3$\\
(.210*.55) $\overline{a}_5$\\
(.86*.55) $\overline{a}_2$\\
(.536*.28) $\overline{a}_1$\\
(.732*.32) $\overline{G}$\\
(.786*.13) $G''$\\
(.017*.25) $a_5''$\\
(1.025*.24) $a_2''$\\
(.50*.03) $a_1''$\\
\endSetLabels
%\ShowGrid
\AffixLabels{%
\includegraphics{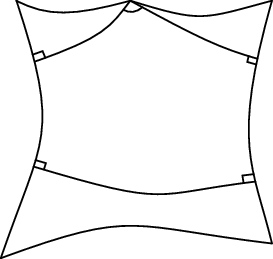}%[width=0.25\linewidth]
}
\end{center}
\caption{Construction of $\overline{G}$}
\label{fig:figura3}
\end{figure}
%%% end Figure 1 %%%
%

\begin{proof}
Denote by $a_3,a_4,a_5$ the three sides of $G$ such that $a_1,a_2,a_3,a_4,a_5$ are consecutive sides.
By Lemma \ref{l:broken}, there exist minimal geodesics $a_1',a_2',a_3',a_4',a_5'$ in $G$
such that $a_i'$ and $a_i$ have the same endpoints for each $1 \le i \le 5$.
Since these arcs have minimal length, they do not intersect each other, except for the common vertices (Lemma \ref{l:nocut}).
Hence, $a_1',a_2',a_3',a_4',a_5'$ are the sides of a geodesic pentagon $G' \subseteq G$
with four angles less than or equal to $\pi/2$ and an angle $0 < \a' \le \a < \pi$.
Thus, $L(a_i') \le L(a_i)$ for each $1 \le i \le 5$.
We triangulate $G'$ with two diagonals of minimal length in their homotopy classes.
Toponogov's comparison theorem \ref{t:top}, which we detail in the next section, provides a geodesic pentagon $G''$ in the hyperbolic plane (with curvature $-1$) with
sides $a_i''$ satisfying $L(a_i'') = L(a_i')$ for each $1 \le i \le 5$,
four angles less than or equal to $\pi/2$ and an angle $0 < \a'' \le \a' \le \a < \pi$.

Let us consider the vertex $A=a_3'' \cap a_4''$ in $G''$ with angle $\a''$.
Consider two minimizing geodesics $\overline{a}_3$, $\overline{a}_4$, from $A$ to $a_2''$, $a_5''$, respectively,
and a minimizing geodesic $\overline{a}_1$ joining $a_2''$ and $a_5''$.
By Lemma \ref{l:nocut} again, since these arcs have minimal length, they do not intersect each other.
Let $\overline{a}_2$ be the arc in $a_2''$ joining the endpoints of $\overline{a}_1$ and $\overline{a}_3$ in $a_2''$,
and $\overline{a}_5$ be the arc in $a_5''$ joining the endpoints of $\overline{a}_1$ and $\overline{a}_4$ in $a_5''$.
Therefore, there exists a geodesic pentagon $\overline{G} \subseteq G''$ with
sides $\overline{a}_i$ satisfying $L(\overline{a}_i) \le L(a_i'')$ for each $1 \le i \le 5$,
four angles equal to $\pi/2$ and an angle $0 < \overline{\a} \le \a'' \le \a' \le \a < \pi$
(see c \ref{fig:figura3}).

If $0 < \overline{\a} \le \pi/2$, then the argument in the proof of Lemma \ref{l:sinh} in \cite{Bu} gives
\begin{equation*}
\begin{aligned}
\sinh L(a_1) \sinh  L(a_2)
& \ge \sinh L(a_1') \sinh  L(a_2')
\ge \sinh L(a_1'') \sinh  L(a_2'')
\\
& \ge \sinh L(\overline{a}_1) \sinh  L(\overline{a}_2)
> 1 .
\end{aligned}
\end{equation*}

Assume now that $\pi/2 < \overline{\a} < \pi$.
Since $\overline{G}$ has curvature $-1$, we have, see e.g. \cite[p.87]{Fe},
\begin{equation*}
\frac{\sinh L(\overline{a}_1)}{\sin \overline{\a}}
= \frac{\cosh L(\overline{a}_4)}{\sinh L(\overline{a}_2)}\,,
\qquad
\sinh L(\overline{a}_1) \sinh L(\overline{a}_2) = \cosh L(\overline{a}_4) \sin \overline{\a} > \sin \overline{\a} ,
\end{equation*}
and we conclude
\begin{equation*}
%\begin{aligned}
\sinh L(a_1) \sinh  L(a_2)
\ge \sinh L(\overline{a}_1) \sinh  L(\overline{a}_2)
> \sin \overline{\a}
\ge \sin \a .
%\end{aligned}
\end{equation*}
\end{proof}

\section{Additional comparison arguments}
\label{sec:Additional}

In this section we apply Toponogov's comparison theorem to closed and open-ended Riemannian polygons. By the latter we mean the following.

\begin{definition}\label{d:Riempoly}
A \emph{Riemannian polygon} is a simply connected metrically complete bordered domain $P$ in an unbordered (not necessarily complete) Riemannian surface $R$, where the interior of $P$ is non-empty and the boundary of $P$ is a connected simple curve consisting of finitely many geodesic arcs and at most two rays, called the \emph{sides} of $P$.

A side $s$ of $P$ is called \emph{straight}, more precisely \emph{straight with respect to $P$}, if for any pair of points $p, q \in s$ the arc on $s$ from $p$ to $q$ is the shortest connection from $p$ to $q$ in $P$ (in $R$ shorter connections are allowed to exist).
\end{definition}

When $P$ is compact all sides are geodesic arcs and we say that $P$ is \emph{closed}. When $P$ is non-compact exactly two sides are geodesic rays and we say that $P$ is \emph{open-ended}; the two rays may or may not have the same ideal endpoint at infinity. The cases occurring in this paper are closed and open-ended Riemannian triangles, quadrangles and pentagons.

The Toponogov triangle comparison theorem is usually stated in the realm of $n$-dimensional unbordered complete Riemannian manifolds. In dimension 2, however, the Riemannian surface $R$ that contains the triangle does not have to be complete. For clarity we restate the theorem for a closed Riemannian triangle. For the proof we refer to \cite{K} or \cite{M}, where no use of the ambient space is made. To keep the statement short we assume a negative lower curvature bound.

\begin{theorem}[Toponogov] \label{t:top}
Let $T$ with sides $a$, $b$, $c$ and opposite angles $\a$, $\b$, $\g$ be a closed Riemannian triangle. Assume that all sides are straight and that the sectional curvature $K$ has the lower bound $K \geq -1$. Then the following hold.
\begin{enumerate}
\item[A)] There exists a geodesic triangle $\widetilde{T}$ in the hyperbolic plane with sides $\tilde{a}$, $\tilde{b}$, $\tilde{c}$ and opposite angles $\tilde{\a}$, $\tilde{\b}$, $\tilde{\g}$ satisfying
\begin{equation*}
L(\tilde{a}) = L(a), L(\tilde{b}) = L(b),  L(\tilde{c}) = L(c) \quad \text{and} \quad \tilde{\a} \leq \a, \tilde{\b} \leq \b, \tilde{\g} \leq \g.
\end{equation*}

\item[B)] There exists a triangle $\overline{T}$ in the hyperbolic plane with sides $\bar{a}$, $\bar{b}$, $\bar{c}$ and opposite angles $\bar{\a}$, $\bar{\b}$, $\bar{\g}$ satisfying
\begin{equation*}
L(\bar{a}) = L(a), L(\bar{b}) = L(b),  \bar{\g} =  \g  \quad \text{and} \quad  L(\bar{c}) \geq L(c).\end{equation*}
\end{enumerate}
\end{theorem}
The so called \emph{comparison triangles} $\widetilde{T}$ and $\overline{T}$ in the theorem are uniquely determined by the equalities up to isometry. In B) there are no general inequalities for the angles adjacent to side $c$. This is, however different in the case of right angles where we have the following particular properties.

\begin{proposition} \label{p:CompRightTriangles}
Let $D$ be a closed Riemannian triangle with curvature $K \geq -1$ and straight sides $a$, $b$, $c$. Assume that $D$ has a right angle between $a$ and $b$ and assume that the angles $\a$, $\b$ opposite to $a$ and $b$ satisfy $\a \leq \pi/2$, $\b \leq \pi/2$. Consider for comparison a possibly open-ended right-angled geodesic triangle $D'$ in the hyperbolic plane with sides $a'$, $b'$, $c'$, the right angle between $a'$ and $b'$ and angles $\a'$, $\b'$ opposite to $a'$ and $b'$ (if $a'$ or $b'$ is a geodesic ray the angle at its endpoint at infinity is by definition 0).

Then the following comparisons hold, where we abbreviate $L(a)=a$, $L(b)=b$, etc.:
\begin{equation*}
\renewcommand{\arraystretch}{1.3}
\arraycolsep=2.3pt
\begin{array}{lrccccccccccccccr}
(1)
&
\textit{If}
&
a'=a
&
\textit{and}
&
b'=b,
&
&
\textit{then}
&
&
(i)
&
\a' \leq \a
&
&
(ii)
&
\b' \leq \b
&
&
(iii)
&
c'\geq c,
&
\rule{120pt}{0pt}
\\
%%end line 1
(2)
&
\textit{if}
&
a'=a
&
\textit{and}
&
\b'=\b,
&
&
\textit{then}
&
&
(i)
&
\a' \leq \a
&
&
(ii)
&
b' \geq b
&
&
(iii)
&
c'\geq c,
&
\\
%%end line 2
(3)
&
\textit{if}
&
b'=b
&
\textit{and}
&
\b'=\b,
&
&
\textit{then}
&
&
(i)
&
\a' \leq \a
&
&
(ii)
&
\phantom{.} a' \leq a.
&
&
&
&
%%end line 3
\end{array}
\end{equation*}
\end{proposition}

In (3) there is no general comparison between $c'$ and $c$. We shall, however, prove the inequality $c' \leq c$ under additional hypotheses in Proposition \ref{p:CompRightTriangle2}.
\begin{proof}
Inequality (1)$(iii)$ holds by part B) of Toponogov's theorem. For the other inequalities in (1) we consider the auxiliary geodesic triangle $D^*$ in the hyperbolic plane with sides $a^* =a$, $b^* =b$, $c^* =c$, opposite vertices $p^*, q^*, r^*$ and respective angles $\a^*, \b^*, \g^*$. By part A) of Toponogov's theorem we have
\begin{equation*}
\a^* \leq \a, \quad \b^* \leq \b,  \quad \g^* \leq \frac{\pi}{2}.
\end{equation*}
Now rotate side $a^*$ around vertex $r^*$ away from $D^*$ (Figure \ref{fig:CompRightTriangles}) until it becomes orthogonal to $b^*$ so as to get the right-angled geodesic triangle $D'$ in the hyperbolic plane with sides $a' = a^*$, $b' = b^*$, $c'$, the opposite vertices $p' = p^*$, $q'$, $r' = r^*$ and respective angles $\a', \b', \pi/2$. Since $\b^* \leq \b \leq \pi/2$ side $c'$ from $p^*$ to $q'$ crosses side $a^*$ of $D^*$, and for this reason we have $\a' \leq \a^* \leq \a$. This is (1)$(i)$ and in the same way we get (1)$(ii)$.

%
%%% Figure 2 %%%
\begin{figure}[h!]
\begin{center}
\leavevmode
\SetLabels
(.238*.94) $p$\\
(-.011*.245) $q$\\
(.235*.245) $r$\\
(.14*.5) $D$\\
(.205*.815) $\alpha$\\
(.032*.33) $\beta$\\
(.115*.24) $a$\\
(.242*.56) $b$\\
(.07*.61) $c$\\
(.635*.07) $\mathbb{H}$\\
(.969*.07) $\mathbb{H}$\\
(.537*.625) $\alpha^*$\\
(.528*.75) $\alpha'$\\
(.504*.08) $\beta'$\\
(.610*.325) $\gamma^*$\\
(.411*.325) $\beta^*$\\
(.494*-.05) $q'$\\
(.673*.24) $r'{=}r^*$\\
(.564*.99) $p'{=}p^*$\\
(.355*.25) $q^*$\\
(.433*.6) $c^*$\\
(.483*.47) $c'$\\
(.584*.11) $a'$\\
(.618*.6) $b'{=}b^*$\\
(.538*.24) $a^*$\\
(.77*.21) $q'{=}q^\#$\\
(1.013*.21) $r'{=}r^\#$\\
(1.019*.99) $p'$\\
(1.024*.62) $p^\#$\\
(.879*.63) $c'$\\
(.911*.495) $c^\#$\\
(.836*.302) $\beta^\#$\\
(.880*.365) $\beta$\\
(.892*.21) $a'{=}a^\#{=}a$\\
(1.04*.41) $b^\#{=}b$\\
(1.017*.79) $b'$\\
(.978*.83) $\alpha'$\\
(.973*.52) $\alpha^\#$\\
(.551*.47) $D'$\\
(.943*.61) $D'$\\%
\endSetLabels
%\ShowGrid
\AffixLabels{%
%\centerline{%
\includegraphics{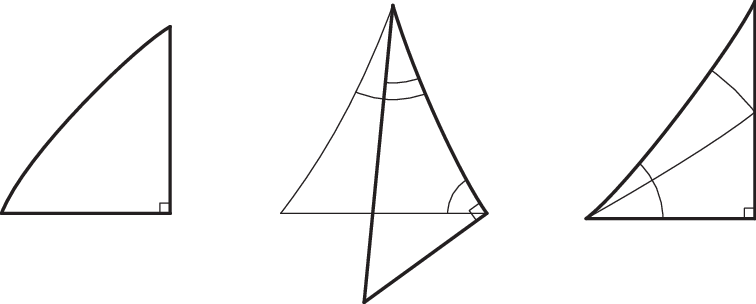}
}
%}
\end{center}
\caption{Triangle $D$ of variable curvature and comparison triangles $D', D^*$ and $D', D^\#$ in the hyperbolic plane.}
\label{fig:CompRightTriangles}
\end{figure}
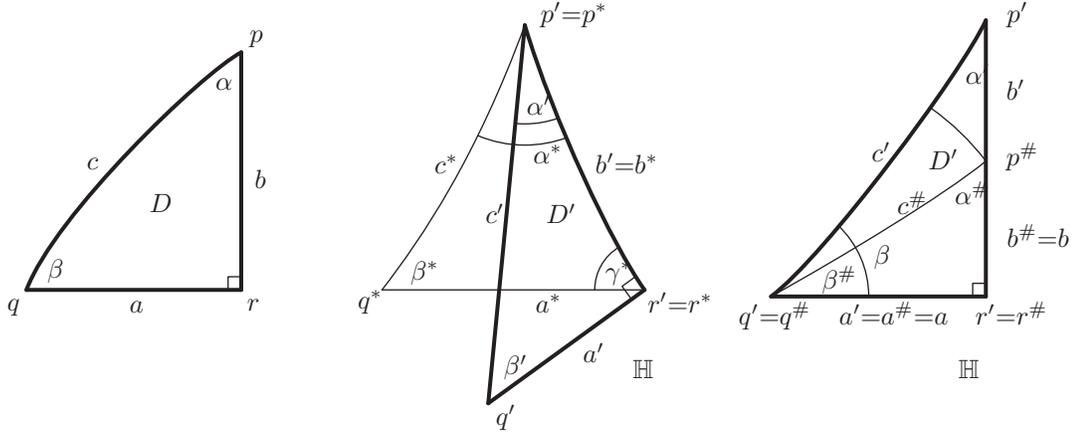
%
%%% end Figure 2 %%%

For the inequalities in (2) we take the auxiliary right-angled triangle $D^\#$ with sides $a^\# = a$, $b^\# = b$ that form a right angle at vertex $r^\#$, and angles $\a^\#$ at vertex $p^\#$ opposite to side $a^\#$ and $\b^\#$ at vertex $q^\#$ opposite to side $b^\#$ (Figure \ref{fig:CompRightTriangles}). By (1)$(i)$--$(iii)$ we have
\begin{equation*}
\a^\# \leq \a, \quad \b^\# \leq \b, \quad  c^\# \geq c.
\end{equation*}
Now we rotate side $c^\#$ around vertex $q^\#$ away from $D^\#$ until it forms an angle $\b'=\b \geq \b^\#$. If the geodesic prolongation of the so rotated side intersects the geodesic prolongation of side $b^\#$ in some point $p'$ (Figure \ref{fig:CompRightTriangles}) then the points $p', q^\#, r^\#$ form the comparison triangle $D'$ with side $c'$ from $p'$ to $q^\#=q'$, side $a'$ from $q^\#$ to $r^\#=r'$ and side $b'$ from $r^\#$ via $p^\#$ to $p'$. By our construction we have $b' \geq b^\# = b$, and by hyperbolic geometry $\a' \leq \a^\#$, $c' \geq c^\#$. Together with the preceding inequalities this proves (2)$(i)$--$(iii)$ in the intersection case. If the prolongations do not intersect each other, then $D'$ is an open-ended triangle with $c'$, $b'$ infinitely long and, by convention, $\a' = 0$. In this case (2)$(i)$--$(iii)$ are thus trivially true.

To prove (3) we use the preceding $D^\#$ with $a^\# = a$ and $b^\# = b$ a second time, but this time shrink side $a^\#$ keeping vertex $r^\#$ fixed, thereby increasing $\b^\#$ until angle $\b$ is reached at which moment we have the comparison triangle $D'$ with vertices $p' = p^\#$, $q'$ on side $a^\#$ and $r' = r^\#$. Obviously this yields $\a' \leq \a^\# \leq \a$ and $a' \leq a^\# = a$ which is (3)$(i)$--$(ii)$.
\end{proof}

We proceed with the ``missing'' (3)$(iii)$ in the above proposition.

\begin{proposition} \label{p:CompRightTriangle2}
Let $D$ with straight sides $a,b,c$ and opposite vertices $p,q,r$ be a closed Riemannian triangle with curvature $K \geq -1$. Assume that $D$ has a right angle at $r$, an acute angle $\b$ at $q$ and an angle $\a < \pi$ at $p$. Furthermore, assume that side $b$ is a shortest connection in $D$ from $p$ to side $a$.

Consider, for comparison, a right-angled geodesic triangle $D'$ in the hyperbolic plane with sides $a'$, $b'$, $c'$, opposite vertices $p',q',r'$, a right angle at $r'$ and angles $\a', \b'$ at $p'$ and $q'$. Then the following hold:
\medskip

\begin{itemize}
\item[(a)]If $\b' = \b$ and $c' = c$, then $b' \geq b$ and $a' \leq a$,
\item[(b)]if $\b' = \b$ and $b' = b$, then $c' \leq c$ and $a' \leq a$.
\end{itemize}
\end{proposition}

\begin{proof}
The major part is the first inequality in (a).
For a variant proof see Remark \ref{r:varpr}.

We begin by parametrising side $c$ by arc length $t \mapsto c(t)$, $t \in [0, L(c)]$ with $c(0) = q$, $c(L(c)) = p$. For any $t \in [0, L(c)]$ there is a shortest arc $\eta_t$ in $D$ from $c(t)$ to side $a$. This arc exists, by Lemma \ref{l:broken}, because $D$ is metrically complete. By the condition on the angles at the vertices of $D$, $\eta_t$ is a geodesic arc that meets side $a$ orthogonally. If for given $t$ there are more than one such shortest perpendiculars to $a$ we let $\eta_t$ be the one whose foot point on side $a$ is closest to vertex $r$. The hypothesis on side $b$ then implies that $b = \eta_{L(c)}$.

The idea is to use Proposition \ref{p:CompRightTriangles}(3)$(i)$ in order to show that
the corresponding perpendiculars in $D'$ grow faster than the $\eta_t$.
Unfortunately, the possible positive curvature in $D$ does not allow us to work with derivatives of length functions and we have to resort to real analysis type arguments.

%
%%% Figure 3 %%%
%
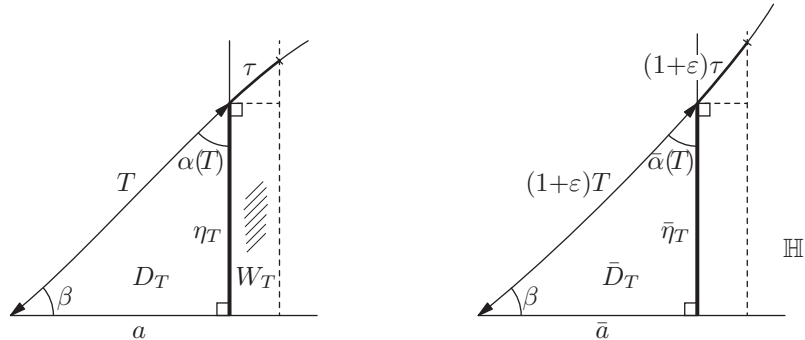
\begin{figure}[b!]
\begin{center}
\leavevmode
\SetLabels
(.15*.4) $T$\\
(.718*.4) $(1{+}\varepsilon)T$\\
(.255*.25) $\eta_T$\\
(.856*.25) $\bar{\eta}_T$\\
(.071*.036) $\beta$\\
(.670*.036) $\beta$\\
(.246*.472) $\alpha(\!T\!)$\\
(.851*.47) $\bar{\alpha}(\!T\!)$\\
(.309*.77) $\tau$\\
(.866*.778) $(1{+}\varepsilon)\tau$\\
(.182*.1) $D_T$\\
(.786*.1) $\bar{D}_T$\\
(.317*.1) $W_T$\\
(.166*-.07) $a$\\
(.762*-.07) $\bar{a}$\\
(1.008*.19) $\mathbb{H}$\\%
\endSetLabels
%\ShowGrid
\AffixLabels{%
%\centerline{%
\includegraphics{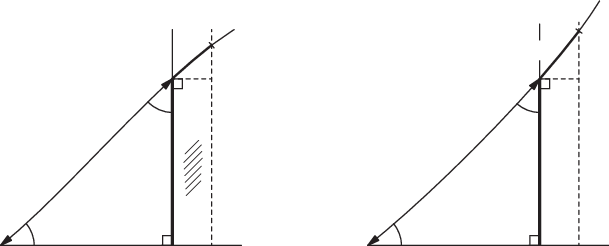}%[width=0.682\textwidth]%
}
%}
\end{center}
\caption{Comparison triangle $\bar{D}_T$ in the hyperbolic plane with angle $\bar{\a}(T) \leq \a(T)$ and height $L(\bar{\eta}_T) = L(\eta_T)$.}
\label{fig:HgrowsFaster}
\end{figure}
%
%%% end Figure 3 %%%
%

Thus, let $\e >0$ be arbitrarily small and let $\bar{D}$ be the slightly bigger right-angled triangle in the hyperbolic plane with angle $\bar{\b} = \b$ at vertex $\bar{q}$ and hypothenuse $\bar{c}$ of length $L(\bar{c}) = (1+ \e) L(c)$. We parametrize $\bar{c}$ with constant speed $(1+\e)$ in the form $t \mapsto \bar{c}(t)$, $t \in [0,L(c)]$, and let for $t \in [0, L(c)]$, $\bar{\eta}_t$ be the perpendicular in $\bar{D}$ from $\bar{c}(t)$ to side $\bar{a}$. Note that $\bar{\eta}_{L(c)}$ is side $\bar{b}$.
For $t$ near 0 we have approximately $L(\bar{\eta}_t) = (1+\e)L(\eta_t)>L(\eta_t)$. Let now $T$ be the maximum of all $t \in [0, L(c)]$ that satisfy $L(\bar{\eta_t})\geq L(\eta_t)$. We claim that $T = L(c)$, which then implies that $\bar{b}\geq b$.

Suppose to the contrary that $T < L(c)$. By continuity then $L(\bar{\eta}_T) = L(\eta_T)$. For the triangles $D_T$ cut away from $D$ by $\eta_T$ and $\bar{D}_T$ cut away from $\bar{D}$ by $\bar{\eta}_T$, we claim that we then have
\begin{equation}\label{eq:aT}
\bar{\a}(T) \leq \a(T),
\end{equation}
where $\a(T)$ and $\bar{\a}(T)$ are respectively the angle of $D_T$ at vertex $c(T)$ and the angle of $\bar{D}_T$ at vertex $\bar{c}(T)$.
Indeed, if $\a(T) \leq \pi/2$, then inequality \eqref{eq:aT} holds by Proposition \ref{p:CompRightTriangles}(3)(i); if
$\a(T) > \pi/2$, the inequality is trivially true.

We shall estimate the lengths of the neighboring perpendiculars using Fermi coordinates along $\eta_T$.

For this we recall that, by our definition of a Riemannian polygon, there exists a (not necessarily complete) Riemannian surface $R$ without boundary that contains $D$ as a bordered Riemannian subsurface. Outside $D$ the curvature does not need to have the lower bound $-1$, but we take $R$ such that we have an upper bound for the absolute curvature
\begin{equation}\label{eq:upperK}
\abs{K} \leq \kappa \quad \text{on $R$}.
\end{equation}

This allows us to extend $\eta_T$, without changing its name, to a simple unit speed geodesic $\th \mapsto \eta_T(\th)\in R$, $\th \in ({-} \omega, L(\eta_T) + \omega)$, for some $\omega > 0$. Furthermore, there exists $w > 0$ and an open neighborhood $W_T$ of $\eta_T$ in which Fermi coordinates based on $\eta_T$ can be introduced such that in these coordinates
\begin{equation}\label{eq:KT}
W_T = \{(r,\th) \mid  -w < r < w,\; -\omega < \th < L(\eta_T) + \omega \}
\end{equation}
and we have a unit speed parametrization of $\eta_T$ in the form $\eta_T(\th) = (0,\th)$, $\th \in ({-}\omega, L(\eta_T)+\omega)$, with $\eta_T(0)=(0,0)$ the footpoint of $\eta_T$ on side $a$ of $D$, while the curves
$r \mapsto (r,\th)$, $r \in ({-}w,w)$, are unit speed geodesic arcs of length $2w$ intersecting $\eta_T$ orthogonally at $\eta_T(\th)$ (see Figure \ref{fig:HgrowsFaster}).

The expression for the metric tensor $ds^2$ in these coordinates on $W_T$ is
\begin{equation}\label{eq:GT}
ds^2 = dr^2 + G_T(r,\th)^2 d\th^2,
\end{equation}
with a function $G_T$ that satisfies $\cos \kappa r \le G_T(r,\th) \le \cosh \kappa r$, where $\kappa$ is the absolute curvature bound from \eqref{eq:upperK}.
(For a reference, see e.g.\ \cite[p.\ 247, (14)]{C}; we shall also use this metric tensor in the proof of Corollary \ref{c:Dclosed} and it shows up again in \eqref{eq:Geta} for the study of distance curves.)

Hence, there exists a positive constant $M$ such that for small $\tau > 0$ we have
\begin{equation*}
L(\eta_{T+\tau}) \leq L(\eta_{T}) + \tau \cos(\a(T)) + M \tau^2
\end{equation*}
using that the right-angled triangle with hypothenuse of length $\tau$ from $c(T)$ to $c(T+\tau)$ and angle $\pi/2 - \a(T)$ at vertex $c(T)$ has height $\tau \cos(\a(T)) + O(\tau^2)$. The corresponding triangle in $\bar{D}$ from $\bar{c}(T)$ to $\bar{c}(T+\tau)$ has hypothenuse of length $(1+\e)\tau$ and angle $\pi/2 - \bar{\a}(T)$. By hyperbolic geometry we have, without an $O$-term appearing,
\begin{equation*}
L(\bar{\eta}_{T+\tau}) \geq L(\eta_{T}) + (1+\e) \tau \cos(\bar{\a}(T)).
\end{equation*}
By \eqref{eq:aT} $\cos(\a(T)) \leq \cos(\bar{\a}(T))$. If we now take $\tau < \frac{\e}{M}\cos(\bar{\a}(T))$ we get $L(\bar{\eta}_{T+\tau}) >
L(\eta_{T+\tau})$, a contradiction. This finishes the proof that $T = L(c)$ and thus $\bar{b} \geq b$. Letting $\e$ converge to 0, $\bar{b}$ converges to $b'$ and we get $b' \geq b$. This is the first inequality in point (a) of the proposition.

For the second inequality we consider, once again, the comparison triangle $D^\#$ in $\HH$ with $a^\# = a$ and $b^\# =b$ and hypothenuse $c^\#$.
By point B) in Toponogov's Theorem \ref{t:top} (which makes no assumptions about angle $\a$) we have $c^\# \geq c$.
If we had $a' > a$, then since $b' \geq b$ we would get $c' > c^\#$, by hyperbolic geometry, and so $c' > c$, a contradiction. Hence, $a' \leq a$.

For (b) we use the auxiliary right-angled triangle $D''$ in $\HH$ with hypothenuse $c'' = c$, angle $\b'' = \b$ at vertex $q''$ and sides $a''$ at $q''$, $b''$ opposite to $q''$. By (a) it satisfies $a'' \leq a$, $b'' \geq b$. Now we shrink side $a''$ keeping vertex $q''$ and the angle $\b''$ fixed thereby shrinking both, $c''$ and the orthogonal $b''$ until the latter reaches length $b$ at which moment we have obtained the comparison triangle $D'$ with $\b' = \b'' = \b$ the short sides $a' \leq a'' \leq a$, $b' = b$ and hypothenuse $c' \leq c'' = c$. Hence (b). The proof of the proposition is now complete.
\end{proof}

\begin{remark}\label{r:varpr}
\emph{Another proof of the first inequality in (a), similar in length, would be to symmetrize triangle $D$ across side $a$, smooth the metric and apply part (B) of Toponogov's theorem to a triangle that lies arbitrarily closed to the symmetrization of $D$.}
\end{remark}

For the next corollary we recall the following trigonometric formula for a right-angled hyperbolic triangle with short sides $a$, $b$ and opposite angles $\b$, $\a$:

\begin{equation}\label{eq:trig1}
\raisebox{-1.8em}{
%beginFigure
\SetLabels
(.50*.06) $\scriptstyle a$\\
\L(.87*.45) $\scriptstyle L(b) = B^{\HH}(\beta,L(a))$\\
(.42*.48) $\scriptstyle c$\\
(.76*.65) $\scriptstyle \alpha$\\
(.265*.17) $\scriptstyle \beta$\\
\endSetLabels
%\ShowGrid
\AffixLabels{%
{\includegraphics{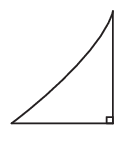}}%
}%endFigure
}
\rule{16ex}{0pt}
\tanh(L(b)) = \tan(\b)\sinh(L(a))
\end{equation}

Accordingly, we define the two functions
\begin{equation}\label{eq:trig2}
B^{\HH}(\varphi,x) := \arctanh(\tan(\varphi)\sinh(x)),
\quad
\b_{\infty}^{\HH}(x) := \arctan\Big(\frac{1}{\sinh(x)}\Big).
\end{equation}
The height over side $a$ of the triangle is then given by $L(b) =B^{\HH}(\b,L(a))$. For $\b \to \b_{\infty}^{\HH}(L(a))$ we have $L(b) \to \infty$ and in the limit the triangle is singular with an ideal vertex at infinity.

\begin{corollary} \label{c:Dopen}
Let $D$ be an open-ended Riemannian triangle with curvature $K \geq -1$. Let $a$ with endpoints $q$, $r$ be the finite side of $D$ and assume that $D$ has an acute angle $\b$ at $q$ and a right angle at $r$. Then $\b \geq \b_{\infty}^{\HH}(L(a))$.
\end{corollary}

\begin{proof}
We first remark that it suffices to prove the corollary under the assumption that $D$ has straight sides. Indeed, we may shrink, if necessary, $D$ to a smaller open ended triangle $D' \subset D$ that satisfies this additional condition by first replacing  the infinite side $c$ at $q$ by a straight geodesic ray $c'$ emanating from $q$ into $D$, then the infinite side $b$ at $r$ by a straight geodesic ray $b'$ emanating orthogonally from $a$ into $D$ without intersecting $c'$ and finally side $a$ by a shortest orthogonal $a'$ from $q$ to $b'$. The straightening procedures that lead from $c$ to $c'$ and from $b$ to $b'$ will be detailed in the remark below. The acute angle $\b'$ of $D'$ at $q$ satisfies $\b' \leq \b$ and for the finite side we have $L(a') \leq L(a)$. Hence, if the claim holds for $D'$ we also have $\b \geq \b' \geq \b_{\infty}^{\HH}(L(a')) \geq \b_{\infty}^{\HH}(L(a))$. We may thus assume that $D$ is $D'$.

We parametrize the infinite side $c$ at $q$ by arc length $t \mapsto c(t)$, $t \in [0, \infty)$, with $c(0) = q$. As in the proof of the preceding proposition there exists for each $t \in [0, \infty)$ a shortest perpendicular $\eta_t$ in $D$ from $c(t)$ to side $a$. Since $c$ is straight we have
\begin{equation*}
L(\eta_t) \geq t - L(a).
\end{equation*}
Let $D'_t$ with sides $a'_t$, $b'_t$, hypothenuse $c'_t$ and angle $\b'_t$ between $a'_t$ and $c'_t$ be the right-angled comparison triangle in $\HH$ that has $\b'_t = \b$ and $L(c'_t) = t$. By Proposition \ref{p:CompRightTriangle2}(a) we have $L(b'_t) \geq L(\eta_t)$ and $L(a'_t) \leq L(a)$. It follows, together with \eqref{eq:trig1}, that for any $t >L(a)$,
\begin{equation*}
\b = \arctan \frac{\tanh L(b'_t)}{\sinh L(a'_t)} \geq \arctan \frac{\tanh L(\eta_t)}{\sinh L(a)}\geq \arctan \frac{\tanh(t-L(a))}{\sinh L(a)}.
\end{equation*}
With $t \to \infty$ we have $\tanh(t-L(a)) \to 1$ and the corollary follows.
\end{proof}

\begin{upremark}[straightening procedures]\label{r:straighten}%
(a) Let $t \mapsto c(t)$, $t \in [0, \infty)$, with $c(0) = q$ be the arc length parametrization of side $c$ of the open ended triangle $D$ as in Corollary \ref{c:Dopen} and assume that $c$ is not straight. By Lemma \ref{l:broken}, there exists for each integer $n \geq 1$ a minimizing geodesic arc $t \mapsto c_n(t) \in D$, $t \in [0, \dist(c(0), c(n))]$, parametrized by arc length, going from $c_n(0) = c(0)=q$ to $c(n)$. Since $c$ is not straight there exists $n_0$ such that for all $n \geq n_0$ the arc $c_n$ is shorter than $n$ and meets side $c$ only at its endpoints. Furthermore, by Lemma \ref{l:nocut}, for $m > n \geq n_0$, the arcs $c_m$ and $c_n$ meet only at $q$. The initial tangent vectors $\bdot{c}_n(0)$ form therefore monotone decreasing angles with the second side $a$ of $D$ at $q$ and converge to a limit unit tangent vector $v$ at $q$ as $n \to \infty$. Hence, the $c_n$ converge uniformly on compact sets to a geodesic ray $c': [0,\infty) \rightarrow D$ which is minimizing because the $c_n$ are minimizing and escapes into the end of $D$.

\smallskip

(b) Let $t \mapsto b(t)$, $t \in [0, \infty)$, with $b(0) = r$ be the unit speed parametrizsation of side $b$ of $D$ which forms a right angle with side $a$ at $r$ and assume that $b$ is not straight. Assume, more generally, that for some $n_0$ and hence, for all integers $n \geq n_0$ the minimal geodesic arc $b_n : [0, \dist(a,b(n))]$ from $a$ to $b(n)$ (Lemma \ref{l:broken}), parametrized by arc length, has length less than $n$. Being minimizing these arcs pairwise do not intersect each other and converge uniformly on compact sets to a minimizing geodesic ray $b' : [0, \infty) \to  D$ emanating perpendicularly from $a$ and escaping into the end of $D$.

\smallskip

(c) Minimizing rays will also be constructed in similar form in the proofs of Lemma \ref{l:openquadri} and Propositions \ref{p:genfunnel}, \ref{p:descomposicion3}.

\end{upremark}

Corollary \ref{c:Dopen} has a compact version:

\begin{corollary} \label{c:Dclosed}
Let $D$ be a compact Riemannian quadrangle of curvature $K \geq -1$, with consecutive sides $a,b,c,d$, having three right angles and an acute angle $\b$ between sides $a$ and $d$. Then $\b \geq \b_{\infty}^{\HH}(L(a))$, $\b \geq \b_{\infty}^{\HH}(L(d))$.
\end{corollary}

\begin{proof}
By symmetry we only have to prove the first inequality. In order to reduce Corollary \ref{c:Dclosed} to Corollary \ref{c:Dopen} we attach a Euclidean strip $E = \{(x,y) \in \RR^2 \mid -\infty < x \leq 0,\; 0 \leq y \leq L(c) \}$ along side $c$ to $D$ in order to get an open ended Riemannian triangle with a right angle at vertex $r$ between $a$ and $b$ and an acute angle $\b$ at vertex $q$ between $a$ and $d$. In the particular case where the standard Euclidean metric tensor $ds^2 = dx^2 + dy^2$ on $E$ matches smoothly with the Riemannian metric tensor on $D$ along the common side $c$ we can apply Corollary \ref{c:Dopen} to $E \cup D$ and are done.

For the general case we apply a smoothing interpolation procedure. We carry this out using Fermi coordinates based on $c$ in some small neighborhood $W = \{z \in D \mid \dist(z, c) < w \}$. In its own coordinates this neighborhood has the description $W = \{(x,y) \mid 0 \leq x < w,\; 0 \leq y \leq L(c) \}$ and the metric tensor, as in \eqref{eq:GT}, is of the form
\begin{equation*}
ds^2 = dx^2 + G(x,y)^2 dy^2.
\end{equation*}
Note that
\begin{equation*}
G(0,y) = 1, \quad  G'(0,y) = 0, \quad K(x,y) = -\frac{G''(x,y)}{G(x,y)},
\end{equation*}
where $K$ is the Gauss curvature and the primes are the partial derivatives with respect to the first variable. The smoothing is performed by multiplying $G(x,y)-1$ with an attenuating function $\phi$ which is equal to 0 in a neighborhood of $x=0$ and equal to 1 in a neighborhood of $x =w$. Since we want to respect the lower curvature bound the choice of $\phi$ needs some care.

We first take any smooth monotone increasing function $\chi : \RR \to [0,1] $ which is constant equal to 0 for $x \leq \frac{1}{4}$ and constant equal to 1 for $x \geq \frac{3}{4}$. Then we define $\phi$ as follows, where $\d>0$ is an arbitrarily small parameter and we set $v := \d^{\frac{1}{\d}}$:
\begin{equation*}
\phi(x) = \chi(\tfrac{x}{v})x^{\d} + \chi(\tfrac{x}{w}) (1-x^{\d}).
\end{equation*}
Note that $\phi(x) = 0$ for $x\leq \frac{1}{4}v$, $\phi(x) = 1$ for $x \geq \frac{3}{4}w$ and $0\leq \phi(x) \leq 1$ in between.

The attenuated metric tensor is defined by $\tilde{G}(x,y) = 1 + \phi(x)g(x,y)$, where $g(x,y) = G(x,y) -1$. Its Gauss curvature $\tilde{K}(x,y)$ becomes
\begin{equation*}
\tilde{K}(x,y) = K(x,y) \phi(x) \frac{G(x,y)}{\tilde{G}(x,y)}     - \frac{2\phi'(x)g'(x,y)}{\tilde{G}(x,y)}     - \frac{\phi''(x)g(x,y)}{\tilde{G}(x,y)}.
\end{equation*}
Let now $a$ be a common upper bound for $\chi'(x)$ and $\abs{\chi''(x)}$ and let $b$ be an upper bound for $\abs{g''(x,y)}$ so that $\abs{g'(x,y)} \leq b x$ and $\abs{g(x,y)} \leq \frac{1}{2}b x^2$. Reducing the size of $w$, if necessary, we may assume that, in addition, $\abs{g(x,y)} < \frac{1}{2}$ and also $w< \frac{1}{4}$. A direct computation then yields the following:
\begin{equation*}
\frac{G(x,y)}{\tilde{G}(x,y)} < 1 + \d \cdot \frac{b}{2}, \quad%
\abs{\phi'(x)g'(x,y)} < \d \cdot \kappa, \quad%
\abs{\phi''(x)g(x,y)} < \d\cdot \kappa,
\end{equation*}
for some constant $\kappa$ just depending on $a,b,w$.
For the verification of the first inequality, note that if $g\le 0$, then
\begin{equation*}
\frac{G(x,y)}{\tilde{G}(x,y)}
= \frac{1+g}{1+\phi g}
\le 1,
\end{equation*}
and if $g\geq 0$, then
\begin{equation*}
\frac{G(x,y)}{\tilde{G}(x,y)}
= \frac{1+g}{1+\phi g}
= 1 + \frac{1-\phi}{1+\phi g}\, g
\le 1 + (1-\phi) g
\le 1 + (1-\phi) \frac12 \,bx^2.
\end{equation*}
For the remaining inequalities one has to deal with the cases $x \leq v$ and $x > v$ separately; for $x < v$ use is made of the fact that $v^{\d}=\d$ and for $x > v$ the expressions for the upper bounds are simplified by the fact that $x(1-x^{\d}) < \frac{\d}{e}$.

It now follows that for $\d \to 0$ the greatest lower bound for $\tilde{K}$ converges to the greatest lower bound for $K$. We may thus apply Corollary \ref{c:Dopen} to the metric defined by $\tilde{G}$ on $E  \cup D$ and get Corollary \ref{c:Dclosed} in the limit as $\d$ goes to zero.
\end{proof}

\begin{figure}[!h]
\begin{center}
\leavevmode
\SetLabels
(.85*1.06) $c$\\
(.85*-.12) $b$\\
(-.005*.48) $a$\\
(-.015*.999) $q$\\
(-.015*-.02) $r$\\
(.195*.90) $\varphi_2$\\
(.404*.90) $\varphi_4$\\
(.578*.90) $\varphi_6$\\
(.682*.90) $\varphi_8$\\
(.064*.068) $\varphi_{1}$\\
(.300*.068) $\varphi_{3}$\\
(.50*.068) $\varphi_{5}$\\
(.629*.068) $\varphi_{7}$\\
(.728*.068) $\varphi_{9}$\\
(.133*.48) $s_1$\\%
(.248*.48) $s_2$\\
(.355*.48) $s_3$\\
(.454*.48) $s_4$\\
(.547*.48) $s_5$\\
(.15*1.06) $p_2$\\
(.364*1.06) $p_4$\\
(.55*1.06) $p_6$\\
(.656*1.06) $p_8$\\
(.242*-.095) $p_3$\\
(.455*-.095) $p_5$\\
(.592*-.095) $p_7$\\
(.70*-.095) $p_9$\\
(.90*.47) $\Delta$\\
\endSetLabels
%\ShowGrid
\AffixLabels{%
%\centerline{%
\includegraphics{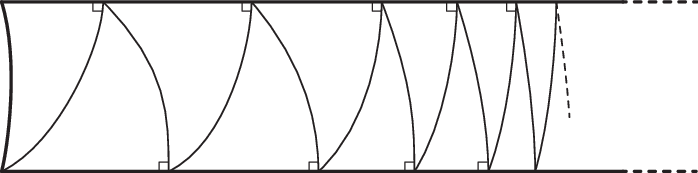}%[width=0.712\linewidth]
}
%}
\end{center}
\caption{Flip flop procedure on an open-ended triangle $\D$ with straight sides.}
\label{fig:flipflop}
\end{figure}
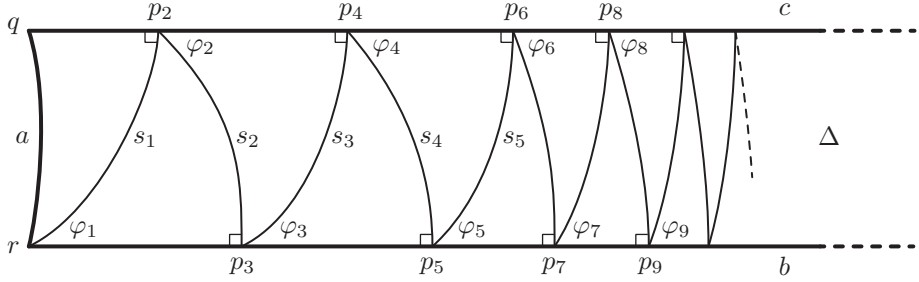

\textbf{Flip flop procedure.} We apply the preceding results to the following  procedure inspired by the Birkhoff curve straightening process. Consider an open-ended Riemannian triangle $\D$ of curvature $K \geq -1$ with straight sides $a$, $b$, $c$, where $b$, $c$ are geodesic rays and $a$ is a geodesic segment meeting $c$ at vertex $q$ and $b$ at vertex $r$. We further assume that the angles at $q$ and $r$ are smaller than $\pi/2$. The procedure is the following. From $r$ we draw the shortest connection $s_1$ in $\D$ to side $c$. If there are several shortest connections we choose the one whose foot point lies closest to $q$. Note that $s_1$ is a geodesic segment meeting $c$ orthogonally. From its foot point we draw a shortest connection $s_2$ to $b$, again the one with the foot point closest to $r$ if there are several such. In this way we continue getting a succession of perpendiculars $s_3$, $s_4$, $s_5$, and so on. We denote by $\varphi_k$
the angle pointing outwards between $s_k$ and $b$ if $k$ is odd, respectively $s_k$ and $c$ if $k$ is even (Figure \ref{fig:flipflop}).

%
%%% Figure 4 %%%
%
%
%%% end Figure 4 %%%
%

\begin{proposition} \label{p:flipflop}
Under the above conditions the angles $\varphi_k$ in the flip flop procedure converge to $\pi/2$.
\end{proposition}

\begin{proof}
For $k=2,3, \dots$ we let $p_k$ be the foot point of $s_{k-1}$ on $c$ respectively, $b$. Since the $s_k$ are minimal connections they are straight and accordingly any $s_k$ intersects the union $s_1 \cup \dots \cup s_{k-1}$ only in $p_k$. This further implies that all $\varphi_k$ are acute and the distances from $p_{2j-1}$ to $r$ along $b$, respectively from $p_{2j}$ to $q$ along $c$ are monotone increasing. We now argue in three different cases with three different arguments.

\emph{Case 1:} the $p_k$ remain within bounded distance from side $a$.
For this case no curvature bounds come into play:  the $p_{2j-1}$ on $b$ converge to a limit point $\bar{r}$ on $b$ and the $p_{2j}$ on $c$ to a limit point $\bar{q}$ on $c$.
The $s_k$ then converge to a geodesic arc $\bar{s}$ that connects $\bar{p}$ with $\bar{q}$.
By continuity $\bar{s}$ meets $b$ and $c$ under right angles and the $\varphi_k$ converge to these right angles.

\emph{Case 2:} the $p_k$ go out to infinity and the lengths of the segments $s_k$ have a positive lower bound. Here we shall use Proposition \ref{p:CompRightTriangle2}.

Since the angles $\varphi_k$ are acute and the $s_k$ are shortest connections the lengths $L(s_k)$ are monotone decreasing and converge to some limit $L >0$. Consider the right-angled triangle $T_k$ with hypothenuse $s_k$, side $s_{k+1}$ and opposite angle $\varphi_k$ at vertex $p_k$. Since $s_{k+1}$ is a shortest perpendicular there exists by Proposition \ref{p:CompRightTriangle2}(a) a right-angled hyperbolic comparison triangle $T'_k$ with hypothenuse $s'_k$, side $s'_{k+1}$ and opposite angle $\varphi'_k$ satisfying
$L(s'_k)= L(s_k)$,  $\varphi'_k = \varphi_k$, $L(s'_{k+1}) \geq L(s_{k+1})$. By the hyperbolic sine formula applied to $T'_k$ it follows that
\begin{equation*}
\sin \varphi_k = \frac{\sinh L(s'_{k+1})}{\sinh L(s_k)} \geq \frac{\sinh L(s_{k+1})}{\sinh L(s_k)}.
\end{equation*}
As $k \to \infty$ the right hand side converges to 1 and hence, $\varphi_k \to \pi/2$.

\emph{Case 3:} the $p_k$ go out to infinity and the lengths of the segments $s_k$ converge to 0.
By Corollary \ref{c:Dopen} with $D$ the open-ended triangle formed by $s_k$ and the parts of $b$, $c$ that go from $s_k$ to infinity and with $\b = \varphi_k$ we have
\begin{equation*}
\varphi_k \geq \b_{\infty}^{\HH}( L(s_k)) = \arctan \frac{1}{\sinh L(s_k)}
\end{equation*}
(see \eqref{eq:trig2}). As $L(s_k)\to 0$ it follows again that $\varphi_k \to \pi/2$.
\end{proof}

\begin{lemma}\label{l:openquadri}
Let $G$ be an open-ended geodesic quadrilateral with curvature $K \geq -1$ and three angles less than or equal to $\pi/2$, and let $a_1$, $a_2$ be the finite sides. Then $\sinh L(a_1) \sinh L(a_2) \geq 1$.
\end{lemma}

\begin{proof}
Let $a_5, a_3$ be the infinite sides adjacent respectively to $a_1$ and $a_2$ (see Figure \ref{fig:openQuadrilateral}). We parametrize $a_5$ by arc-length $t \mapsto a_5(t)$, $t \in [0, \infty)$ and draw, for any $t \in [0, \infty)$ a shortest arc $a_{5,t}$ from $a_5(t)$ to $a_1$. If there are several such arcs we take the one whose foot point on $a_1$ has minimal distance to the common vertex $q$ of $a_1$ and $a_2$. The arcs meet $a_1$ orthogonally and, by the same arguments as in Remark \ref{r:straighten}, they converge uniformly on compact sets to a minimizing geodesic ray $a_5'$ in $G$ orthogonal to $a_1$ with initial point different from $q$. We next let $a_1'$ be a shortest arc from $a_5'$ to $a_2$ in $G$ and then $a_2'$ a shortest arc from $a_1'$ to side $a_3$. In the open-ended quadrilateral $G'$ whose sides are part of $a_5'$, $a_1'$, $a_2'$ and $a_3$ we repeat the construction that lead to $a_5'$ in order to get a minimizing geodesic ray $a_3'$ in $G'$ emanating orthogonally from $a_2'$. In this way we get an open-ended right-angled quadrilateral $\overline{G} \subseteq G$ with finite sides $\bar{a}_1$ on $a_1'$, $\bar{a}_2$ on $a_2'$ and infinite sides $\bar{a}_3=a_3'$, $\bar{a}_5$ on $a_5'$. Note that since $a_1'$, $a_2'$, $a_3'$, $a_5'$ are minimizing, all sides of $\overline{G}$ are straight.

%
%%% Figure 5 %%%
%
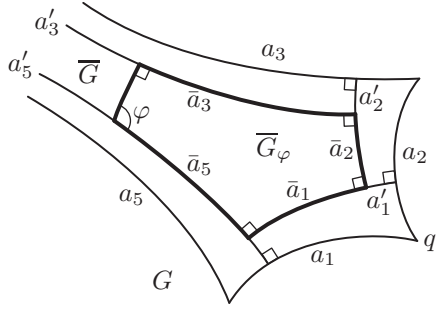
\begin{figure}[t!]
\begin{center}
\leavevmode
\SetLabels
(.75*.14) $a_1$\\
(.891*.316) $a_1'$\\
(.872*.66) $a_2'$\\
(.982*.48) $a_2$\\
(.625*.82) $a_3$\\
(.437*.65) $\bar{a}_3$\\
(.437*.44) $\bar{a}_5$\\
(.265*.34) $a_5$\\
(.797*.49) $\bar{a}_2$\\
(.688*.355) $\bar{a}_1$\\
(.290*.616) $\varphi$\\
(.625*.49) $\overline{G}_{\varphi}$\\
(.156*.73) $\overline{G}$\\
(.343*.05) $G$\\
(-.016*.77) $a_5'$\\
(.052*.91) $a_3'$\\
(1.019*.19) $q$\\
\endSetLabels
%\ShowGrid
\AffixLabels{%
%\centerline{%
\includegraphics{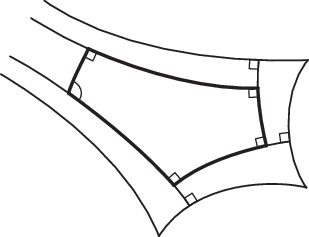}%[width=0.712\linewidth]
}
%}
\end{center}
\caption{Pentagon $\overline{G}$ with straight sides and angle $\varphi > \pi/2$ in the open-ended quadrilateral $G$.}
\label{fig:openQuadrilateral}
\end{figure}
%
%%% end Figure 5 %%%
%

It follows from the construction that
\begin{equation}\label{eq:a1a2}
L(\bar{a}_1) \leq L(a_1') \leq L(a_1), \quad L(\bar{a}_2) \leq L(a_2') \leq L(a_2).
\end{equation}
If $\bar{a}_3$ and $\bar{a}_5$ have a common orthogonal, then the lemma follows from Lemma \ref{l:sinh} and \eqref{eq:a1a2}. Otherwise we apply the flip flop procedure to the couple of sides $\bar{a}_3, \bar{a}_5$ to get a closed geodesic pentagon $\overline{G}_{\varphi}$ (Figure \ref{fig:openQuadrilateral}) with four right angles and an obtuse angle $\varphi$ with $\varphi$ as close to $\pi/2$ as we wish. $\overline{G}_{\varphi}$ shares sides $\bar{a}_1$ and $\bar{a}_2$ with $\overline{G}$ and has an opposite side $\bar{a}_{\varphi}$ orthogonal to $\bar{a}_3$ and forming the angle $\varphi$ with $\bar{a}_5$. By Lemma \ref{l:sinh2} and \eqref{eq:a1a2} we have $\sinh L(a_1) \sinh L(a_2) \geq \sin \varphi$ and the lemma follows with $\varphi \to \pi/2$, by Proposition \ref{p:flipflop}.
\end{proof}

\section{Proof of Theorems \ref{t:a} and \ref{t:b}}
\label{sec:Proof}

Since $k d_\g$, $k L(\g)$, $k L(\mu)$ etc. are invariant when the metric is multiplied by a constant, we may assume that $k=1$, i.e. $K \ge -1$ from now on.

We begin with a particular case that includes the generalized funnels and generalized cusps with geodesic boundary.

\begin{proposition} \label{p:genfunnel}
Let $F$ be a collared end with curvature $K \ge -1$ and simple closed boundary geodesic $\g_0$,
and let us define for $d>0$ the subsets
\begin{equation*}
Z(d) = \big\{ x\in F \,|\, \dist(x, \g_0) < d  \, \big\} .
\end{equation*}
Then, $Z(d)$ is topologically equivalent to $\SS^1 \times [0,1)$ for every
\begin{equation*}
d \le d_0 :=\arcsinh \cosech \Big( \frac{1}2\, L(\g_0)\Big).
\end{equation*}
\end{proposition}

\begin{figure}[h!]
\hspace{-0.25\textwidth}
\begin{minipage}[b]{0.55\textwidth}
\centering
\leavevmode
\SetLabels
(.362*.11) $a$\\
(.518*.099) $a_3$\\
(.615*.108) $a_2$\\
(.713*.130) $a_1$\\
(.283*.500) $\omega_1$\\
(.365*.480) $\omega_2$\\
(.444*.470) $\omega_3$\\
(.495*.72) $b$\\
(.589*.66) $b_3$\\
(.663*.63) $b_2$\\
(.725*.61) $b_1$\\
(.497*.90) $g$\\
(.117*.71) $\gamma_0$\\%
\endSetLabels
%\ShowGrid
\AffixLabels{%
\includegraphics{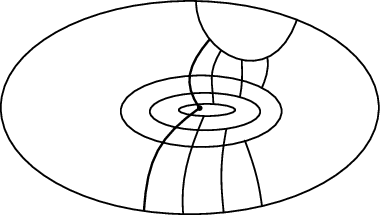}%[width=0.815\linewidth]
}
\caption{Construction of $a$ and $b$}
\label{fig:figura4}
\end{minipage}
\hspace{-0.1\textwidth}
\begin{minipage}[b]{0.55\textwidth}
\centering
\leavevmode
\SetLabels
(.50*.24) $A$\\
(.204*.60) $g_1$\\
(.829*.60) $g_2$\\
(.381*.425) $\alpha(\!t\!)$\\
(.461*.89) $h_A$\\
(.50*-.09) $\alpha(\!\tau\!)$\\
\L(.53*.48) $\alpha(\!t_0\!)$\\
\L(.53*.70) $\alpha(\!t_1\!)$\\
\endSetLabels
%\ShowGrid
\AffixLabels{%
\includegraphics{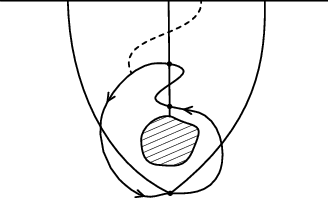}%[width=0.703\linewidth]
}
%\newline
\caption{Construction of $g$}
\label{fig:dArch}
\end{minipage}
\hspace{-0.25\textwidth}
\end{figure}

\begin{proof}
The result is clear if $K \le 0$.
Let us prove it for $K \ge -1$.

Seeking for a contradiction we assume that $Z(d)$ is not topologically equivalent to $\SS^1 \times [0,1)$ for some $d \le d_0$. Then $F \setminus Z(d)$ has an infinite closed connected component $\L_d$ and at least one compact connected component. By Lemma \ref{l:dArch} below there exists a geodesic arc $g$ of length $L(g) < 2d \leq 2d_0$ in $F$ that meets $\g_0$ orthogonally at its endpoints (Figure \ref{fig:figura4}). Thus, $g$ defines a right-angled $2$-gon in $F$ and
\begin{equation}\label{eq:r2}
\sinh \Big( \frac{1}2\, L(\g_0)\Big) \sinh \Big( \frac{1}2\, L(g)\Big)
< \sinh \Big( \frac{1}2\, L(\g_0)\Big) \sinh d_0
= \sinh \Big( \frac{1}2\, L(\g_0)\Big) \cosech \Big( \frac{1}2\, L(\g_0)\Big)
= 1.
\end{equation}

We now take an infinite sequence $\omega_1, \omega_2, \omega_3, \dots $, of simple closed curves in $F$ homotopic to $\g_0$ such that the distance from $\omega_n$ to $\g_0$ is greater than $d_0$ and goes to infinity as $n \to \infty$.

By Lemma \ref{l:broken},
for each $n$, we can consider minimizing geodesics $a_n$ from $\g_0$ to $\omega_n$, and $b_n$ from $g$ to $\omega_n$
(see Figure \ref{fig:figura4}). We let these geodesics be parametrized by arc length with initial points $a_n(0)$, $b_n(0)$ on $\g_0$ and denote by $\bdot{a}_n(0)$, $\bdot{b}_n(0)$ the initial tangent vectors.

By the minimality of $a_n$ and $b_n$, they are pairwise disjoint for each fixed $n$ (Lemma \ref{l:nocut}).
Since $\g_0$ is a compact set, there exist two convergent subsequences $\{a_{n_k}(0)\}_{k=1}^\infty$ and $\{b_{n_k}(0)\}_{k=1}^\infty$;
since the tangent vectors $\bdot{a}_{n_k}(0)$ and $\bdot{b}_{n_k}(0)$ are orthogonal to $\g_0$ for every $k$, $\{a_{n_k}\}_{k=1}^\infty$ (respectively, $\{b_{n_k}\}_{k=1}^\infty$) converges uniformly on compact sets to a geodesic ray $a: [0,\infty) \rightarrow F$ starting orthogonally at $\g_0$ (respectively, $b: [0,\infty) \rightarrow F$ starting orthogonally at $g$).
Since $a_{n_k} \cap b_{n_k} = \emptyset$ for every $k$, we have $a \cap b = \emptyset$.

Hence, $a$ and $b$ belong to two right-angled open-ended geodesic quadrilaterals $Q_i$, $i=1,2$ in $F$,
each having a compact side $u_i$ on $\g_0$ and another compact side $v_i$ on $g$.
Lemma \ref{l:openquadri} implies that $\sinh L(u_i) \sinh L(v_i) \geq 1$.
Corollary \ref{c:B} gives
\begin{equation*}
\sinh \frac{L(u_1)+L(u_2)}2 \,\sinh \frac{L(v_1)+L(v_2)}2
\ge \min \big\{
\sinh L(u_1) \sinh L(v_1),\,
\sinh L(u_2) \sinh L(v_2)
\big\}
\ge 1
\end{equation*}
and so, \eqref{eq:r2} implies
\begin{equation*}
1 \le
\sinh \frac{L(u_1)+L(u_2)}2 \,\sinh \frac{L(v_1)+L(v_2)}2
< \sinh \Big( \frac{1}2\, L(\g_0)\Big) \sinh \Big( \frac{1}2\, L(g)\Big)
< 1,
\end{equation*}
a contradiction.
\end{proof}

%
%%% Figures 6,7 %%%
%
%
%%% end Figures 6,7 %%%
%

\begin{lemma}\label{l:dArch}
Let $F$ be a doubly connected complete Riemannian surface whose boundary is a simple closed geodesic $\g_0$. For given $d>0$ consider the set $Z(d) = \big\{ x\in F \,|\, \dist(x, \g_0) < d  \, \big\}$ and let $\L_d$ be the infinite closed connected component of $F \setminus Z(d)$. Then, for any compact connected component $A$ of $F \setminus Z(d)$ there exists a geodesic arc $g$ in $Z(d)$ of length $L(g) < 2d$ meeting $\g_0$ orthogonally at its endpoints that separates $A$ and $\L_d$.
\end{lemma}

\begin{proof}
Since $A$ has positive distance from $\L_d$ there exists $d' < d$ and a compact connected component $A'$ of $F \setminus Z(d')$ that contains $A$. Only finitely many further compact connected components of $F \setminus Z(d')$, say $A_1, \dots, A_m$ exist that are not completely contained in $Z(d)$ because each such $A_i$ contains some distance ball of radius $d-d'$, the areas of these balls have a lower bound depending only on $d-d'$ and they sum up to less than the area of $F \setminus \L_{d'}$. We can now find  a simple closed curve $\a : [0,1] \to F \setminus (\L_{d'} \cup A_1 \cup \dots \cup A_m)$ that goes once around $A'$ but not around $\L_{d'}$. Since all components of $F \setminus Z(d')$ different from $\L_{d'}, A_1, \dots, A_m$ are ``swept over'' by $Z(d)$, $\a$ is contained in $Z(d)$ and it goes once around $A$ but not around $\L_d$.

Let $h_A$ be a minimal geodesic from $A$ to $\g_0$ and parametrize $\a : [0,1] \to Z(d)$ such that $\a(0) = \a(1) \in h_A$ and such that for some $t_1 \in (0,1]$ the restriction of $\a$ to $[0,t_1]$ goes once around $A$ as illustrated in Figure \ref{fig:dArch} with $\a(t_1) \in h_A$ and such that $\a(t) \notin h_A$ for all $t \in (0,t_1)$, (if $\a$ intersects $h_A$ only once, then of course $t_1 = 1$).

For any $t \in (0,t_1)$ we let $H_t$ be the set of all minimal geodesics from $\a(t)$ to $\g_0$. Then, unlike the dashed arc in Figure \ref{fig:dArch}, no $h \in H_t$ intersects $h_A$ because both, $h$ and $h_A$ are minimal connections to $\g_0$. For $t$ close to $0$ the $h \in H_t$ lie near $h_A$ on the left hand side of $h_A$ (speaking with Figure \ref{fig:dArch}) while for $t$ close to $t_1$ the $h \in H_t$ lie near $h_A$ on the right hand side of $h_A$. By the continuity of the lengths of the $h \in H_t$ as a function of $t$ there exists $\tau \in (0,t_1)$ such that $H_{\tau}$ contains two members $g_1$ and $g_2$ that lie on either side of $h_A$ thereby forming a simple arc $g' =g_1 \cup g_2$ in $Z(d)$ that has both endpoints on $\g_0$ and separates $A$ from $\L_d$. Since $\a(\tau) \in Z(d)$ we have $L(g') < 2d$. On the surface $F \setminus (\L_d \cup A)$ the homotopy class of $g'$ with endpoints moving on $\g_0$ contains a rectifiable curve $g$ of minimal length. Since the distances from $A$ to $\g_0$ and from $\L_d$ to $\g_0$ are equal to $d$ while $L(g) < 2d$, the curve $g$ does not meet $A \cup \L_d$ and is therefore a geodesic meeting $\g_0$ orthogonally at its endpoints. Finally, since $L(g) < 2d$, $g$ is contained in $Z(d)$. \end{proof}

Proposition \ref{p:genfunnel} proves Theorem \ref{t:a} for the case where $S$ is doubly connected.
Our next particular case are the generalized Y-pieces.

\begin{proposition} \label{p:descomposicion3}
Let $Y$ be a generalized Y-piece with curvature $K \ge -1$, boundary geodesics $\{\g_i\}_{i=1}^n$, $n \in \{1,2,3\}$, and $n-3$ generalized cusps.
For each boundary geodesic set
\begin{equation*}
Z_i(d) = \big\{ x\in Y \,|\, \dist(x, \g_i) < d \, \big\} .
\end{equation*}
Then $Z_i(d)$ is topologically equivalent to $\SS^1 \times [0,1)$ for every
\begin{equation*}
0 < d \le d_i :=\arcsinh \cosech \Big( \frac{1}2\, L(\g_i)\Big).
\end{equation*}
Also, the sets $\{Z_i(d_i)\}_{i=1}^n$ are pairwise disjoint.
\end{proposition}

\begin{proof}
If $Y$ is an ordinary Y-piece, then the result is proved in \cite{Bu}.

Hence, we may assume that $Y$ contains one or two generalized cusps.

Let us assume first that $Y$ contains two generalized cusps.
So, we may assume that $\g_1$ is a simple closed geodesic and that $\g_2$ and $\g_3$ are generalized cusps.

In a generalized cusp there is no minimizing simple closed geodesic, but it is possible to have
non minimizing simple closed geodesics generating the fundamental group of the cusp.
If this is the case, we can cut the generalized Y-piece by this geodesic and so,
we are in the case where $Y$ is either a Y-piece or a generalized Y-piece with a single generalized cusp.
Hence, we may assume that there is no simple closed geodesic in each cusp.

%
%%% Figures 8,9 %%%
%
\begin{figure}[b!]
\begin{minipage}[b]{0.48\textwidth}
\centering
\leavevmode
\SetLabels
(.077*.6) $\scriptstyle \alpha_1^2$\\
(.132*.5) $\scriptstyle \alpha_2^2$\\
(.212*.45) $\scriptstyle \alpha_3^2$\\
(.15*.90) $\gamma_1$\\
(.31*.53) $\gamma_2$\\
(.759*.52) $\gamma_3$\\
(.462*.85) $g$\\
(.38*.08) $Y$\\
(.55*.66) $\scriptstyle \alpha_2$\\
(.52*.543) $\scriptstyle \alpha_3$\\
(.515*.45) $\scriptstyle \alpha_1$\\
(.88*.25) $\scriptstyle \alpha_1^3$\\
(.795*.305) $\scriptstyle \alpha_2^3$\\
(.68*.337) $\scriptstyle \alpha_3^3$\\
\endSetLabels
%\ShowGrid
\AffixLabels{%
\includegraphics{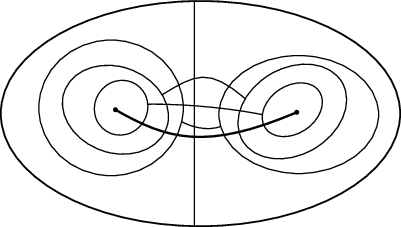}%[width=0.858\linewidth]%
}
\caption{Construction of $\a$
}
\label{fig:figura1}
\end{minipage}
\hspace{0.02\textwidth}
\begin{minipage}[b]{0.48\textwidth}
\centering
\leavevmode
\SetLabels
(.14*.21) $\beta$\\%%%
(.43*.72) $\beta$\\%%%
(.143*.375) $\scriptstyle \beta_3$\\
(.19*.40) $\scriptstyle \beta_2$\\
(.24*.435) $\scriptstyle \beta_1$\\
(.874*.43) $\scriptstyle \alpha_1^3$\\
(.824*.46) $\scriptstyle \alpha_2^3$\\
(.773*.49) $\scriptstyle \alpha_3^3$\\
(.72*.57) $ \sigma$\\
(.831*.82) $\scriptstyle \sigma_1$\\
(.838*.70) $\scriptstyle \sigma_2$\\
(.779*.65) $\scriptstyle \sigma_3$\\
(.17*.92) $\gamma_1$\\
(.268*.53) $\gamma_2$\\
(.74*.43) $\gamma_3$\\
(.50*.93) $a_1$\\
(.09*.695) $a_2$\\
(.35*.445) $\alpha$\\
(.33*.81) $G$\\
(.55*.44) $\scriptstyle x_1$\\
(.605*.438) $\scriptstyle x_2$\\
(.650*.436) $\scriptstyle x_3$\\
\endSetLabels
%\ShowGrid
\AffixLabels{%
\includegraphics{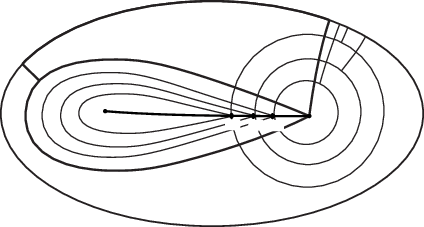}%[width=0.908\linewidth]
}
\caption{Construction of $\b$ and $\s$ }
\label{fig:figura2}
\end{minipage}
%\caption{Construction of $\overline{G}$}
\end{figure}
%
%%% end Figures 8,9 %%%
%

Consider a minimizing geodesic $g$ from $\g_1$ to $\g_1$ and separating $\g_2$ and $\g_3$ in $Y$
(see Figure \ref{fig:figura1}).
For each $i =2,3,$
consider a sequence of pairwise disjoint simple closed curves $\{\a_n^i\}_{n=1}^\infty \subset Y \setminus g$
such that each of them generates the fundamental group of the cusp $\g_i$ and $\lim_{n \to \infty} \dist(\a_n^i,\a_1^i) = \infty$.
For each $n$, there exists a minimizing geodesic $\a_n$ from $\a_n^2$ to $\a_n^3$.
Note that $g \cap \a_n \neq \emptyset$ for every $n$,
since $g$ and separates $\g_2$ and $\g_3$.
Since these arcs have minimal length, the intersection $g \cap \a_n$ is a single point.
We may assume that $\a_n:[-u_n,v_n] \rightarrow Y$ is parametrized by arc-length and $\a_n(0) = g \cap \a_n$ for every $n$.
Since $g$ and $\SS^1$ are compact sets, there exists a subsequence $\{n_k\}_{k=1}^\infty$ of $\NN$ such that
$\{\a_{n_k}(0)\}_{k=1}^\infty$ and the tangent vectors $\{\bdot{\a}_{n_k}(0)\}_{k=1}^\infty$ are convergent.
By using the Arzel\`a-Ascoli theorem, one can check that $\{\a_{n_k}\}_{k=1}^\infty$ converges (uniformly on compact sets)
to a geodesic line $\a:\RR \rightarrow Y$ from $\g_2$ to $\g_3$ such that $\a(0) \in g$ and $\a$ is minimizing
(see Figure~\ref{fig:figura1}).

For each $n$, choose $x_n \in \a \cap \a_n^3$,
a minimizing geodesic loop $\b_n$ with base point at $x_n$ and ``surrounding" $\g_2$,
and a minimizing geodesic $\s_n:[0,w_n] \rightarrow Y$ from $\g_1$ to $\a_n^3.$
By the minimality of $\b_n$ and $\s_n$, they do not intersect except perhaps at $x_n$.
By the minimality of $\a$ and $\b_n$, we have $\a \cap \b_n = \{x_n\}$.
Thus, $x_n$ is surrounded by $\b_k$ for every $k>n$ and so, the minimality of $\b_n$ gives that
$\b_n$ is surrounded by $\b_k$ for every $k>n$.
Hence, the domains in $Y$ bounded by $\b_n$ increase with $n$ and so, $\{\b_{n}\}_{n=1}^\infty$
(with appropriate arc-length parametrizations)
converges (uniformly on compact sets)
to a geodesic line $\b:\RR \rightarrow Y$ from $\g_3$ to $\g_3$ and surrounding $\g_2$.
Since $\g_1$ is a compact set, there exists a convergent subsequence $\{\s_{n_k}(0)\}_{k=1}^\infty$;
since $\s_{n_k}'(0)$ is orthogonal to $\g_1$ for every $k$, $\{\s_{n_k}\}_{k=1}^\infty$ converges (uniformly on compact sets)
to a geodesic ray $\s: [0,\infty) \rightarrow Y$ from $\g_1$ to $\g_3$.
Since $\b_{n_k} \cap \s_{n_k} \subseteq \{x_{n_k}\}$ for every $k$, we have $\b \cap \s = \emptyset$
(see Figure \ref{fig:figura2}).

Let $a_2$ be a minimizing geodesic from $\g_1$ to $\b$.
The points $a_2 \cap \g_1$ and $\s \cap \g_1$ divide $\g_1$ into two geodesics, denote by $a_1$ one of these geodesics satisfying $L(a_1) \le L(\g_1)/2$.
Let $G$ be the right-angled open-ended geodesic quadrilateral contained in $Y$ with finite sides $a_1$, $a_2$, and infinite sides $\s$ and a subset of $\b$.
Lemma \ref{l:openquadri} gives
\begin{equation*}
\sinh L(a_1) \sinh L(a_2) \ge 1.
\end{equation*}
Consequently,
\begin{equation*}
\begin{aligned}
\sinh \Big( \frac{1}2\, L(\g_1) \Big) \sinh L(a_2)
& \ge \sinh L(a_1) \sinh L(a_2)
\ge 1,
\\
L(a_2)
& \ge \arcsinh \cosech \Big( \frac{1}2\, L(\g_1) \Big) = d_1.
\end{aligned}
\end{equation*}
Hence, we have proved that $Z_1(d_1)$ is contained in a doubly connected domain $W$ whose boundary is $\g_1 \cup \b$.
It remains to prove that $W \setminus Z_1(d)$ is connected for all $d \leq d_1$.

For this we proceed as in the proof of Proposition \ref{p:genfunnel}. Assume that $W \setminus Z_1(d)$ is disconnected for some $d \leq d_1$. By Lemma \ref{l:dArch} there exists a geodesic arc $g$  in $W$ that meets $\g_1$ orthogonally at its endpoints whose length satisfies $L(g) < 2d \leq 2d_1$ and thus
\begin{equation}\label{eq:r3}
\sinh \Big( \frac{1}2\, L(\g_1)\Big) \sinh \Big( \frac{1}2\, L(g)\Big)
< 1.
\end{equation}

%\begin{figure}[h!]
%\begin{minipage}{0.45\textwidth}
%\centering
%\includegraphics[width=\linewidth]{imgs/Dibujo5.png}
%\caption{Construction of $a^2$, $a^3$, $b^2$ and $b^3$}
%\label{fig:figura5}
%\end{minipage}
%\end{figure}

%
%%% Figures 10, 11 %%%
%
\begin{figure}[b!]
\begin{minipage}[b]{0.48\textwidth}
\centering
\leavevmode
\SetLabels
(.046*.5) $\scriptstyle \alpha_1^2$\\
(.105*.5) $\scriptstyle \alpha_2^2$\\
(.175*.5) $\scriptstyle \alpha_3^2$\\
(.964*.5) $\scriptstyle \alpha_1^3$\\
(.904*.5) $\scriptstyle \alpha_2^3$\\
(.84*.5) $\scriptstyle \alpha_3^3$\\
(.382*.14) $\scriptstyle a_1^2$\\
(.307*.15) $\scriptstyle a_2^2$\\
(.252*.17) $\scriptstyle a_3^2$\\
(.889*.295) $\scriptstyle a_1^3$\\
(.83*.24) $\scriptstyle a_2^3$\\
(.782*.194) $\scriptstyle a_3^3$\\
(.34*.895) $\scriptstyle b_1^2$\\
(.376*.84) $\scriptstyle b_2^2$\\
(.42*.76) $\scriptstyle b_3^2$\\
(.786*.815) $\scriptstyle b_1^3$\\
(.692*.81) $\scriptstyle b_2^3$\\
(.632*.76) $\scriptstyle b_3^3$\\
(.19*.18) $a^2$\\
(.66*.18) $a^3$\\
(.49*.520) $b^2$\\
(.561*.62) $b^3$\\
(.15*.90) $\gamma_1$\\
(.253*.565) $\gamma_2$\\
(.778*.568) $\gamma_3$\\
(.495*.84) $g$\\
\endSetLabels
%\ShowGrid
\AffixLabels{%
\includegraphics{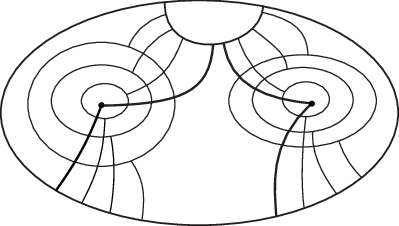}%[width=0.856\linewidth]
}
\caption{Construction of $a^2$, $a^3$, $b^2$ and $b^3$}
\label{fig:figura5}
\end{minipage}
\hspace{0.015\textwidth}
\begin{minipage}[b]{0.48\textwidth}
\centering
\leavevmode
\SetLabels
(.06*.79) $\gamma_1$\\
(.299*.53) $\gamma_2$\\
(.744*.427) $\gamma_3$\\
(.50*.92) $a_1$\\
(.192*.794) $a_2$\\
(.53*.44) $\alpha$\\
(.42*.17) $\beta$\\
(.764*.29) $\scriptstyle \alpha_n^3$\\
\endSetLabels
%\ShowGrid
\AffixLabels{%
\includegraphics{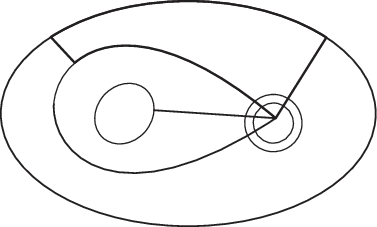}%[width=0.806\linewidth]
}
\caption{Geodesic $\b$ separating $Z_1(d_1)$ and $Z_2(d_2)$}
\label{fig:uniquecusp}
\end{minipage}
%\caption{Construction of $\overline{G}$}
\end{figure}
%
%%% end Figures 10, 11 %%%
%

By Lemma \ref{l:broken},
for each $n$ and $i=2,3,$ we can consider minimizing geodesics $a_n^i$ from $\g_1$ to $\a_n^i$, and $b_n^i$ from $g$ to $\a_n^i$
(see Figure \ref{fig:figura5}).
By the minimality of $a_n^2,a_n^3,b_n^2,b_n^3,$ they are pairwise disjoint for each fixed $n$ (Lemma \ref{l:nocut}, as always).
Since $\g_1$ is a compact set, one can check, as in the previous arguments, that there exists a subsequence $\{n_k\}_{k=1}^\infty$
such that $\{a_{n_k}^i\}_{k=1}^\infty$ (respectively, $\{b_{n_k}^i\}_{k=1}^\infty$) converges uniformly on compact sets to a geodesic ray $a^i$ starting orthogonally at $\g_1$ (respectively, $b^i$ starting orthogonally at $g$), for $i=2,3$.
Since
$a_{n_k}^2,a_{n_k}^3,b_{n_k}^2,b_{n_k}^3$ are pairwise disjoint for each fixed $k$, we conclude that $a^2,a^3,b^2,b^3$ are pairwise disjoint.

Hence, $a^i$ and $b^i$ belong to a right-angled open-ended geodesic quadrilateral $Q_i$, $i=2,3$ in $Y$,
each having a compact side $u_i$ on $\g_1$ and the other compact side $v_i$ on $g$.
Lemma \ref{l:openquadri} implies that $\sinh L(u_i) \sinh L(v_i) \geq 1$.
As in the proof of Proposition \ref{p:genfunnel}, Corollary \ref{c:B} gives
\begin{equation*}
\sinh \frac{L(u_2)+L(u_3)}2 \,\sinh \frac{L(v_2)+L(v_3)}2
\ge \min \big\{
\sinh L(u_2) \sinh L(v_2),\,
\sinh L(u_3) \sinh L(v_3)
\big\}
\ge 1
\end{equation*}
and so, \eqref{eq:r3} implies
\begin{equation*}
1 \le
\sinh \frac{L(u_2)+L(u_3)}2 \,\sinh \frac{L(v_2)+L(v_3)}2
< \sinh \Big( \frac{1}2\, L(\g_1)\Big) \sinh \Big( \frac{1}2\, L(g)\Big)
< 1,
\end{equation*}
a contradiction.

\medskip

Finally, let us assume that $Y$ has a unique generalized cusp.
So, we may assume that $\g_1$ and $\g_2$ are simple closed geodesics and that $\g_3$ is a generalized cusp.
A similar argument to the previous one gives the desired results. The only difference in the construction of $\b$ which now surrounds the geodesic $\g_2$ (Figure \ref{fig:uniquecusp}) is that we consider in this case a minimizing geodesic $\a_n$ from $\g_2$ to $\a_n^3$
(instead of from $\a_n^2$ to $\a_n^3$). As before, $Z_1(d_1)$ is doubly connected and is contained in the doubly connected domain with boundaries $\g_1$, $\b$. By symmetry, $Z_2(d_2)$ is contained in the domain with boundaries $\g_2$, $\b$. Hence, the collars $Z_1(d_1)$ and $Z_2(d_2)$ are separated by $\b$ and therefore disjoint.
%We just mention two different details in the argument:
%In this case it is necessary to check that the collars $Z_1$ of $\g_1$ and $Z_2$ of $\g_2$ are disjoint;
%this is obtained following the argument in \cite[pp.353-354]{Bu}.
\end{proof}

We finish the proofs of Theorems \ref{t:a} and \ref{t:b}. Let $S$ (without boundary) and $\g$ be as in the theorems.
We first deal with the cases where Theorem \ref{t:b} is void.
If $S$ is a cylinder then, as already mentioned, Theorem \ref{t:a} follows from Proposition \ref{p:genfunnel}.
If $S$ is a sphere with three ends and exactly one of them has a collared neighborhood with geodesic boundary, then $\g$ separates $S$ into such a neighborhood and a generalized Y-piece with two cusps and Theorem \ref{t:a} follows from Propositions \ref{p:genfunnel} and \ref{p:descomposicion3}. If $S$ is a sphere with four ends and none of them has a collared neighborhood with geodesic boundary, then $\g$ separates $S$ into two generalized Y-pieces and Theorem \ref{t:a} follows from Proposition \ref{p:descomposicion3}. In these cases Theorem \ref{t:b} is void because any other simple closed geodesic is either homotopic to or intersecting $\g$.

Assume now that $S$ is not in one of these cases. Then, since furthermore $S$ is by hypothesis not a closed torus, Theorem \ref{t:descomposicion2} implies that that there exist homotopically nontrivial simple closed geodesics in $S$ that are not homotopic to and not intersecting $\g$, and we let $\mu$ be one of them, as in the hypothesis of Theorem \ref{t:b}.

Cut $S$ along $\g$ and $\mu$.
Thus, $S$ splits into either one, two or three connected components $\{S_i\}_{i=1}^m$ ($m \in \{1,2,3\}$).
We have to estimate the widths of the collars at $\g$ and $\mu$ on the boundaries of these components.
Fix some $S_i$. Since $\g$ and $\mu$ are homotopically nontrivial simple closed geodesics in $S$, one gets that $S_i$ is not simply connected.

Assume first that $S_i$ is doubly connected.
Since $\g$ and $\mu$ are not homotopic in $S$, it is not possible to have $\p S_i = \g \cup \mu$.
By symmetry we may assume that $\p S_i = \g$.
Then, $S_i$ is a generalized funnel or a generalized cusp with geodesic boundary and the estimate is given by Proposition \ref{p:genfunnel}.

Assume now that $S_i$ is neither simply nor doubly connected.
Thus, Theorem \ref{t:descomposicion2} gives that $S_i$
is the union (with pairwise disjoint interiors) of generalized Y-pieces, generalized funnels and halfplanes.
Furthermore, if $\g$ (respectively, $\mu$) is contained in $\p S_i$,
then $\g$ (respectively, $\mu$) is contained in the border of one of these generalized Y-pieces. The estimates of the width and the lower bound for the distance of $\g$ and $\mu$ are thus given by Proposition \ref{p:descomposicion3}. This completes the proofs of Theorems \ref{t:a} and \ref{t:b}.  \hfill \qedsymbol

\section{Proof of Theorem \ref{t:d}}\label{Proofd}

Assume first that $S$ is simply connected.
Since $S$ is not homeomorphic to a sphere, it is homeomorphic to $\RR^2$.
Hence, $S\setminus D_\g$ is a generalized funnel or a generalized cusp with geodesic boundary, and Proposition \ref{p:genfunnel} gives the result.

Assume now that $S$ is not simply connected.
Then, $S\setminus D_\g$ is a complete orientable geodesically bordered Riemannian surface which is neither simply nor doubly connected.
By Theorem \ref{t:descomposicion2}, $S\setminus D_\g$
is the union of generalized Y-pieces, generalized funnels and halfplanes;
also, $\g$ is contained in the border of one of these generalized Y-pieces.
Therefore, Proposition \ref{p:descomposicion3} implies the result.

\section{Improvements of Theorem \ref{t:c}}
\label{sec:Improvements}

In this section we improve Theorem \ref{t:c} in several ways.

\smallskip

Theorem \ref{t:c} follows from Theorem \ref{t:a}, since $\eta$ goes through the collar of $\g$ at least once.
If $S \setminus \g$ is not connected, then $\eta$ goes through the collar of $\g$ at least twice;
this fact gives the following result:

\begin{theorem} \label{t:c2}
Let $S$ be a complete orientable Riemannian surface with curvature $K \ge -k^2$.
If $\g$ and $\eta$ are homotopically nontrivial simple closed geodesics in $S$,
$S \setminus \g$ is not connected,
and $\eta$ is not homotopic to but intersecting $\g$, then
\begin{equation*}
\sinh \left( \frac{k}2\, L(\g) \right) \sinh \left( \frac{k}4\, L(\eta) \right)
\ge 1 .
\end{equation*}
\end{theorem}

Let us consider genus zero surfaces now.
This type of surface plays an important role in Geometric Function Theory, since the Poincar\'e metric of any domain in the complex plane (with more than one boundary point) is complete and has constant curvature $-1$,
and every holomorphic function between two of these domains is $1$-Lipschitz with respect to their respective Poincar\'e metrics.

\smallskip

Since $S \setminus g$ is not connected for any simple closed curve $g$ in any genus zero surface $S$,
Theorem \ref{t:c2} and the symmetry between the two geodesics give the following result:

\begin{theorem} \label{t:c3}
Let $S$ be a genus zero complete orientable Riemannian surface with curvature $K \ge -k^2$.
If $\g$ and $\eta$ are homotopically nontrivial simple closed geodesics in $S$,
and $\eta$ is not homotopic to but intersecting $\g$, then
\begin{equation*}
\min \left\{
\sinh \left( \frac{k}2\, L(\g) \right) \sinh \left( \frac{k}4\, L(\eta) \right),\,
\sinh \left( \frac{k}4\, L(\g) \right) \sinh \left( \frac{k}2\, L(\eta) \right)
\right\}
\ge 1 .
\end{equation*}
\end{theorem}

If we choose
\begin{equation*}
u_1 = \frac{k}2\, L(\g),
\quad
v_1 = \frac{k}4\, L(\eta),
\quad
u_2 = \frac{k}4\, L(\g),
\quad
v_2 = \frac{k}2\, L(\eta),
\end{equation*}
then Corollary \ref{c:B} and Theorem \ref{t:c3} have the following consequence:

\begin{theorem} \label{t:c4}
Let $S$ be a genus zero complete orientable Riemannian surface with curvature $K \ge -k^2$.
If $\g$ and $\eta$ are homotopically nontrivial simple closed geodesics in $S$,
and $\eta$ is not homotopic to but intersecting $\g$, then
\begin{equation*}
\sinh \left( \frac{3k}8\, L(\g) \right) \sinh \left( \frac{3k}8\, L(\eta) \right)
> 1 .
\end{equation*}
\end{theorem}

\section{Thin collared ends}
\label{sec:Thin}

In this section $C$ is a complete doubly connected Riemannian surface of curvature $K \geq - k^2$ ($k>0$) whose boundary is a simple closed geodesic $\g_0$. Thus, $C$ is either a generalized funnel or a generalized cusp with a geodesic boundary. The presence of positive curvature allows $C$ to have bifurcating arms like in the illustration of Figure \ref{fig:cactus1} that locally give an impression of bigger topology. This is certainly a hazy remark, but we shall catch an element of it that allows a precise definition. As a result we get a collar theorem for thin ends.

\begin{definition}
A point $p \in C$ is called \emph{of type 1} if there exists $\d_p > 0$ such that any path $\s$ that goes from $p$ to infinity satisfies
\begin{equation*}
\dist(\s, \g_0) < \dist(p,\g_0) - \d_p.
\end{equation*}
Points in $C$ that are not of type 1 are called \emph{of type 0}. We let $C_1, C_0$ be the subsets of points of type 1 and type 0, respectively. It follows from the triangle inequality that $C_1$ is an open subset of $C$. Any open connected component of $C_1$ is called a \emph{cactus arm}.
\end{definition}

The botanical nomenclature ``cactus arm'' is inspired by Figure \ref{fig:cactus1} and is not to be taken too seriously. In this figure, $p$ and $q$ are of type 1 while $r$ and $s$ are of type 0.

There is also another way to define $C_0$: For any $d \geq 0$ let $\omega_d = \{x \in C \mid \dist(x, \g_0) = d\}$, define $\O_d^0$ to be the infinite open connected component of $C \setminus \omega_d$ and let $\omega_d^0 = \partial \O_d^0$ be the boundary of $\O_d^0$. Then $\omega_d^0$ is a subset of $\omega_d$ (not necessarily a full closed connected component of $\omega_d$)  and is a simple closed Lipschitz curve homotopic in $C$ to $\g_0$ (Theorem \ref{t:rectA}).

We now have
\begin{equation}\label{eq:C0}
C_0 = \overline{\bigcup_{d\geq 0} \omega_d^0} \, ,
\end{equation}
where the overline denotes the closure of a set.

\begin{proof}
Let $C_0'$ be the right hand side in \eqref{eq:C0}. For any $d \geq 0$ and any $p \in \omega_d^0$ there exists a path $\s$ from $p$ to infinity that, except for its initial point, is contained in $\O_d^0$. Hence, any such $p$ and thus all $\omega_d^0$ are contained in $C_0$. Since $C_0$ is closed it follows that $C_0' \subseteq C_0$.

Conversely, consider $p \in C_0$.
Let $d=\dist(p,\g_0)$ and take $\e > 0$.
Since $p$ is not of type 1 there exists a path $\s$ from $p$ to infinity such that all of its points are at distance $> d-\e$ from $\g_0$.
This path thus lies in $\O_{d - \e}^0$ and hence, $p \in \O_{d - \e}^0$. Let now $\eta$ be a shortest connection from $p$ to $\g_0$. It has length $d$ and splits into a part $\eta'$ from $p$ to $\omega_{d-\e}^0$ and a part $\eta''$ from there to $\g_0$.
The latter has length $\geq d - \e$ (actually $= d - \e$, by minimality) and thus $L(\eta') \leq \e$.
It follows that $\dist(p, \omega_{d-\e}^0) \leq \e$ and, hence, $\dist(p,C_0') \leq \e$. As this holds for any $\e >0$ and $C_0'$ is closed we get $p \in C_0'$.
This concludes the proof of \eqref{eq:C0}.
\end{proof}

\begin{remark}\label{r:cactusarm} Each cactus arm $A$ is one of the connected components of $C \setminus \overline{Z(d_A)}$, where $d_A = \inf \{\dist(x,\g_0) \mid x \in A \}$. Its boundary $\p A$ is a simple closed Lipschitz curve of constant distance to $\g_0$.
\end{remark}
\begin{proof}
Select $p \in A$ and let $A'$ be the connected component of $C \setminus{\overline{Z(d_A)}}$ that contains $p$. We show that $A = A'$. If $x \in A$, then there exists a path $\tau$ in $A$ from $p$ to $x$. By the definition of $d_A$ all points on $\tau$ have distance $> d_A$ to $\g_0$ which implies that $\tau$ does not intersect $\p A'$ and so $x \in A'$. This proves $A \subseteq A'$.

Before we proceed to the inverse inclusion we first show that $A'$ is not the infinite component $\O_{d_A}^0$. Assume the contrary. Then there exists a path $\zeta$ in $\O_{d_A}^0$ going from $p$ to infinity and we set $\e = \frac{1}{2}\dist(\zeta, \omega_{d_A}^0)$, $d_A' = d_A + \e$. Then, like $\omega_{d_A}^0$, the curve $\omega_{d_A'}^0$ is simple closed and homotopic to $\g_0$. Furthermore, it runs within distance $\e$ from $\omega_{d_A}^0$ and does therefore not intersect $\zeta$. Hence, $\zeta$ is contained in $\O_{d'_A}^0$. Since $p$ is a point of $\zeta$ we thus have $p \in A\cap \O_{d'_A}^0$. It follows that on either side of $\omega_{d'_A}^0$ there are points of $A$. Since $A$ is connected, $\omega_{d'_A}^0$ intersects $A$. But all points on $\omega_{d'_A}^0$ are of type $0$, whereas the points of $A$ are of type 1, a contradiction.

We now show that $A' \subseteq A$. Take $x \in A'$. Since $A'$ is not $\O_{d_A}^0$, any path $\s$ from $x$ to infinity must cross $\p A'$, and so $x$ is of type 1 with $\d_x = \dist(x, \p A')$. Hence, $A' \subset C_1$. Since $A'$ is connected it is contained in one of the connected components of $C_1$, and since $p \in A'$ this component is $A$. Thus, $A' \subseteq A$.

The statement about $\p A$ holds by Theorem \ref{t:rectA} because we now know that $A=A'$.
\end{proof}

There exists at least one minimising geodesic ray emanating orthogonally from $\g_0$ going to infinity and we let $h_0$ be one of them which we fix in the following. Since $h_0$ is minimising, all points on $h_0$ are of type 0. Furthermore, any $\O_d^0$ intersects $h_0$ and hence also all $\omega_d^0$ intersect $h_0$. This shows that the union of all $\omega_d^0$ and, hence, $C_0$ \emph{is connected}. We call $C_0$ the \emph{trunk} of $C$.

We now show that arms can bifurcate from the trunk only at places where it is thick.

\begin{proposition}\label{p:cactusarm}
Let $A$ be a cactus arm of $C$ and $p \in A$. Then any closed curve $\g$ passing through $p$ that is freely homotopic to $\g_0$ has length $L(\g) \geq 2\arcsinh(1)/k$.
\end{proposition}

\begin{proof}
Again, without loss of generality we assume $k=1$. By Remark \ref{r:cactusarm} the boundary $\p A$  is a simple closed curve.
We denote by $d_A$ the common distance of the points of $\p A$ to $\g_0$.
By Lemma \ref{l:dArch} there exists a geodesic arc $g$ of length $L(g) \leq 2d_A$ with both endpoints on $\g_0$ and meeting $\g_0$ there orthogonally that splits $C$ into a doubly connected domain $U$ that contains $h_0$ and a simply connected domain $V$ that contains $A$.

Let now $m \in g$ be the midpoint of $g$ and $\mu$ a minimising geodesic ray in $U$ going from $m$ to infinity. Then $U$ is further split by $h_0$ and $\mu$ into two open ended quadrangles $U_1$, $U_2$. For $i=1,2$ we select in $U_i$ a minimal geodesic arc $a_i$ from $m$ to $h_0$. These arcs meet $h_0$ orthogonally on opposite sides. The strategy of the proof is to show first that $L(\g) \geq L(a_1) + L(a_2)$, and in a second step that $L(a_1) + L(a_2)$ is bounded from below by the constant in the proposition.

We may assume that $\g$ is of minimal length and thus is a geodesic loop with base point $p$ that intersects $h_0$ exactly once. It splits into a part $\g'$ from $p$ to $h_0$, say arriving at $h_0$ on the side of $U_1$ and a part $\g''$ leaving $h_0$ on the side of $U_2$ returning to $p$. In the homotopy class of $\g'$ where the endpoint $p$ is fixed and the other endpoint is allowed to move freely on $h_0$ there exists an arc $c_1$ of minimal length which is thus a minimising geodesic from $p$ to $h_0$ arriving orthogonally at $h_0$ on the side of $U_1$. There is an analogous minimising arc $c_2$ arriving orthogonally at $h_0$ on the side of $U_2$ and we have $L(\g) = L(\g') + L(\g'') \geq L(c_1) + L(c_2)$. For the first part of the proof it remains to show that
\begin{equation}\label{eq:Lca}
L(c_i) \geq L(a_i), \quad i=1,2.
\end{equation}
We carry it out for $i=1$. Let $x$ be the intersection point of $c_1$ and $g$. Then $c_1$ is split into a segment $v_1$ from $p$ to $x$ and a segment $u_1$ from $x$ to $h_0$. Also, the half of $g$ that contains $x$ (a priori, both halves are possible and Figure \ref{fig:cactus2} illustrates the two cases by the two dashed curves) is split into an arc $w$ from $m$ to $x$ and an arc $w'$ from $x$ to $\g_0$. Since $L(v_1) + L(w') \geq d_A$ while $L(w) + L(w') = \frac{1}{2} L(g) \leq d_A$ we have $L(v_1) \geq L(w)$ and so $L(c_1) = L(v_1) + L(u_1) \geq L(w) + L(u_1)$. Since the path formed by $w$ from $m$ to $x$ and $u_1$ from $x$ to $h_0$ is in the homotopy class of $a_1$ we have $L(w) + L(u_1) \geq L(a_1)$. This proves \eqref{eq:Lca} for $i=1$ and for $i=2$ the proof is the same. Hence, $L(\g) \geq L(a_1) + L(a_2)$.

%
%%% Figures 12, 13 %%%
%
\begin{figure}[t!]
\begin{minipage}[b]{0.48\textwidth}
\centering
\leavevmode
\SetLabels
(.75*.0) $\gamma_0$\\
(.063*.72) $q$\\
(.418*.475) $s$\\
(.744*.435) $p$\\
(.566*.065) $r$\\
(.59*.26) $\sigma$\\
\endSetLabels
%\ShowGrid
\AffixLabels{%
\includegraphics{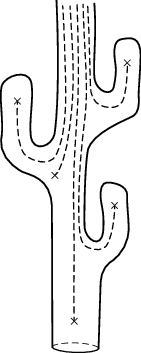}
}
\caption{Collared end with arms.}
\label{fig:cactus1}
\end{minipage}
\hspace{-0.02\textwidth}
\begin{minipage}[b]{0.48\textwidth}
\centering
\mbox{}\vspace{0.5cm}\\
\leavevmode
\SetLabels
(.93*.16) $\gamma_0$\\
(.57*.10) $U_1$\\
(.50*.85) $U_2$\\
(.85*.70) $V$\\
(.53*.33) $a_1$\\
(.40*.755) $a_2$\\
(.84*.575) $p$\\
(.30*.125) $u_1$\\
(.40*.39) $u_1$\\
(.82*.35) $v_1$\\
(.085*.555) $h_0$\\
(.35*.54) $\mu$\\
(.73*.22) $g$\\
(.73*.80) $g$\\
(.73*.50) $m$\\
(.73*.38) $x$\\
(.73*.65) $x$\\
(.645*.42) ${\scriptstyle \varphi_1'}$\\
(.575*.466) $\scriptstyle \varphi_1''$\\
(.645*.60) $\scriptstyle \varphi_2'$\\%%
(.575*.54) $\scriptstyle \varphi_2''$\\
\endSetLabels
%\ShowGrid
\AffixLabels{%
\includegraphics{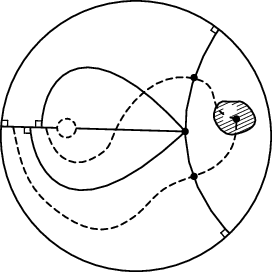}%Dibujo2new
}
\mbox{}\vspace{0.5cm}\\
\caption{Proof of Proposition \ref{p:cactusarm}}
\label{fig:cactus2}
\end{minipage}
%\caption{Construction of $\overline{G}$}
\end{figure}
%
%%% end Figures 12, 13 %%%
%

For the second part of the proof we note that for $i = 1,2$, the orthogonal $a_i$ from $m$ to $h_0$ divides $U_i$ into a compact quadrangle $U_i'$ with three right angles and some angle $\varphi_i'$ at $m$, and an open ended triangle $U_i''$ with a right angle at the vertex on $h_0$ and some angle $\varphi_i''$ at $m$. We let $\varphi_i$ be the minimum of $\varphi_i'$ and $\varphi_i''$. Since the sum of all angles at $m$ equals $\pi$ we get $\varphi_1+\varphi_2 \leq \pi/2$.

By Corollaries \ref{c:Dopen} and \ref{c:Dclosed} (and \eqref{eq:trig2}) we have $L(a_i) \geq \arcsinh( \cot \varphi_i)$. Since the function involved here is convex and monotone decreasing in the interval $[0,\pi/2]$  it follows that
\begin{equation}\label{eq:Laa}
L(a_1) +  L(a_2)\geq 2\arcsinh\Big( \cot \Big(\frac{\varphi_1 + \varphi_2}{2}\Big)\Big) \geq 2\arcsinh\Big( \!\cot \frac{\pi}{4}\Big) = 2\arcsinh(1).
\end{equation}
With \eqref{eq:Lca} an \eqref{eq:Laa} the proof of the proposition is complete.
\end{proof}

To formulate the next theorem we make a definition.

\begin{definition}
$C$ is called \emph{$\l$-thin} (for given $\l >0$) if for any $d \geq 0$ the boundary $\omega_d^0$ of the infinite connected component of $C \setminus \overline{Z(d)}$ has length $L(\omega_d^0) \leq \l$.
\end{definition}
Note that in this definition there is no condition on the thicknesses or the lengths of the cactus arms. The announced collar theorem is this:
\begin{theorem} \label{t:cactus1}
If $C$ is $\l$-thin with $\l < 2\arcsinh(1)/k$, then for any $d > 0$ the boundary part $\p Z(d) \setminus \g_0$ is a simple closed Lipschitz curve homotopic to $\g_0$. In particular, $Z(d)$ is doubly connected for any $d > 0$.
\end{theorem}
\begin{proof} We first show that $C$ has no cactus arms. Assume the contrary. Then $\p C_1 \neq \emptyset$. Take $p \in \p C_1$. By Proposition \ref{p:cactusarm} any simple closed curve $\g$ freely homotopic to $\g_0$ passing through $p$ has length $L(\g) \geq 2 \arcsinh(1)/k$. But we also have $p \in \p C_0$, and so, by hypothesis, one can construct such a $\g$ out of a nearby $\omega_d^0$ that has length $L(\g) <  2 \arcsinh(1)/k$, a contradiction.

Assume now that for some $d > 0$ there exists a point $p \in C$ with $\dist(p, \g_0) = d$, that does not belong to $\omega_d^0$. We shall lead this to a contradiction by modifying the Riemannian metric arbitrarily little in a small neighborhood of $p$ so as to create a cactus arm. To this end we let $B_p(\rho) = \{x \in C \mid \dist(x,p) < \rho \}$ with $0 < \rho < \frac{1}{4}\dist(p, \omega_d^0)$ be a small coordinate neighborhood of $p$. We furthermore take $\rho$ so small that $\l + 2\rho < 2 \arcsinh(1)/k$. Since $C$ has no cactus arms we have $C = C_0$ and so, by \eqref{eq:C0}, the curves $\omega_{d'}^0$, $d'>0$,  are dense in $C$. We can thus, furthermore, adjust $\rho$ such that there exists $\hat{d}$ with the property that $\dist(p,\omega_{\hat{d}}^0) = \rho$.

Now we take a smooth function $\psi : B_p(\rho) \to [0,1]$ that satisfies $\psi(p) = 1$ and $\psi(x) = 0$ for $x \in B_p(\rho) \setminus B_p(\rho/2)$. For small $t>0$ we let $ds'^2$ be the Riemannian metric on $C$ which coincides with the given metric $ds^2$ on $C \setminus B_p(\rho/2)$ and is defined by $ds'^2 = (1+t \psi(x))^2 ds^2$ for $x \in B_p(\rho)$. Note that $\omega_d^0$ is not altered when we replace $ds^2$ by $ds'^2$.

As $t \to 0$ the curvature of $ds'^2$ converges to the curvature of $ds^2$. We can thus take $t>0$ so small that the curvature of $ds'^2$ has a lower bound $- k'^2$ with $k'$ so close to $k$ that still $\l+2\rho < 2\arcsinh(1)/{k'}$.

With respect to $ds'^2$ point $p$ sits on a cactus arm, because its distance to $\g_0$ is bigger than $d$ and any curve from $p$ to infinity must cross $\omega_d^0$ which implies that $p$ is of type 1. By Proposition \ref{p:cactusarm} applied to $C$ endowed with $ds'^2$ any simple closed curve $\g$ freely homotopic to $\g_0$ passing through $p$ has length $L(\g) \geq 2 \arcsinh(1)/k'$. On the other hand the curve $\hat{\g}$ that goes from $p$ along a minimal connecting  arc to $\omega_{\hat{d}}^0$ then once around $\omega_{\hat{d}}^0$ and back to $p$ has length $L(\hat{\g}) = 2\rho + L(\omega_{\tilde{d}}^0) \leq 2\rho + \l < 2\arcsinh(1)/{k'}$, a contradiction.
\end{proof}

\section{Appendix: Rectifiability}
\label{sec:Rect}

There are many results in the literature about the rectifiability of distance sets in Riemannian manifolds and more general metric spaces. In the previous section, we made significant use of the following variant (Theorem \ref{t:rectA}) for which we include a proof because we could not find one in the literature. A brief comparison with other variants will be given in Remark \ref{r:biblio}.

In the theorem, $M$ has a boundary curve. By this we mean that $M$ is a subsurface of a larger (not necessarily complete) Riemannian surface $M'$ that contains $M$ together with the boundary in its interior.

\begin{theorem} \label{t:rectA}
Let $M$ be an orientable complete bordered Riemannian surface of genus zero whose boundary is a simple closed curve $\sC$. For $d > 0$ consider the set $Z(d) = \{x \in M \mid \dist(x,\sC) < d\}$. Then, for any open connected component $\O$ of $M \setminus \overline{Z(d)}$ the boundary $\omega = \p\O$ is a simple closed Lipschitz curve.
\end{theorem}

The Lipschitz property will come out naturally in the following form:

\begin{definition}\label{d:Lipcurve}
A simple curve $c : [a,b] \to M$ shall be said to be \emph{Lipschitz}, more precisely, a \emph{one dimensional Lipschitz submanifold} of $M$ if there exists a decomposition $a = a_0 < a_1 < a_2 < \dots < a_n = b$, and for each $i = 1, \dots, n$, a coordinate neighborhood $V_i$ of $M$ such that the restriction of $c$ to $[a_{i-1},a_i]$ is contained in $V_i$ and admits a change of parameter so that with respect to the coordinates of $V_i$ it is of the form $t \mapsto (t, f_i(t))$, $t \in [a_{i-1},a_i]$, where $f_i : [a_{i-1},a_i] \to \RR$ is a Lipschitz function.
\end{definition}

Lipschitz curves are, of course, rectifiable. There are many variant concepts of Lipschitz manifolds and Lipschitz submanifolds, we refer the reader to \cite{P} for a thorough discussion.

We shall actually prove a stronger version of Theorem \ref{t:rectA} which we formulate separately so as to keep the theorem shaped to the needs of the preceding section.
In this version, $M$ and $\sC$ have arbitrary topology.
Here, however, the distance curves may have several components and there are limit cases where two or more such components intersect each other, albeit ``almost not'': we shall say that closed curves $c_i : \SS^1 \to M$, $i=1,2$, are \emph{quasi disjoint} if there are only finitely many intersection points and if for any $\d > 0$ there exist disjoint curves $c_i^\d : \SS^1 \to M$ satisfying $\dist(c_i^{\d}(t), c_i(t)) < \d$, for all $t \in \SS^1$, $i=1,2$.
One may also rephrase this by saying that the (finitely many) self intersections are \emph{non transversal}. The statement is now

\begin{theorem} \label{t:rectB}
Let $\sM$ be an orientable or non orientable complete unbordered Riemannian surface and $\sC \subset \sM$ a nonempty compact subset. Let $M \subset \sM$ be an open connected component of $\sM \setminus \sC$ and consider, for $d >0$ the sets $Z(d) = \{x \in M \mid \dist(x,\sC) < d\}$. Then, for any open connected component $\O$ of $M \setminus \overline{Z(d)}$ the boundary $\omega = \p\O$ is a finite union of pairwise quasi disjoint simple closed Lipschitz curves.
\end{theorem}

\begin{proof}[Proof of Theorems \ref{t:rectA} and \ref{t:rectB}]
For the case of Theorem \ref{t:rectA} we shall show that for each $p \in \omega$ there exists an open coordinate neighborhood $U_p \subset M$ of $p$ and a parametrized simple Lipschitz curve $w_p : [a_p, b_p] \to M$ with $a_p < b_p$ such that $w_p(a_p) \in \p U_p$,  $w_p(b_p) \in \p U_p$ and
\begin{equation}\label{eq:local}
\omega \cap U_p = \{w_p(r) \mid a_p < r < b_p \}.
\end{equation}
Since $\omega$ is compact this implies that any closed connected component of $\omega$ is a simple closed curve. Furthermore, since $M$ is orientable and has genus zero, $\omega$ can have only one component. Hence, once \eqref{eq:local} is achieved the proof of Theorem \ref{t:rectA} will be complete. For Theorem \ref{t:rectB}, which is proved at the same time, we shall have the additional situation where at $p$ instead of one curve $w_p$ satisfying \eqref{eq:local} there is a pair of such curves that intersect each other exactly in  $p$ and the intersection is non transversal. (The nature of the intersection is depicted in Figure \ref{fig:Rectif088} with the two touching branches $w_p^{\succ}$ and $w_p^{\prec}$).

The hypothesis that $\omega = \p\O$ is the boundary of a \emph{connected component} enters the proof in the paragraph following \eqref{eq:obs}.

Thus, consider $p \in \omega$ and let $\eta$ be a minimal connection from $p$ to $\sC$. Its length is $L(\eta)=d$. We denote by $q$ the endpoint of $\eta$ on $\sC$. As in \eqref{eq:KT} we introduce Fermi coordinates in an open neighborhood $W_{\eta}$ of $\eta$ such that in these coordinates
\begin{equation*}%\label{eq:KT}
W_{\eta} = \{(r,\th) \mid  -\rho < r < \rho,\; -\rho < \th < d + \rho \},
\end{equation*}
where $\rho > 0$ is an appropriate (small) constant satisfying $\rho < d/2$, and the parametrization of $\eta$ by arc length has the coordinate expression $\eta(t) = (0,d-t)$, $t \in [0, d]$, with $\eta(0)=(0,d) = p$ on $\omega$ and $\eta(d) = (0,0)=q$ on $\sC$. Here we use that $M \subset M'$ for some Riemannian surface $M'$ ($=\sM$ in the case of Theorem \ref{t:rectB}) that contains $\sC$ in its interior. We furthermore take $\rho$ smaller than  half the convexity radius of $M$ at $p$ (i.e. any pair $x,y \in B_p(2\rho)$ is joined by a unique minimizing geodesic and this geodesic is contained in $B_p(2\rho)$).

For any $0<\e< \rho$ we denote by $W_{\eta}^{\e}$ the smaller neighborhood
\begin{equation*}%\label{eq:KT}
W_{\eta}^{\e} = \{(r,\th) \in W_{\eta} \mid -\e < r < \e, \ -\e < \th \leq d \} \cup B_p(\e),
\end{equation*}
where $B_p(\e) = \{x \in M \mid \dist(x,p) < \e \}$. As in \eqref{eq:GT} again, the metric tensor in $W_{\eta}$ is of the form
\begin{equation}\label{eq:Geta}
ds^2 = dr^2 + G_{\eta}(r,\th)^2 d\th^2,
\end{equation}
where $G_{\eta}$ satisfies
\begin{equation}\label{eq:boundG}
1- \kappa^2 r^2 \leq G_{\eta}(r,\th) \leq 1 + \kappa^2 r^2
\end{equation}
for some constant $\kappa >0$. Consequently, there exist constants $0 < \e_{\eta} \leq \rho$ and $\l > 0$ (we may take $\l = \frac{1}{d}+\kappa^2d$) such that the following holds:
\begin{flalign}\label{eq:lessd}
\text{\emph{ For all $x = (r,\th) \in W_{\eta}^{\e_{\eta}}$ with $\th < d-\l r^2$ we have $\dist(x, \sC) \leq \dist(x,q) < d$. }}&&
\end{flalign}
For technical reasons we furthermore take $\e_{\eta} > 0$ so small that, in addition,
\begin{equation}\label{eq:epp}
\kappa^2 \e_{\eta}^2 < \frac{1}{100}, \quad \l \e_{\eta}^2 < \frac{1}{4}\e_{\eta}.
\end{equation}
Hence, all $x \in W_{\eta}^{\e_{\eta}}$ that lie below (speaking with Figure \ref{fig:Rectif01}) the curve $c_{\eta} : (-\e_{\eta}, \e_{\eta}) \to W_{\eta}^{\e_{\eta}}$ defined by
\begin{equation*}
c_{\eta}(t) = (t,d-\l t^2)
\end{equation*}
belong to $Z(d)$. Since $c_{\eta}$ is orthogonal to $\eta$ at $p$ this has the following implication: Assume that $\eta_p^{\sL}$ and $\eta_p^{\sR}$ are two minimal geodesics from $p$ to $\sC$ forming some angle $\a_p < \pi$ at $p$ as illustrated in Figure \ref{fig:Rectif02}. Then in some disk $B_p(r_p)$ centered at $p$ the two regions below the curves $c^{\sL} := c_{\eta_p^{\sL}}$ and $c^{\sR} := c_{\eta_p^{\sR}}$ overlap and so
\begin{flalign}\label{eq:Bp}
\begin{minipage}{0.9\textwidth}
\emph{%
all $x \in B_p(r_p) \setminus \{p\}$ that belong to $M \setminus Z(d)$ lie in the curved angular sector of angle $\pi - \a_p$ at $p$ bounded by parts of $c^{\sL}$ and $c^{\sR}$
}
\end{minipage}
\end{flalign}
(shaded area in Figure \ref{fig:Rectif02}). Since in any neighborhood of $p$ there are points of $\O$ we conclude that there exist two extremal minimal connections $\eta_p^{\sL}$ and $\eta_p^{\sR}$ from $p$ to $\sC$ forming an angular sector $A_p$ at $p$ with \emph{extremal angle} $\a_p \leq \pi$ (if $\a_p = \pi$ then $A_p$ is one of the two half disk sectors), where the extremal property is that \emph{any other} minimal geodesic $\eta$ from $p$ to $\sC$ has its initial part $\eta \cap B_p(r_p)$ in $A_p$. (If there would exist three minimal connections, say $\eta_1$, $\eta_2$, $\eta_3$ from $p$ to $\sC$ whose tangent vectors at $p$ do not lie in a shared half plane, then for some disk $B_p(r_p)$, $p$ would be the only point in $B_p(r_p)$ at distance $\geq d$ from $\sC$, contradicting the fact that $p \in  \p\O$).

Thus, speaking with Figure \ref{fig:Rectif02} (which from now on illustrates the extremal case),
\begin{flalign}\label{eq:LRmost}
\begin{minipage}{0.9\textwidth}
\emph{%
$\eta_p^{\sL}$ is the left most and $\eta_p^{\sR}$ the right most shortest geodesic from $p$ to $\sC$.}
\end{minipage}
\end{flalign}
Of course, the generic case is that $\eta_p^{\sL} = \eta_p^{\sR}$, but we do not make use of this fact.

We distinguish two cases for the extremal angle: $\a_p < \pi$ and $\a_p = \pi$.

%\includegraphics[width=0.2\linewidth]{imgs/a_Rectif03.eps}%[width=0.526\linewidth]
%\caption{Strip bounded by minimal connections}
%\label{fig:Rectif03}
%\end{minipage}
%\hspace{-0.15\textwidth}%%%%%
%\end{figure}

\emph{Case a)} Assume that $\a_p < \pi$. In this case (which concludes after \eqref{eq:Up5} with the desired Lipschitz arc through $p$) property \eqref{eq:Bp} is valid, where from now on $\eta_p^{\sL}$ and $\eta_p^{\sR}$ are the extremal minimal connections.  We set
\begin{equation*}
\e_p = \min\{\e_{\eta_p^{\sL}},\e_{\eta_p^{\sR}},\rho/4\}
\end{equation*}
and for any $\e \leq \e_p$,
%%%
\begin{align*}
&\begin{alignedat}{4}
W_p^{\sL,\e} &= \{(r,\th) \in W_{\eta_p^{\sL}}^{\e}%
\mid
-\e < r \leq 0,\ \th \geq d-\rho \},
\quad&
W_p^{\sR,\e} &= \{(r,\th) \in W_{\eta_p^{\sR}}^{\e}%
\mid%
0 \leq r < \e,\ \th \geq d-\rho \}\\
b_p^{\sL,\e} &= \{(r,\th) \in W_p^{\sL,\e}%
\mid
\th = d-\rho \},
\quad&
b_p^{\sR,\e} &= \{(r,\th) \in W_p^{\sR,\e}%
\mid
\th = d-\rho \},
\end{alignedat}\\
&\rule{0pt}{12pt}W_p^{\e} = W_p^{\sL,\e} \cup W_p^{\sR,\e}.
\end{align*}
%%%
%
Reducing the size of $\e_p$, if necessary, we have that
\begin{flalign}\label{eq:outside}
\begin{minipage}{0.9\textwidth}
\emph{%
there are points in $\O$ that lie outside $W_p^{\e_p}$.
}
\end{minipage}
\end{flalign}
Since $\rho$ is smaller than half the convexity radius of $M$ at $p$ we then have that
\begin{flalign}\label{eq:onlyp}
\begin{minipage}{0.9\textwidth}
\emph{%
the parts $\{(r,\th) \in W_p^{\sL,\e_p} \mid d-\rho \leq \th \leq d\}$ and $\{(r,\th) \in W_p^{\sR,\e_p} \mid d-\rho \leq \th \leq d\}$ of $W_p^{\e}$ intersect each other only in $p$.
}
\end{minipage}
\end{flalign}
This is illustrated in Figure \ref{fig:Rectif02}. We now claim that
\begin{flalign}\label{eq:inW}
\begin{minipage}{0.9\textwidth}
\emph{%
there exists $\e_p'$ with $0 < \e_p' \leq \frac{1}{2}\e_p$ such that for all $x \in W_p^{\e_p'} \setminus \{p\}$ any shortest connection $\eta$ from $x$ to $\sC$  has an initial segment $\breve{\eta} \subset \eta$ contained in $W_p^{\e_p}$ that goes from $x$ to $b_p^{\sL,\e_p}$ or $b_p^{\sR,\e_p}$.
}
\end{minipage}
\end{flalign}

%
%%% Figures 14, 15, 16 %%%
%
\begin{figure}[b!]
\hspace{-0.20\textwidth}%%%%
%%First Figure
\begin{minipage}[b]{0.42\textwidth}
\centering
\leavevmode
\SetLabels
(.45*.82) $p$\\
(.45*.15) $q$\\
(.45*.50) $\eta$\\
(.15*.70) $c_{\eta}$\\
(.85*.70) $c_{\eta}$\\
(.20*.12) $\sC$\\
(.923*.227) $r$\\
(.565*.93) $\theta$\\
(.91*.40) $W_{\eta}$\\
\endSetLabels
%\ShowGrid
\AffixLabels{%
\includegraphics{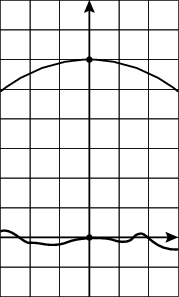}%[width=0.445\linewidth]
}
\caption{A Fermi coordinate neighborhood}
\label{fig:Rectif01}
\end{minipage}
\hspace{-0.1\textwidth}
%%Second Figure
\begin{minipage}[b]{0.42\textwidth}
\centering
\leavevmode
\SetLabels
(.245*.38) $W_p^{\hspace{-0.5pt}\sL\hspace{-1.2pt},\e\hspace{-1.2pt}{}_p}$\\
(.746*.38) $ W_p^{\hspace{-0.5pt}\sR\hspace{-1.2pt},\e\hspace{-1.2pt}{}_p}$\\
(.493*.62) $ \alpha_p$\\%%
(.49*.470) $ B\hspace{-1.1pt}{}_p\hspace{-0.8pt}(\!r \hspace{-1.3pt}{}_p\!)$\\
(.5*.048) $ B\hspace{-0.7pt}{}_p\hspace{-0.8pt}(\rho)$\\
(.385*.27) $ \eta_p^{\sL}$\\
(.595*.27) $ \eta_p^{\sR}$\\
(.32*.740) $c^{\sL}$\\
(.66*.740) $c^{\sR}$\\
(.90*.85) $W_p^{\e_p}$\\
(.10*-.01) $b_p^{\hspace{-0.5pt}\sL\hspace{-1.2pt},\e\hspace{-1.2pt}{}_p}$\\
(.90*-.01) $b_p^{\hspace{-0.5pt}\sR\hspace{-1.2pt},\e\hspace{-1.2pt}{}_p}$\\
(.493*.825)$ p$\\
\endSetLabels
%\ShowGrid
\AffixLabels{%
\includegraphics{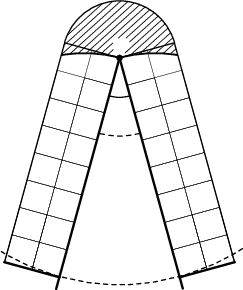}%[width=0.552\linewidth]
}
\caption{Left and right most minimal connections}
\label{fig:Rectif02}
\end{minipage}
\hspace{-0.1\textwidth}
%%Third Figure
\begin{minipage}[b]{0.42\textwidth}
\centering
\leavevmode
\SetLabels
(.10*-.025) $b_p^{\hspace{-0.5pt}\sL\hspace{-1.2pt},\e\hspace{-1.2pt}{}_p}$\\
(.90*-.025) $b_p^{\hspace{-0.5pt}\sR\hspace{-1.2pt},\e\hspace{-1.2pt}{}_p}$\\
(.399*.3) $ \eta_p^{\sL}$\\%
(.600*.3) $\eta_p^{\sR}$\\
(.15*.45) $\hat{\eta}^{\sL}$\\
(.87*.45) $\hat{\eta}^{\sR}$\\
(.202*.195) $S$\\
(.545*.88) $p$\\
(.6*1.0) $\hat{x}$\\
\endSetLabels
%\ShowGrid
\AffixLabels{%
\includegraphics{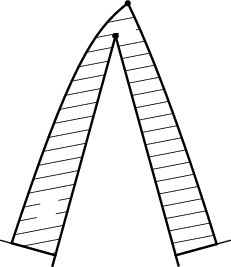}%[width=0.526\linewidth]
}
\caption{Strip bounded by minimal connections}
\label{fig:Rectif03}
\end{minipage}
\hspace{-0.20\textwidth}
\end{figure}
%
%%% end Figures 14, 15, 16 %%%
%

For the proof assume that this is not the case. Then there exists a sequence $\{x_k\}_{k=1}^{\infty}$ of points in $W_p^{\e_p}$, where the distance from $x_k$ to $(\eta_p^{\sL} \cup \eta_p^{\sR}) \cap B_p(\rho)$ goes to zero and where for each $k$ some minimal connection $\eta_k$ from $x_k$ to $\sC$ leaves $W_p^{\e_p}$ at some point $y_k \in \G := \p W_p^{\e_p} \setminus (b_p^{\sL, \e_p} \cup \eta_p^{\sL} \cup \eta_p^{\sR} \cup b_p^{\sR, \e_p})$. There exists then a sequence of indices $\{k_n\}_{n=0}^{\infty}$ such that the $x_{k_n}$ converge to some point $\tilde{x} \in (\eta_p^{\sL} \cup \eta_p^{\sR}) \cap \overline{B_p(\rho)}$, the $y_{k_n}$ to some point $\tilde{y} \in \G$ and the $\eta_{k_n}$ to a minimal connection $\tilde{\eta}$ from $\tilde{x}$ via $\tilde{y}$ to $\sC$. By the minimality of $\tilde{\eta}$, $\eta_p^{\sL}$ and $\eta_p^{\sR}$ the only possible location for $\tilde{x}$ is $\tilde{x}=p$. But, since $\tilde{\eta}$ leads from $p$ to a point on $\G$ this contradicts the fact that $\eta_p^{\sL}$ is leftmost and $\eta_p^{\sR}$ rightmost.

The next step is to compare shortest connecting curves to $\sC$ with the Fermi coordinate lines $\{r = \text{constant}\}$ in $W_p^{\sL,\e_p}$ and $W_p^{\sR,\e_p}$. We claim that for any $\nu > 0$ there exists $0 < \e_p^{\nu} \leq \e_p'$ with the following property, for $\sQ = \sL$ and $\sQ = \sR$ (here the exponents of $\e_p$ are indices, not powers).
\begin{flalign}\label{eq:angleb}
\begin{minipage}{0.9\textwidth}
\emph{%
Let $x \in W_p^{\sQ,\e_p^{\nu}}$, let $\eta$ be a shortest connection from $x$ to $\sC$ and assume that it intersects $b_p^{\sQ,\e_p}$. Then the initial segment $\breve{\eta}$ of $\eta$ from $x$ to $b_p^{\sQ,\e_p}$ is contained in $W_p^{\sQ,\e_p}$
and for any intersection of $\breve{\eta}$ with a coordinate line $\{r = constant\}$ the angle $\zeta$ at the intersection point (the smaller of the two) satisfies $\abs{\zeta} < \nu$.}
\end{minipage}
\end{flalign}
%
%\footnotetext[1]{by \eqref{eq:inW} we already know that $\breve{\eta}$ is contained in $W_p^{\e_p}$, but here we add that it stays in the part $W_p^{\sQ,\e_p}$}
%
For the proof we take $\sQ = \sR$, for $\sQ = \sL$ the arguments are the same. We first show that if $\e_p^{\nu} \leq \e_p'$ then $\breve{\eta}$ stays in $W_p^{\sR,\e_p}$. We already know by \eqref{eq:inW} that it stays in $W_p^{\e_p} = W_p^{\sL,\e_p} \cup W_p^{\sR,\e_p}$. Hence, we only have to show that $\breve{\eta}$ cannot enter  $W_p^{\sL,\e_p} \setminus W_p^{\sR,\e_p}$. But the latter is only possible by crossing the geodesic boundary segment $\{(r,\th) \in \p W_p^{\sR,\e_p} \mid r = 0, \ d \leq \th \leq d + \e_p \}$, and since $\breve{\eta}$ must return to $W_p^{\sR,\e_p}$ it would have to intersect this segment twice which is impossible because $\e_p$ is smaller than the convexity radius of $M$ at $p$.

For the statement about the angles we can now proceed in the same way as in the proof of \eqref{eq:inW}: If it were false,
we could find an infinite sequence of $\eta_{k_n}$ converging to a minimal connection $\tilde{\eta}$ from some point $\tilde{x} \in \eta_p^{\sR} \cap W_p^{\sR,\e_p}$ to $\sC$ that forms a positive angle with $\eta_p^{\sR}$ at $\tilde{x}$. Again, since $\eta_p^{\sR}$ and $\tilde{\eta}$ are both minimal connections to $\sC$ the only possible location for $\tilde{x}$ is $\tilde{x} = p$ and we get a contradiction with the fact that $\eta_p^{\sR}$ is rightmost. \eqref{eq:angleb} is now established.

We are getting closer to the intended parametrization of $\omega$ in a neighborhood of $p$. Let us fix $\nu \leq \frac{1}{100}$ and $\e_p^{\nu}$ as in \eqref{eq:angleb}, where we recall the various epsilons
\begin{equation*}
\e_p^{\nu} \leq \e_p' \leq \frac{1}{2}\e_p = \frac{1}{2} \min\{\e_{\eta_p^{\sL}},\e_{\eta_p^{\sR}},\rho/4\}.
\end{equation*}
We begin with the part in $W_p^{\sR,\e_p}$. Fix $r \in [0,  \frac{1}{2}\e_p^{\nu}]$. On the partial parameter line $h_r : [d-\rho, d+\e_p^{\nu}] \to W_p^{\sR,\e_p}$ defined by
\begin{equation*}
h_r(\th) = (r, \th), \quad d-\rho \leq \th \leq d+ \e_p^{\nu},
\end{equation*}
we seek an intersection with $\omega$. By \eqref{eq:epp} and because $\e_p^{\nu} < \e_{\eta_p^{\sR}}$ we have $\l r^2 < \l (\e_p^{\nu})^2 < \frac{1}{4} \e_p^{\nu}$. Hence, the point $h_r(d-\frac{1}{4}\e_p^{\nu})$ lies below the curve
$t \mapsto c^{\sR}(t) = c_{\eta_p^{\sR}}(t) = (t,d - \l t^2)$, $t \in [0, \e_p^{\nu}]$,
and so, by the remark subsequent to \eqref{eq:epp} \emph{the distance from $h_r(d-\frac{1}{4}\e_p^{\nu})$ to $\sC$ is less than $d$}.

For the  $\th$-s in the interval $[d - \rho, d+ \e_p^{\nu}]$ there are two cases.

\emph{Case 1)} \ There exists $\acute{\th} \in [d - \rho, d+ \e_p^{\nu}]$ such that some minimal geodesic $\acute{\eta}$ from $\acute{x} = h_r(\acute{\th})$ to $\sC$ intersects $b_p^{\sL,\e_p}$.  Then $\acute{\th} \geq d - 2\e_p \geq d -\frac{1}{2}\rho $, because otherwise $\dist(\acute{x}, \sC) \leq r + d-2\e_p < d - \e_p$ while by \eqref{eq:inW}, $\acute{\eta}$ intersects the coordinate line $\{\th = d \}$ in $W_p^{\sR,\e_p}$ at some point $x'=(r',d)$ with $\dist(x',\sC) > d - r' > d-\e_p$.

Let thus $\hat{\th} \in [d - \frac{1}{2}\rho, d+\e_p^{\nu}]$ be the infimum of all such $\acute{\th}$. Then we find a decreasing convergent sequence $\acute{\th}_k \to \hat{\th}$ such that the corresponding $\acute{\eta}_k$ converge to some minimal geodesic $\hat{\eta}^{\sL}$ from $\hat{x} = h_r(\hat{\th})$ to $\sC$ that intersects $b_p^{\sL,\e_p}$. At the same time there exists a monotone increasing convergent sequence $\grave{\th}_k \to \hat{\th}$ with corresponding minimal geodesics $\grave{\eta}_k$ converging to a minimal geodesic $\hat{\eta}^{\sR}$ from $\hat{x}$ to $\sC$ that intersects $b_p^{\sR,\e_p}$. The situation is depicted in Figure \ref{fig:Rectif03}. Since both geodesics are minimal we have $L(\hat{\eta}^{\sL}) = L(\hat{\eta}^{\sR})$. The observation is now that
\begin{equation}\label{eq:obs}
L(\hat{\eta}^{\sL}) = L(\hat{\eta}^{\sR}) > d.
\end{equation}
Assume, for a contradiction, that the common length is $\leq d$. Then $\hat{\eta}^{\sL} \cup \hat{\eta}^{\sR}$ is contained in $\overline{Z(d)}$, in the same way as $\eta_p^{\sL} \cup \eta_p^{\sR}$, $b_p^{\sL,\e_p}$ and $b_p^{\sR,\e_p}$. Hence, the boundary of the simply connected  strip $S \subset W_p^{\e_p}$ bounded by
$(\hat{\eta}^{\sL} \cup \hat{\eta}^{\sR}) \cap W_p^{\e_p}$, $(\eta_p^{\sL} \cup \eta_p^{\sR}) \cap W_p^{\e_p}$and parts of $b_p^{\sL,\e_p}$ and $b_p^{\sR,\e_p}$ (shaded domain in Figure \ref{fig:Rectif03}) does not intersect $\O$. Since $p \in \p\O$ and $\O$ is connected it follows that $\O \subset S \subset W_p^{\e_p}$ in contradiction to \eqref{eq:outside}.

Since by \eqref{eq:obs} the distance from $\hat{x} = h_r(\hat{\th})$ to $\sC$ is greater than $d$ it follows that  $\hat{x}$ lies above the coordinate line $\{ \th = d - \frac{1}{4}\e_p^{\nu}\}$, i.e. we have $\hat{\th} > d - \frac{1}{4}\e_p^{\nu}$. At the same time, because we are in Case 1), we have $\hat{\th} \leq d +  \e_p^{\nu}$.
By the continuity of the distance function it follows that there exists a real number $f_p(r) \in (d-\frac{1}{4}\e_p^{\nu}, d+ \e_p^{\nu})$ such that
\begin{equation}\label{eq:fpr}
\dist(h_r(f_p(r)), \sC) = d
\end{equation}
and we let $f_p(r)$ be the smallest value in $(d-\frac{1}{4}\e_p^{\nu}, d+ \e_p^{\nu})$ with this property. That $h_r(f_p(r))$ is on the boundary of $\O$ and not of some other component of $M \setminus \overline{Z(d)}$ will be shown later (see \eqref{eq:Up}).

\emph{Case 2)} \ For all $\th \in [d-\rho, d+ \e_p^{\nu}]$ any minimal geodesic from $h_r(\th)$ to $\sC$ intersects $b_p^{\sR,\e_p}$. Here we use a monotonicity argument.

Take  $x=h_r(\th_x)$, $y=h_r(\th_y)$ with $d-\frac{1}{4}\e_p^{\nu} \leq \th_x < \th_y \leq d+ \e_p^{\nu}$ and let $\eta_x$, $\eta_y$ be minimal geodesics from $x$ and $y$ to $\sC$ Figure \ref{fig:Rectif04}). By \eqref{eq:angleb} $\eta_y$ has an initial part in $W_p^{\sR,\e_p}$ running from $y$ to $b_p^{\sR,\e_p}$. It therefore intersects the coordinate line $\{\th = \text{constant} = \th_x \}$ in some point $z$. By \eqref{eq:angleb} the angle $\zeta$ at $y$ between $\eta_y$ and $h_r$ satisfies $\abs{\zeta} < \nu \leq \frac{1}{100}$. Using that $\dist(x,\sC) \leq \dist(x,z) + \dist(z,\sC)$, that $\dist(z,\sC) + \dist(y,z) = \dist(y, \sC)$, and that by \eqref{eq:boundG} and \eqref{eq:epp} the metric tensor \eqref{eq:Geta} is close to the Euclidean metric $ds_{\EE}^2 = dr^2 + d\th^2$ we get%
\begin{equation}\label{eq:dmono}
\dist(y,\sC) - \dist(x,\sC) \geq \dist(y,z)-\dist(x,z) > \frac{9}{10}(\th_y - \th_x).
\end{equation}
This shows that along $h_r$ the distance to $\sC$ is monotone increasing. By the triangle inequality $\dist(h_r(d),\sC) \geq \dist(p,\sC) -\dist(p,h_r(d)) = d-r \geq d - \frac{1}{2} \e_p^{\nu}$. Hence, by \eqref{eq:dmono} $\dist(h_r(d+\e_p^{\nu}),\sC) > d$. Therefore, as in Case 1), there exists a minimal (in this case, in fact, unique) value $f_p(r) \in (d-\frac{1}{4}\e_p^{\nu}, d+ \e_p^{\nu})$ for which \eqref{eq:fpr} holds.

%
%%% Figures 17, 18, 19 %%%
%
\begin{figure}[b!]
\hspace{-0.15\textwidth}
%%First Figure
\begin{minipage}[b]{0.42\textwidth}
\centering
\leavevmode
\SetLabels
(-.05*.37) $p$\\
\R(.113*.93) $\theta$\\
\R(-.005*.78) $\theta_y$\\
\R(-.005*.51) $\theta_x$\\
\R(-.005*.21) $\eta_p^{\sR}$\\
(.307*.21) $\eta_x$\\
(.901*.30) $\eta_y$\\
(.65*.48) $x$\\
(.65*.82) $y$\\
(.86*.48) $z$\\
(.35*.54) $r$\\
(.35*.82) $r$\\
(.665*.05) $h_r$\\
(.645*.68) $\zeta$\\
\endSetLabels
%\ShowGrid
\AffixLabels{%
\includegraphics{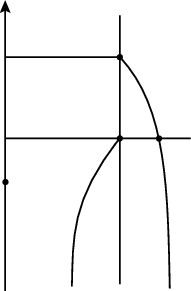}
}
\caption{Small angles with coordinate lines}
\label{fig:Rectif04}
\end{minipage}
\hspace{-0.08\textwidth}
%%Second Figure
\begin{minipage}[b]{0.42\textwidth}
\centering
\leavevmode
\SetLabels
\R(.168*.93) $\theta$\\
\R(.06*.42) $f_p(r_1)$\\
\R(.06*.74) $f_p(r_2)$\\
(.20*.308) $w_p$\\
(.52*.02) $h_{r_1}$\\
(.82*.02) $h_{r_2}$\\
(.395*.455) $x_1$\\
(.693*.775) $x_2$\\
(.71*.456) $z$\\
(.30*.18) $d$\\
(.885*.18) $d$\\
\endSetLabels
%\ShowGrid
\AffixLabels{%
\includegraphics{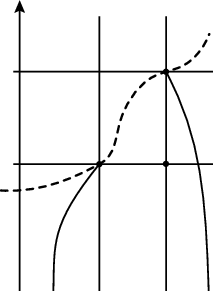}
}
\caption{Proof of the Lipschitz continuity}
\label{fig:Rectif05}
\end{minipage}
\hspace{-0.08\textwidth}
%%Third Figure
\begin{minipage}[b]{0.42\textwidth}
\centering
\leavevmode
\SetLabels
(.09*.93) $\theta$\\
\R(-.02*.074) $d{-}\rho$\\
\R(-.02*.656) $d$\\
(.10*.30) $\eta_p^{\sR}$\\
(.30*.105) $b_p^{\hspace{-0.5pt}\sL\hspace{-1.2pt},\e\hspace{-1.2pt}{}_p}$\\
(.10*.30) $\eta_p^{\sR}$\\
(.63*.01) $h_r$\\
(.83*.01) $A$\\
(.835*.30) $\eta_x$\\
(.71*.41) $a$\\
(.40*.61) $b$\\
(.41*.41) $c$\\
(.818*.10) $z$\\
(.715*.71) $y$\\
(.635*.825) $x$\\
(.625*.605) $\gamma$\\
\R(-.02*.825) $f_p(\!r\!)$\\
(.076*.625) $p$\\
\endSetLabels
%\ShowGrid
\AffixLabels{%
\includegraphics{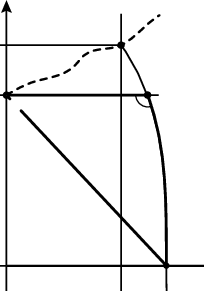}
}
\caption{For the right-hand differentiability at $p$}
\label{fig:Rectif066}
\end{minipage}
\hspace{-0.15\textwidth}
\end{figure}
%
%%% end Figures 17, 18, 19 %%%
%

Bringing Cases 1) and 2) together we have obtained a function $f_p :  [0, \frac{1}{2}\e_p^{\nu}] \to (d-\frac{1}{4}\e_p^{\nu},d+\e_p^{\nu})$ (which a few lines hereafter will be shown to be Lipschitz) and with it a curve  $w_p : [0, \frac{1}{2}\e_p^{\nu}] \to \p Z(d) \cap W_p^{\sR,\e_p}$ defined
by
\begin{equation}\label{eq:curvwr}
w_p(r) = (r,f_p(r)), \quad 0\leq r \leq \tfrac{1}{2} \e_p^{\nu}.
\end{equation}
Note that by the definition of $f_p$ we have the following properties:
\begin{flalign}\label{eq:inWp}
\begin{minipage}{0.9\textwidth}
\emph{%
$w_p(0) = p$. Any point in the closure of $W_p^{\sR,\frac{1}{2}\e_p^{\nu}}$ that lies below $w_p$ is closer to $\sC$ than $d$. Furthermore, for all $r \in (0,\tfrac{1}{2} \e_p^{\nu}]$ any minimal geodesics $\eta$ from $w_p(r)$ to $\sC$ has an initial part $\breve{\eta}  \subset W_p^{\sR,\e_p}$ going from $w_p(r)$ to $b_p^{\sR,\e_p}$.
}
\end{minipage}
\end{flalign}
As a consequence of \eqref{eq:inWp}, inequality \eqref{eq:dmono} also holds in Case 1) as long as $x$ and $y$ are below or on the curve $w_p$.

We now prove that $f_p$ is a Lipschitz function. Let $0 \leq r_1 < r_2 \leq \frac{1}{2}\e_p^{\nu}$ and consider $x_1 = w_p(r_1) = (r_1, f_p(r_1))$, $x_2 = w_p(r_2) = (r_2, f_p(r_2))$. We introduce the auxiliary point $z=(r_2, f_p(r_1))$ on the coordinate line $h_{r_2}$ through $x_2$ (Figure \ref{fig:Rectif05}). By the triangle inequality,
\begin{equation*}
\abs{\dist(z,\sC)-d} = \abs{\dist(z,\sC)-\dist(x_1,\sC)} \leq r_2 - r_1,
\end{equation*}
and by \eqref{eq:dmono} (which is applicable, as mentioned after \eqref{eq:inWp}) we have
\begin{equation*}
\abs{\dist(z,\sC)-d} = \abs{\dist(z,\sC)-\dist(x_2,\sC)} \geq \frac{9}{10}\abs{f_p(r_2)-f_p(r_1)}.
\end{equation*}
With the two inequalities together,
\begin{equation}\label{eq:fLip}
\abs{f_p(r_2) - f_p(r_1)} \leq 1.2 \abs{r_2 - r_1}.
\end{equation}
Hence, $f_p$ and with it the curve $w_p : [0, \frac{1}{2}\e_p] \to \p Z(d) \cap W_p^{\sR,\e_p} $ are Lipschitz with Lipschitz constant $1.2$.
For $w_p$ the Lipschitz constant has a meaning only with respect to the chosen Fermi coordiante system.)
In particular, $w_p$ satisfies Definition \ref{d:Lipcurve}.

In the part in $W_p^{\sL,\e_p}$ we proceed in the same way, where we reduce the sizes of $\e_p$, $\e_p'$ and $\e_p^{\nu}$, if necessary, so that all the preceding statements hold in analogous form and with the same constants when we replace $\sR$ by $\sL$. We thus get a simple Lipschitz curve
\begin{equation}\label{eq:wr}
w_p : [-\tfrac{1}{2}\e_p^{\nu}, \tfrac{1}{2}\e_p^{\nu}] \to \p Z(d).
\end{equation}
That the combination of the two simple parts corresponding to $r \in [-\frac{1}{2}\e_p^{\nu}, 0]$, respectively $r \in [0,\frac{1}{2}\e_p^{\nu}]$ is itself simple follows from the fact that by \eqref{eq:inWp}, for $r < 0$ the minimal geodesics from $w_p(r)$ to $\sC$ intersect $b_p^{\sL,\e_p}$ but not $b_p^{\sR,\e_p}$, while for $r > 0$ the minimal geodesics from $w_p(r)$ to $\sC$ intersect $b_p^{\sR,\e_p}$ but not $b_p^{\sL,\e_p}$. Hence, the only common point of the two parts is $p$.

We now show that for some open neighborhood $U_p$ of $p$ we have $w_p(-\frac{1}{4}\e_p^{\nu}) \in \p U_p$,  $w_p(\frac{1}{4}\e_p^{\nu}) \in \p U_p$ and
\begin{equation}\label{eq:Up}
 \omega \cap U_p = \{w_p(r) \mid -\tfrac{1}{4}\e_p^{\nu} < r < \tfrac{1}{4}\e_p^{\nu} \},
\end{equation}
which then is \eqref{eq:local} for $w_p$ restricted to the interval $[a_p,b_p] = [-\frac{1}{4}\e_p^{\nu},\frac{1}{4}\e_p^{\nu}]$.

To define $U_p$ we set $w_p^{\sL}=\{w_p(r) \mid -\frac{1}{4}\e_p^{\nu} \leq r \leq 0\}$, $p^{\sL} =w_p(-\frac{1}{4}\e_p^{\nu})$, $w_p^{\sR}=\{w_p(r) \mid 0\leq r \leq \frac{1}{4}\e_p^{\nu}\}$, $p^{\sR} =w_p(\frac{1}{4}\e_p^{\nu})$ and let
\begin{equation*}
\d = \frac{1}{4}\min\{\dist(p^{\sL}, w_p^{\sR}), \dist(p^{\sR}, w_p^{\sL}), r_p, \e_p^{\nu}\},
\end{equation*}
where $r_p$ is the radius of the disk $B_p(r_p)$ in \eqref{eq:Bp} (see Figure \ref{fig:Rectif02}). With this we set
\begin{equation*}
\setlength{\arraycolsep}{1.5pt}
\begin{array}{ccccc}
U_p^{\sL} &=& \{ (r,\th) \in W_p^{\sL,\e_p} \mid & -\tfrac{1}{4}\e_p^{\nu} < r < 0, &\abs{f_p(r)-\th} < \d\},
\\
U_p^{\sR} &=& \{ (r,\th) \in W_p^{\sR,\e_p} \mid &  0 < r < \tfrac{1}{4}\e_p^{\nu} , &\abs{f_p(r)-\th} < \d\},
\\
A_p &=& B_p(\d) \setminus \text{int} W_p^{\e_p},\vphantom{W_p^{\sL,\e_p}}
\end{array}
\end{equation*}
where $\text{int}$ denotes the interior of a subset of $M$, and now define (see Figure \ref{fig:Rectif077}, where a case with $\a_p > \frac{1}{2} \pi$ is shown),
\begin{equation*}
U_p= U_p^{\sL} \cup U_p^{\sR} \cup A_p.
\end{equation*}
 By construction we have
\begin{equation}\label{eq:Up1}
w_p(-\tfrac{1}{4}\e_p^{\nu}) \in \p U_p, \quad w_p(\tfrac{1}{4}\e_p^{\nu}) \in \p U_p, \quad w_p(r) \in U_p \text{ for all } r \in (-\tfrac{1}{4}\e_p^{\nu},\tfrac{1}{4}\e_p^{\nu}).
\end{equation}
The path
\begin{equation*}
\mathfrak{w}_p := \{w_p(r) \mid -\tfrac{1}{4}\e_p^{\nu} < r < \tfrac{1}{4}\e_p^{\nu} \}
\end{equation*}
(thickened in Figure \ref{fig:Rectif077}) separates $U_p$ into two open subsets $\invbreve{U}_p$ and $\breve{U}_p$ defined as follows,
\begin{equation*}
\setlength{\arraycolsep}{1.5pt}
\begin{array}{ccl}
\invbreve{U}_p &=& \{(r,\th) \in U_p^{\sL} \mid \th < f_p(r) \} \cup \{(r,\th) \in U_p^{\sR} \mid \th < f_p(r) \} \cup (A_p \setminus \{p\}),
\\[3pt]
\breve{U}_p &=& U_p \setminus (\breve{U}_p \cup \mathfrak{w}_p).
\end{array}
\end{equation*}
(The part $\breve{U}_p$ is shaded Figure \ref{fig:Rectif077}.) By \eqref{eq:Bp}, since $\d < r_p$, we have $A_p \setminus \{p\} \subset Z(d)$. Together with \eqref{eq:inWp} and the analogue of \eqref{eq:inWp} for $W_p^{\sL,\e_p}$ we get
\begin{equation}\label{eq:Up2}
\invbreve{U}_p \subset Z(d).
\end{equation}
For the upper part we claim that
\begin{equation}\label{eq:Up3}
\breve{U}_p  \subset M \setminus \overline{Z(d)}.
\end{equation}
For the proof we consider $x \in \breve{U}_p$ and let $\eta$ be a shortest connection from $x$ to $\sC$. We distinguish four cases. 1°) If $x = (r_x,\th_x) \in W_p^{\sR,\e_p}$ and $\eta$ intersects $b_p^{\sR,\e_p}$, then $\th_x > f_p(r_x)$ and by \eqref{eq:angleb} $\eta$ very closely follows the coordinate line $\{r = \text{constant} = r_x\}$ and therefore intersects the trace $\{w_p(r) \mid 0\leq r<\frac{1}{2}\e_p^{\nu}\}$, which implies that $\dist(x,\sC) = L(\eta) > d$. 2°) $x \in W_p^{\sL,\e_p}$ and $\eta$ intersects $b_p^{\sL,\e_p}$. This is the mirror case of 1° with the same conclusion. 3°) Assume that $x  \in W_p^{\sR,\e_p}$ and $x \notin W_p^{\sL,\e_p}$ but $\eta$ intersects $b_p^{\sL,\e_p}$ (when the angle $\a_p$ is big as in Figure \ref{fig:Rectif077} this case cannot occur). Then, since by \eqref{eq:inW} the initial segment $\breve{\eta}$ of $\eta$ from $x$ to $b_p^{\sL,\e_p}$ is contained in $W_p^{\e_p}$, $\eta$ enters $W_p^{\sL,\e_p}$ at some point $\acute{y} = (0,\acute{\th})$ (coordinate system of $W_p^{\sL,\e_p}$) and then crosses the coordinate line $\{(r,d) \in W_p^{\sL,\e_p} \mid -\e_p < r \leq 0 \}$ at some point $\acute{z}$. Since $\eta$ is a minimal connection we have $\dist(x,\sC) = \dist(x,\acute{y}) + \dist(\acute{y},\acute{z}) + \dist(\acute{z},\sC)$; furthermore, we have $d = \dist(p,\sC) \leq \dist(p,\acute{z}) + \dist(\acute{z},\sC)$ and $\dist(\acute{y},\acute{z}) \geq \dist(p,\acute{z})$ (by the properties of the Fermi coordinates). Bringing this together we get
\begin{equation*}
\dist(x,\sC) \geq \dist(x,\acute{y}) + d > d.
\end{equation*}
4°) $x  \in W_p^{\sL,\e_p}$ and $x \notin W_p^{\sR,\e_p}$ but $\eta$ intersects $b_p^{\sR,\e_p}$. This is the mirror case of 3° with the same conclusion. There are no other cases and the proof of \eqref{eq:Up3} is complete.

Combining \eqref{eq:Up2} and \eqref{eq:Up3} we get, given that $\mathfrak{w}_p$ lies on $\p Z(d)$
\begin{equation}\label{eq:Up4}
\p Z(d) \cap U_p =  \mathfrak{w}_p.
\end{equation}
Since $p \in \omega = \p\O$ there are points of $\O$ arbitrarily close to $p$. Since $\invbreve{U}_p \subset Z(d)$ these points lie in $\breve{U}_p$. Thus, $\breve{U}_p \cap \O \neq \emptyset$. Since $\breve{U}_p$ does not intersect $\p Z(d)$ it does not intersect $\p\O$. Since, furthermore, $\breve{U}_p$ is connected and contains points of $\O$ it follows that $\breve{U}_p \subset \O$. Now every $x \in \mathfrak{w}_p$ is adjacent to $\breve{U}_p$ and thus adjacent to $\O$ without, however, being member of $\O$. Therefore we have
\begin{equation}\label{eq:Up5}
\mathfrak{w}_p \subset \omega .
\end{equation}
By \eqref{eq:Up1}, \eqref{eq:Up4} and \eqref{eq:Up5} and since $\omega \subset \p Z(d)$ we now have
\begin{equation*}
\omega \cap U_p \subset \p Z(d) \cap U_p = \mathfrak{w}_p = \mathfrak{w}_p \cap U_p \subset \omega \cap U_p.
\end{equation*}
This proves \eqref{eq:Up} and accomplishes the proof of \eqref{eq:local}  in Case a).

\emph{Case b)} Assume now that $\eta_p^{\sL}$ and $\eta_p^{\sR}$ form an angle $\a_p = \pi$. In this special case the regions ``below'' the curves $c^{\sL} = c_{\eta_p^{\sL}}$ and $c^{\sR} = c_{\eta_p^{\sR}}$ do not overlap and \eqref{eq:Bp}
has to be modified.

%
%%% Figures 20, 21 %%%
%
\begin{figure}[t!]
\hspace{-0.2\textwidth}
\begin{minipage}[b]{0.45\textwidth}
\centering
\leavevmode
\SetLabels
(.255*.385) $ p$\\
(.205*.305) $ A_p$\\
(.36*.10) $ \eta_p^{\sR}$\\
(.16*.59) $ \eta_p^{\sL}$\\
(0.95*.34) $ w_p^{\sR}$\\
(0.72*.90) $ w_p^{\sL}$\\
%(0.82*.90) $\scriptstyle w_p^{\sL}$\\%
(0.72*.76) $ p^{\sL}$\\
(0.86*.44) $ p^{\sR}$\\
(0.5*.355) $\scriptstyle \invbreve{U}_p$\\
(0.457*.596) $\scriptstyle \invbreve{U}_p$\\
(0.57*.435) $\scriptstyle \breve{U}_p$\\
\endSetLabels
%\ShowGrid
\AffixLabels{%
\includegraphics{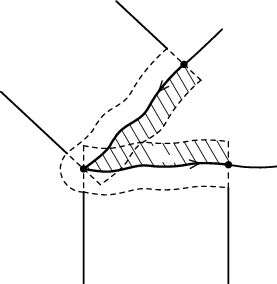}
}
\caption{~Distance curve crossing a neighborhood}
\label{fig:Rectif077}
\end{minipage}
\hspace{-0.0\textwidth}
%%Second picture
\begin{minipage}[b]{0.57\textwidth}
\centering
\leavevmode
\SetLabels
(.33*.62) $A_p^{\succ}(r)$\\
(.67*.62) $A_p^{\prec}(r)$\\
(.03*.72) $w_p^{\succ}$\\
(.03*.22) $w_p^{\succ}$\\
(.99*.22) $w_p^{\prec}$\\
(.99*.72) $w_p^{\prec}$\\
(.52*.44) $p$\\
\L(.12*.225) $c^{\sR}$\\
\L(.14*.70) $c^{\sL}$\\
\R(.84*.225) $c^{\sR}$\\
\R(.86*.70) $c^{\sL}$\\
(.455*.88) $\eta_p^{\sL}$\\
(.53*.08) $\eta_p^{\sR}$\\
(.14*.47) $\Omega$\\
(.84*.47) $\Omega$\\
\endSetLabels
%\ShowGrid
\AffixLabels{%
\includegraphics{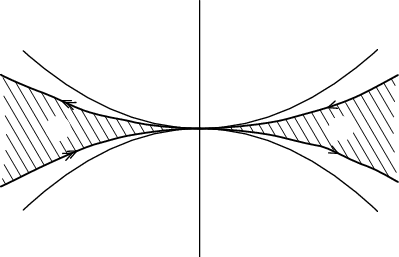}
}
\caption{~Cuspidal double point of $\p\O$ at boundary point $p$}
\label{fig:Rectif088}
\end{minipage}
\hspace{-0.2\textwidth}
\end{figure}
%
%%% end Figures 20, 21 %%%
%

We let again $\e_p = \min\{\e_{\eta_p^{\sL}},\e_{\eta_p^{\sR}},\rho/4\}$, but now replace $B_p(r_p)$  by $B_p(\e_p)$.
The geodesic $\eta_p^{\sL} \cup \eta_p^{\sR}$ separates $B_p(\e_p)$ into two half disks and \emph{in the case of Theorem \ref{t:rectA} it also disconnects $M$} (by the Jordan curve theorem). Instead of
\eqref{eq:Bp} we now have
\begin{flalign}\label{eq:App}
\begin{minipage}{0.9\textwidth}
\emph{%
all $x \in B_p(\e_p)$ that belong to $\O$ lie in the union of the two curved angular sectors $\sA_p^{\succ}$, $\sA_p^{\prec}$ of angle $0$ at $p$ in the two half disks bounded by parts of $c^{\sL}$ and $c^{\sR}$ (see Figure \ref{fig:Rectif088}).
}
\end{minipage}
\end{flalign}
We now distinguish two cases.

\emph{Case b1)} $p \in \p (\O \cap A_p^{\prec})$ and $p \notin \p (\O \cap A_p^{\succ})$. Taking $\e_p$ sufficiently small we may assume for this case that $\O \cap A_p^{\succ}  = \emptyset$. This situation necessarily occurs in the case of Theorem 1 because  in that case $\O$ lies entirely on one side of $\eta_p^{\sL} \cup \eta_p^{\sR}$ in $M$, given that $\O$ is connected and $\eta_p^{\sL} \cup \eta_p^{\sR}$ separates $M$. Reducing the size of $\e_p$ further, if necessary, we may assume that $\O$ contains points outside $B_p(\e_p)$.

The arguments are now the same as in Case a), though with some simplifications and up to some obvious adaptions: We extend $W_p^{\sR,\e}$ to $\{(r,\th) \mid 0 \leq r < \e, \ d-\rho \leq \th \leq d + \rho \}$, let $h_r$ be the coordinate lines $h_r(\th) = (r,\th)$, $\th \in [d-\rho, d + \rho]$, and translate the monotonicity based on \eqref{eq:dmono} into the statement that along $h_r$ between $c^{\sR}$ and $c^{\sL}$ the distance to the lower part of $\sC$ in $W_p^{\sR,\e_p}$ is continuously increasing while the distance to the upper part of $\sC$ in $W_p^{\sR,\e_p}$ is continuously decreasing. At the point $\hat{x} =(r, \hat{\th})$ of equality these distances, in the same way as in \eqref{eq:obs}, must be greater than $d$ for otherwise the coordinate line $h_r$ would separate $\O$. Hence, there are exactly two values $\th_r^{\sR}=: f_p^{\sR}(r) < \th_r^{\sL} =: f_p^{\sL}(r)$ for which these distances are equal to $d$ for the first, respectively for the last time. Similarly to \eqref{eq:fLip} the functions $f_p^{\sR}$, $f_p^{\sL} : [0,\frac{1}{2}\e_p^{\nu}] \to \RR$ are 1.2-Lipschitz continuous and the definition of $w_p$ can now be given as follows, where we modify the notation:
\begin{equation}\label{eq:wpPi}
w_p^{\prec}(r) = (-r,f_p^{\sL}(-r)), \, r \in [-\tfrac{1}{2}\e_p^{\nu},0], \quad w_p^{\prec}(r) = (r,f_p^{\sR}(r)), \, r \in [0,\tfrac{1}{2}\e_p^{\nu}].
\end{equation}
For $U_p$ one may now proceed as earlier, but it is not difficult to see that in the present situation the requested property \eqref{eq:local} is also achieved with the following
\begin{equation*}
U_p =  \{(r,\th) \mid -\tfrac{1}{4}\e_p^{\nu} < r < \tfrac{1}{4}\e_p^{\nu}, \ d-\e_p < \th < d+\e_p\}.
\end{equation*}
The proof of Theorem \ref{t:rectA} is now complete, and for Theorem \ref{t:rectB}, Case b1) is achieved.

\emph{Case b2)} $p \in \p (\O \cap A_p^{\succ})$ and $p \in \p (\O \cap A_p^{\prec})$. In this case we have, in addition to \eqref{eq:wpPi}, a curve
\begin{equation}\label{eq:wpPii}
w_p^{\succ}(r) = (r,g_p^{\sL}(r)), \, r \in [-\tfrac{1}{2}\e_p^{\nu},0], \quad w_p^{\succ}(r) = (-r,g_p^{\sR}(-r)), \, r \in [0,\tfrac{1}{2}\e_p^{\nu}].
\end{equation}
with analogous Lipschitz functions $g_p^{\sL}, g_p^{\sR} : [-\frac{1}{2}\e_p^{\nu},0] \to \RR$. The two curves meet in $p$ and the intersection can be removed with arbitrarily small homotopies in a neighborhood of $p$.

Since $\omega$ is compact Case b2) can occur only at finitely many points along $\omega$. Furthermore, by its compactness $\omega$ can be covered with finitely many neighborhoods $U_p$ for which $U_p \cap \omega$ is either an arc $w_p$ satisfying  \eqref{eq:local} or a pair of arcs $w_p^{\succ}$, $w_p^{\prec}$ as in \eqref{eq:wpPi}, \eqref{eq:wpPii}.
For each pair $w_p^{\succ}$, $w_p^{\prec}$ the intersection of $w_p^{\succ}$ and $w_p^{\prec}$ at $p$ is non transversal. By symbolically doubling $p$ in each pair, they become pairs of disjoint simple arcs and $\omega$ becomes a union of finitely many pairwise disjoint simple closed Lipschitz curves. Reversing the doublings again, these curves become quasi disjoint. The proof of Theorem \ref{t:rectB} is now complete.
\end{proof}

We add a remark concerning differentiability. For the estimate of the Lipschitz constant of $f_p$ in \eqref{eq:fLip} (see Figure \ref{fig:Rectif05}) we made use of the triangle inequality. Using an argument based on Toponogov's theorem instead, the estimate of the constant in \eqref{eq:fLip} can be significantly improved. We shall detail this for the growth of $f_p$ at $p$. As a result we shall get the following theorem.
In its statement we use the symbol $\frac{1}{\#}$ for the normalisation of a vector
\begin{equation}\label{eq:normalize}
\frac{1}{\#}X := \frac{1}{||X||} X.
\end{equation}
Although the use of the symbol $\frac{1}{\#}$ makes the denominators in \eqref{eq:RLtang} superfluous, we keep them for later discussion.
\begin{theorem}\label{t:rectC} Let $M$, $\sC$ and the Lipschitz curve $\omega : \SS^1 \to M$ be as in Theorems \ref{t:rectA} and \ref{t:rectB}. Let $t_p \in \SS^1$ and $p = \omega(t_p)$. Then there exist the left-hand and right-hand tangent vectors
\begin{equation}\label{eq:RLtang}
\omega'^{\sL}(t_p) =  \lim_{t\to t_p^-} \frac{1}{\#}\frac{\tilde{\omega}(t)-\tilde{\omega}(t_p)}{t-t_p},
\quad
\omega'^{\sR}(t_p) = \lim_{t\to t_p^+} \frac{1}{\#}\frac{\tilde{\omega}(t)-\tilde{\omega}(t_p)}{t-t_p},
\end{equation}
where $\tilde{\omega}$ is the lift of $\omega$ in a small neighborhood of $p$ to the tangent plane of $M$ at $p$ via the inverse of the exponential map.

Furthermore, let $\eta_p^{\sL}$ and $\eta_p^{\sR}$ be the left and right most minimal geodesics from $p$ to $\sC$ as in \eqref{eq:LRmost}. Then
\begin{equation}\label{eq:RLorth}
\omega'^{\sL}(t_p) \perp \eta_p^{\sL},
\quad
\omega'^{\sR}(t_p) \perp \eta_p^{\sR}.
\end{equation}
In particular, if $\eta_p^{\sL} = \eta_p^{\sR}$, then $\omega$ is differentiable at $t_p$.
\end{theorem}
\begin{proof}
We pick up at the point in the proof of Theorems \ref{t:rectA} and \ref{t:rectB} where the Lipschitz property of $f_p$ is investigated, but now replace the configuration of Figure \ref{fig:Rectif05} by the one in Figure \ref{fig:Rectif066}, where, in order to adapt the notation, we assume that $t_p = 0$ (so that $f_p(t_p) = f_p(0) = d$) and write $r$ instead of $t$:

Let $0 < r < \frac{1}{2}\e_p^{\nu}$ and $x=(r,f_p(r))$.
First, we are going to prove that
\begin{equation*}
\limsup_{r \to 0^+}\frac{f_p(r) - f_p(0)}{r} \leq 0.
\end{equation*}
For this we shall establish \eqref{eq:fprr} below. If $f_p(r) \leq d$, then \eqref{eq:fprr} holds trivially. Let us therefore assume that $f_p(r) > d$.

Let $\eta_x$ be a minimal connection from $x$ to $\sC$. Since $\eta_x$ runs below $w_p$ it intersects the coordinate line $\{ \th = \text{constant} = d\}$ in some point $y = (b,d)$ and the coordinate line $b_p^{\sR,\e_p} = \{ \th = \text{constant} = d - \rho\}$ in some point $z$ (see Figure \ref{fig:Rectif066}). We recall that $\rho$ is smaller than half the convexity radius of $M$ at $p$. Drawing the minimal geodesic $c$ from $p$ to $z$ we get a convex geodesic triangle with vertices $p$, $z$, $y$ and sides $a$ from $z$ to $y$, $b$ from $y$ to $p$, and $c$. We denote by $\g$ the angle opposite $c$. The points $y$ and $z$ decompose $\eta_x$ into three segments: $u$ from $x$ to $y$, then $a$ from $y$ to $z$ and $A$ from $z$ to $\sC$. For the lengths we have
\begin{equation*}
L(\eta_x) = d = u + a + A,
\end{equation*}
where we again abbreviate $L(u) = u$, etc. Since $d = \dist(p,\sC) \leq b+a+A$, we have the rough estimate $u \leq b$. Since by \eqref{eq:epp} the metric tensor \eqref{eq:Geta} is close to the $ds_{\EE}^2 = dr^2 + d\th^2$ and since $\eta_x$ forms angles $\leq \nu \leq \frac{1}{100}$ with the coordinate lines $\{r = \text{constant}\}$ we have $b \leq r + 2 \nu u \leq r + 2 \nu b$, from which follows that
\begin{equation*}%\label{eq:bnu}
b \leq \frac{r}{1-2\nu}  < r (1 + 3 \nu).
\end{equation*}
By Lemma \ref{l:trigineq} below, and since $\rho \leq c < 2\rho$ and $\g \leq \frac{\pi}{2} + \nu$, there exists a constant $\sK$ depending only on $\rho$ (which is fixed) and on the lower bound of the Gauss curvature on $B_p(\rho)$ such that $c-a \leq \sK (b \sin(\nu) + b^2)$.

Since $b < r (1 + 3 \nu)$ we get, being generous with factors,
\begin{equation}\label{eq:cma}
c-a \leq 2\sK (r \sin(\nu) + r^2).
\end{equation}
Now $u+a+A = d =\dist(p,\sC) \leq c + A$, and so $u \leq c - a$. Since $u$ is almost the same as $f_p(r) - d = f_p(r) - f_p(0)$, it follows from \eqref{eq:cma} that
\begin{equation}\label{eq:fprr}
f_p(r) - f_p(0) \leq 3 \sK (r \sin(\nu) + r^2).
\end{equation}
Since this holds for any $r \in (0, \frac{1}{2}\e_p^{\nu})$ and for arbitrarily small $\nu >0$ we get
\begin{equation*}
\limsup_{r \to 0^+}\frac{f_p(r) - f_p(0)}{r} \leq 0.
\end{equation*}
On the other hand, by \eqref{eq:lessd} we have
\begin{equation*}
\liminf_{r \to 0^+}\frac{f_p(r) - f_p(0)}{r} \geq 0.
\end{equation*}
Hence, the right-hand derivative of $f_p$ at $t_p$ exists and is equal to 0.
It follows that the right hand tangent vector
\begin{equation*}
w_p'^{\sR}\vphantom{w}(0) = \lim_{r \to 0^+} \frac{1}{r}(w_p(r) -w_p(0))
\end{equation*}
of $w_p$ at $p$ exists and is orthogonal to the coordinate line $r=0$, i.e.\ orthogonal to $\eta_p^{\sR}$. Going from the Fermi coordinates to geodesic normal coordinates in a neighborhood of $p$ via the inverse exponential map $\exp_p^{-1}$ we get the following, where as usual the tangent plane $T_pM$ of $M$ at $p$ and its own tangent plane $T_0T_pM$ at the origin are identified so that for tangent vectors at $p$ the differential map $(\exp_p^{-1})^*$ operates as the identity. We also remove the simplification $t_p =0$.
\begin{equation}\label{eq:wpr}
w_p'^{\sR}(0) = (\exp_p^{-1})^*(w_p'^{\sR}(0))
=
\lim_{r \to 0^+} \frac{\exp_p^{-1}\!w_p(r) - \exp_p^{-1}\!w_p(0)}{r} = \lim_{t \to t_p^+} \frac{\tilde{\omega}(t) - \tilde{\omega}(t_p)}{t-t_p}.
%\lim_{r \to 0^+} \frac{1}{r}(\exp_p^{-1}\!w_p(r) - \exp_p^{-1}\!w_p(0)) = \lim_{t \to t_p^+} \frac{\tilde{\omega}(t) - \tilde{\omega(t_p)}}{t-t_p}.
\end{equation}
%$
In the last limit we have a special parametrization of $\tilde{\omega}$ at $p$ coming from $w_p$. Normalizing the quotient with $\frac{1}{\#}$ we get the expression for the tangent vector in \eqref{eq:RLtang} which is independent of the choice of the parametrization. \eqref{eq:RLtang} and \eqref{eq:RLorth} are now proved for $\omega'^{\sR}(t_p)$. For $\omega'^{\sL}(t_p)$ the proof is the same.
\end{proof}

We append some additional properties which follow rather directly from relation \eqref{eq:angleb} in the proof of Theorems \ref{t:rectA}, \ref{t:rectB}.

\begin{corollary} \label{c:omegadiff}\label{c:rectC} Let $M$, $\sC$ and the Lipschitz curve $\omega : \SS^1 \to M$ be as in Theorems \ref{t:rectA},  \ref{t:rectB}, \ref{t:rectC}. Then:
\begin{enumerate}
\item For all $t_p \in \SS^1$ we have the following continuity of the one-sided tangent vectors, where $\Pmove{t}{t_p}$ denotes the Riemannian parallel transport along $\omega$ of tangent vectors at $\omega(t)$ to the tangent plane at $\omega(t_p)$.
\begin{equation*}
\arraycolsep=1.8pt
\begin{array}{ccccc}
\omega'^{\sR}(t_p) &=& \displaystyle \lim_{t \to t_p^+}\Pmove{t}{t_p}(\omega'^{\sR}(t)) &=& \displaystyle\lim_{t \to t_p^+}\Pmove{t}{t_p}(\omega'^{\sL}(t)),
\\
\omega'^{\sL}(t_p) &=& \displaystyle\lim_{t \to t_p^-}\Pmove{t}{t_p}(\omega'^{\sR}(t)) &=& \displaystyle\lim_{t \to t_p^-}\Pmove{t}{t_p}(\omega'^{\sL}(t)).
\end{array}
\end{equation*}
\item There exists a countable subset $\sN \subset \SS^1$ such that
\begin{equation*}
\omega'^{\sL}(t) = \omega'^{\sR}(t) \quad \text{for all $t \in \SS^1 \setminus \sN$}.
\end{equation*}
\item
If $\omega$ is parametrized in proportion to arc length then \eqref{eq:RLtang} can be replaced by
\begin{equation*}
\omega'^{\sL}(t_p) =  \frac{1}{\#}\lim_{t\to t_p^-}\frac{\tilde{\omega}(t)-\tilde{\omega}(t_p)}{t-t_p},
\quad
\omega'^{\sR}(t_p) = \frac{1}{\#}\lim_{t\to t_p^+}\frac{\tilde{\omega}(t)-\tilde{\omega}(t_p)}{t-t_p}.
\end{equation*}
\end{enumerate}
\end{corollary}

\begin{proof}
We resume the Fermi coordinates used in the proof of Theorem \ref{t:rectC}. For $\omega(t)$ within the coordinate neighborhood we denote by $\mathbf{e}^1(t)$,
$\mathbf{e}^2(t)$ the unit tangent vectors at $\omega(t)$ parallel, respectively to the coordinate lines $\th = \text{constant}$ and $r = \text{constant}$. We use the symbol $\sphericalangle$ to denote the absolute angular measure between tangent vectors.

Let now $\nu>0$ be arbitrarily small and consider $\e_p^{\nu}$ as in \eqref{eq:angleb}.
We further reduce the size of $\e_p^{\nu}$, if necessary, so that for $\abs{t - t_p} < \e_p^{\nu}$ we have $\sphericalangle(\Pmove{t}{t_p}(\mathbf{e}^i(t)), \mathbf{e}^i(t_p))< \nu$, $i=1,2$ (that this is possible can be deduced from the Gauss-Bonnet theorem). For $\abs{t - t_p} < \frac{1}{2}\e_p^{\nu}$ \eqref{eq:angleb} then says that the left and right most minimal connections $\eta_{\omega(t)}^{\sL}$, $\eta_{\omega(t)}^{\sL}$ form angles smaller than $\nu$ with $-\mathbf{e}^2(t)$ at $\omega(t)$. By \eqref{eq:RLorth}, with $t_p$ replaced by $t$,  we have therefore
%%
%\begin{equation}\label{eq:omegadiff1}
%\begin{alignedat}{3}
%&\sphericalangle(\omega'^{\sL}(t),\mathbf{e}^1(t)) < \nu,
%\sphericalangle(\omega'^{\sR}(t),\mathbf{e}^1(t)) < \nu,
%\\
%&\sphericalangle(\omega'^{\sL}(t),\omega'^{\sR}(t)) < 2\nu,
%&& \quad
%\text{for all} \ s \in (t- \tfrac{1}{2}\e_p^{\nu},t) \cup (t,t+ \tfrac{1}{2}\e_p^{\nu}).
%\end{alignedat}
%\end{equation}
%%
%
\begin{equation}\label{eq:omegadiff1}
\sphericalangle(\omega'^{\sL}(t),\mathbf{e}^1(t)) < \nu, \quad
\sphericalangle(\omega'^{\sR}(t),\mathbf{e}^1(t)) < \nu,
 \quad
\text{for} \ s \in (t- \tfrac{1}{2}\e_p^{\nu},t) \cup (t,t+ \tfrac{1}{2}\e_p^{\nu}).
\end{equation}
By the bound for the angle of rotation under parallel transport we get
\begin{equation}\label{eq:omegadiff2}
\sphericalangle(\omega'^{\sL}(t_p),\Pmove{t}{t_p}\omega'^{\sL}(t)) < 2\nu, \ \text{for} \ t \in (t_p- \tfrac{1}{2}\e_p^{\nu},t_p),\quad
\sphericalangle(\omega'^{\sR}(t_p),\Pmove{t}{t_p}\omega'^{\sR}(t)) < 2\nu, \ \text{for} \ t \in (t_p,t_p+ \tfrac{1}{2}\e_p^{\nu}).
\end{equation}
Property (1) is now a direct consequence of \eqref{eq:omegadiff1} and \eqref{eq:omegadiff2} because all vectors have length 1. For (2) we note that by compactness, $\SS^1$ is finitely covered by intervals $(t_p- \tfrac{1}{2}\e_p^{\nu},t_p+\tfrac{1}{2}\e_p^{\nu})$ and that by \eqref{eq:omegadiff1}, in each of these intervals we have $\sphericalangle(\omega'^{\sL}(t), \omega'^{\sR}(t)) < 2\nu$, except possibly for $t = t_p$. Hence, for any $\nu > 0$, only finitely many couples $\omega'^{\sL}(s)$, $\omega'^{\sR}(s)$ form an angle greater than $2\nu$, and thus the relation $\omega'^{\sL}(s) \neq \omega'^{\sR}(s)$ is possible only for countably many $s$. For (3) we restrict ourselves to a sketch: in the coordinate neighborhood $W_p^{\sR,\e}$ the small angles in \eqref{eq:omegadiff2} imply that the parametrization of $\omega$ by $w_p$ in \eqref{eq:wpr} satisfies
\begin{equation*}
\lim_{t \to t_p^+} \frac{L(\omega |_{[t_p,t]})}{t-t_p}
=
\lim_{r \to 0^+} \frac{L(w_p |_{[0,r]})}{r} = 1.
\end{equation*}
Together with \eqref{eq:wpr} this implies the right hand side in (3) and the left hand side is obtained in the same way.
\end{proof}

\begin{lemma}\label{l:trigineq}
Let $T$ be a Riemannian triangle with sides $a$, $b$, $c$ and angle $\g$ opposite $c$, with $a \le a_0$, $b \le b_0$, $c \ge c_0$.
Assume that the Gauss curvature on $T$ satisfies $K \geq - k^2$, $k>0$, and that $\g \leq \frac{\pi}{2} + \zeta$ for some $\zeta \in [0, \frac{\pi}{2})$.
Then, there exists a constant $C$ which just depends on $a_0$, $b_0$, $c_0$ and $k$, such that
\begin{equation*}
c-a \leq C(b \sin(\zeta) + b^2 ).
\end{equation*}
\end{lemma}

\begin{proof}
Since, up to the value of the constant $C$, the inequality remains invariant if the metric is scaled by a factor we may assume that $k=1$. By Toponogov's comparison theorem it then suffices to prove the lemma for the case where $T$ is a triangle in the hyperbolic plane. In this latter case we have
\begin{equation*}
\begin{aligned}
\cosh{c} & = \cosh{a} \cosh{b} - \sinh{a} \sinh{b} \cos(\g) \leq \cosh{a} \cosh{b} + \sinh{a} \sinh{b} \sin(\zeta),
\\
\cosh{c} - \cosh{a} & \leq \cosh{a} (\cosh{b} -1) + \sinh{a} \sinh{b} \sin(\zeta).
\end{aligned}
\end{equation*}
Since $(\cosh t -1)/t^2$ and $(\sinh t)/t$ are increasing functions on $(0,\infty)$, we obtain
\begin{equation*}
\begin{aligned}
\cosh{a} (\cosh{b} -1)
& = \cosh{a} \, \frac{\cosh{b} -1}{b^2} \, b^2
\le \cosh{a_0} \, \frac{\cosh{b_0} -1}{b_0^2} \, b^2 ,
\\
\sinh{a} \sinh{b} \sin(\zeta)
& = \sinh{a} \, \frac{\sinh{b}}{b} \, b \sin(\zeta)
\le \sinh{a_0} \, \frac{\sinh{b_0}}{b_0} \, b \sin(\zeta) ,
\\
\cosh{c} - \cosh{a}
& = \, \frac{\cosh{c} - \cosh{a}}{c-a} \, (c-a)
= \, \frac{2 \sinh \frac{c+a}{2} \sinh \frac{c-a}{2}}{c-a} \, (c-a)
\\
& \ge (c-a) \sinh \frac{c+a}{2}
\ge (c-a) \sinh \frac{c_0}{2} \,,
\end{aligned}
\end{equation*}
and the desired inequality follows.%
\end{proof}

\begin{upremark}\label{r:biblio}
We mention a few variants of Theorems \ref{t:rectA} and \ref{t:rectB}. Further bibliography may be found in \cite{RZ}.

The earliest result is by Fiala \cite{Fi}, who uses distance sets in the study of the isoperimetric problem. The setting is that of a real analytic surface $M$ homeomorphic to $\RR^2$ endowed with a complete real analytic metric and $\sC : \SS^1 \to M$ is a simple closed real analytic unit speed curve. Under these conditions it is shown (in 9.2, page 325) that on either side of $\sC$ for any $d > 0$ the distance set $S(d) = \p Z(d) \setminus \sC$, if not empty, consists of finitely many real analytic arcs. The arguments make use of the Gaussian parallel curve $\sC_d(t) = \exp(d \cdot \dot{\sC}(t)^{\perp})$, where $\dot{\sC}(t)^{\perp}$ is the normal vector field along $\sC$ pointing towards the chosen side of $\sC$. The result implies Theorem \ref{t:rectA} in the analytic setting.

Hartman, in \cite{Ha}, replaces the analyticity condition by class $C^2$ differentiability.  Under this hypothesis he shows (in Proposition 6.1, page 717) that for almost all $d$ the distance set $S(d)$ is a union of finitely many $C^2$- differentiable arcs. The approach is by an analysis of the cut locus of $\sC$ on the chosen side of $\sC$ and its intersections with the Gaussian parallel curves $\sC_d$. The focus of \cite{Ha} is to prove Fiala's isoperimetric results (in stronger form) in the $C^2$ setting.

For the same setting as in \cite{Ha} Shiohama proves in \cite[Theorem B, page 127]{Sh} that if, in addition, $M$ has total curvature, then for some $d_0>0$ all $S(d)$ with $d>d_0$ are homeomorphic to a circle.

That $S(d)$ is piecewise smooth for almost all $d$ was later proved by Shiohama and Tanaka in \cite{ShT1} for complete surfaces $M$ of any topology and $\sC$ a smooth simple closed curve. The approach in \cite{ShT1} is a generalization of the analysis of the cut locus in \cite{Ha}. In \cite{ShT2}, Shiohama and Tanaka analyse the cut locus in an entirely different way that uses only metric arguments. As a result they show (in Theorem B, page 334) that if $M$ is a complete Alexandrov surface without boundary whose curvature is bounded from below and if $\sC \subset M$ is a compact subset, then for almost all $d$ (in the range of distances to $\sC$) the distance set $S(d)$ is a finite union of pairwise disjoint simple closed rectifiable curves.

In \cite{RZ} Rataj and Zaj{\'i}{\v c}ek follow another line (with its history given in their introduction) using that for a Riemannian manifold the distance $d(x, \sC)$ to a compact subset $\sC$ is locally semi convex, where the concept of local semi convexity is inherited from that in $\RR^n$ via the coordinate maps and is itself not related to the Riemannian metric. Then, from this property alone they show (Theorem 5.3, page 458) that for the case where $M$ is a complete unbordered two dimensional Riemannian manifold and $\sC \subset M$ a compact subset, there exists a set $\sE \subset \RR$ of exceptional values that has Hausdorff measure $\sH^{1/2}(\sE) = 0$ such that for any distance $d \notin \sE$ the distance set $S(d)$ is a one dimensional Lipschitz manifold. Their concept of Lipschitz manifold (Definition 2.3, page 447) differs slightly from the one we use.

In \cite{BMO} Blokh, Misiurewicz and Oversteegen prove (Theorem 1.1, page 734) that if $M$ is a  possibly bordered topological 2-manifold endowed with a proper geodesic metric, and if $\sC \subset M$ is a compact subset, then, except for countably many $d$ each connected component of $S(d)$ is either a point, a simple closed curve, or a simple arc with its endpoints on the boundary of $M$. The proof is by showing that $S(d)$ is a Suslin set for any $d$.

\end{upremark}

\end{document}